\documentclass[leqno,12pt,a4paper,twoside]{article}
\usepackage{xypic}
\usepackage{babel}
\usepackage{amssymb}
\usepackage{amsmath}
\def\N{{\mathbb N}} 
\def\P{{\mathbb P}} 
 
\def\Z{{\mathbb Z}} 
 
\def\R{{\mathbb R}} 
\def\C{{\mathbb C}}

\begin{document}
\baselineskip=0.52cm

\thispagestyle{empty}

\pagestyle{myheadings}

\markboth{{\small REMOVABLE SINGULARITIES FOR 
INTEGRABLE CR FUNCTIONS}}{{\rm J. {\sc Merker} and E. {\sc Porten}}}

                         $\:$
                        \bigskip
                         \begin{center}
                         ON REMOVABLE SINGULARITIES \\
                         FOR INTEGRABLE CR FUNCTIONS 
                         \end{center}

                         \bigskip
                         \begin{center}
                         {\bf J. Merker} and {\bf E. Porten}
                         \end{center}

\bigskip
\bigskip

The subject of removable singularities for the boundary values of holomorphic
functions of several variables has been intensively studied in recent years, 
see works: \cite{AC},  \cite{HAPO}, \cite{JO1}, \cite{JO2}, \cite{JOMZ},
\cite{JO4},
\cite{KY}, \cite{KYREA}, \cite{LUP}, \cite{LUST}, \cite{ME2}, \cite{PO1}, 
\cite{ROST}, and the excellent surveys \cite{CS}, \cite{ST}.

The present paper is devoted
to a further step into this direction. Its main intention is
to illustrate the way how, on a CR-manifold M of arbitrary codimension, 
removable singularity theorems may be understood in terms
of the analysis of CR-orbits of M.

Let $M$ be a locally embeddable
CR manifold, ${\rm dim}_{CR} M=m$, 
${\rm codim}_{\C^n} \ M=n$, ${\rm dim}_{\R} \ M=d=2m+n$.
Locally, $M\subset \C^{m+n}$.
Let $\Phi\subset M$ be a closed set. 

We give various conditions in order that $\Phi$ is $L^1$-removable, {\it i.e.}
                      \begin{equation}
                      L^1_{loc}(M) \cap L_{loc,CR}^1(M\backslash\Phi
                      )=L_{loc,CR}^1(M)
                      \end{equation} 
Let us recall here that $\Phi\subset \subset M$ is called
$(L^{\rm p},\overline{\partial}_{\rm b})$-removable if $L_{loc,CR}^{\rm
p}(M\backslash \Phi) \cap L_{loc}^{\rm p}(M)= L_{loc,CR}^{\rm
p}(M)$, {\it i.e.} any $f\in L_{loc}^{\rm p}(M)$
satisfying $\int_{M\backslash \Phi} f\overline{\partial}
\varphi=0$ for each $(m+n,m-1)$-form $\varphi\in {\cal
C}^{\infty}_c(\C^{m+n}\backslash \Phi)$ satisfies $\int_{M}
f \overline{\partial}\psi=0$ for each $(m+n,m-1)$-form $\psi \in {\cal
C}^{\infty}_c(\C^{m+n})$. 

Let $H^{\kappa}$ denote
$\kappa$-dimensional Hausdorff measure. Our main results are the following.
Since $L^{\rm p}_{loc}$ embeds in $L^1_{loc}$ for
${\rm p}\geq 1$, we state them for ${\rm p}= 1$.
We give explanations below.

\smallskip   {\sc Theorem 1.} (\cite{JO3}, \cite{PO1}) 
{\it If $M$ is ${\cal C}^3$,
a function $f\in L_{loc}^1(M)$ is CR if and only if $f|_{{\cal O}}$
belongs to $L_{loc,CR}^1({\cal O})$ for almost every CR orbit ${\cal
O}$.}

\smallskip   {\sc Corollary 1.} {\it If $\Phi=\cup_{a\in A}
{\cal O}_a$ is of zero $d$-dimensional measure, then $(1)$ is
satisfied.}

\smallskip Theorem 1 reduces the problem to the case where $M$ is a
single CR orbit, {\it i.e.} $M$ is {\it globally minimal} \cite{ME1}.
All the results also hold for ${\rm p}\geq 1$ equally.

\smallskip   {\sc Theorem 2.}   {\it Let $M$ be ${\cal
C}^{2,\alpha}$, $0< \alpha <1$, ${\rm dim}_{CR} M =m\geq 1$. Every
closed subset $E$ of $M$ such that $M$ and $M\backslash E$ are globally
minimal and such that $H_{loc}^{d-3}(E)<\infty$ is $L^1$-removable.}

\smallskip The notion of wedge- (${\cal W}$-) removability is defined
here in higher codimension (see below).

\smallskip   {\sc Theorem 3.} (\cite{MP2})   {\it Let $M$ be ${\cal
C}^{\omega}$, ${\rm dim}_{CR}M=m\geq 1$. Every closed set $E\subset M$
such that $M$ and $M\backslash E$ are globally minimal and such that
$H^{d-2}(E)=0$ is ${\cal W}$- and $L^1$-removable.}

\smallskip   {\sc Theorem 4.}   {\it Let $M$ be ${\cal
C}^{2,\alpha}$, $0< \alpha <1$, 
$m\geq 1$, and let $N$ be a connected ${\cal C}^2$
submanifold of $M$ such that $M$ and $M \backslash N$ are globally minimal.

{\rm (i)} \ If ${\rm codim}_M N \geq 3$, then $N$ is ${\cal
W}$- and $L^1$-removable{\rm ;}

{\rm (ii)} \ Every closed set $\Phi \subset N$ is ${\cal
W}$- and $L^1$-removable if $\Phi \neq N$, 
${\rm codim}_M N =2$ and $m\geq 1${\rm ;}

{\rm (iii)} \ $N$ is ${\cal W}$- and $L^1$-removable if $N$ is generic at
one point, ${\rm codim}_M N =2$ and $m\geq 2$.}

\smallskip 
{\it Remark.} One of the main feature of our results is that
they are stated in analogy with known results in the hypersurface
case ($n=1$) with variations on this theme. Also, they become
classical and easy in case $n=0$, {\it i.e.} $M$ 
is an open set $\Omega$ in $\C^m$, $m\geq 2$. Finally, 
there is obstruction to removability
in case $N$ is not generic at any point, similar
to $N$ being a complex hypersurface in $\Omega$  if $n=0$.

\smallskip 
{\it Remark.} After the reduction to $M$ being globally minimal,
the assumption that $M\backslash \Phi$ is globally 
minimal too is essential and cannot be dropped (see below).

\smallskip 
{\it Remark.} Theorem 3 is treated in \cite{MP2} and
the following will appear in \cite{PO2}: {\it

{\rm (iv)} \ If $m\geq 2$, ${\rm codim}_M N =1$ and $N$ is generic,
then every closed set $K\subset N$  which does not contain
any CR orbit of $N$ is $L^1$-removable.}

\smallskip Finally, using the theory of CR orbits, we
extend also \cite{KYREA}.  
A set $S\subset M$ is called a ${\cal C}^{\lambda}$ peak
set, $0<\lambda <1$, if there exists a {\it nonconstant} function $\varpi
\in {\cal C}^{\lambda}_{CR}(M)$ such that $S=\{\varpi=1\}$ and
$|\varpi| \leq 1$.

\smallskip   {\sc Theorem 5.}   {\it Let $M$ be ${\cal
C}^{2,\alpha}$ globally minimal. Then every ${\cal C}^{\lambda}$ peak
set $S$ satisfies $H^{d}(S)=0$ and is $L^1$-removable.}

\smallskip   {\sc Corollary 2.}   {\it Let $M$ be ${\cal
C}^3$. Then a ${\cal C}^{\lambda}$ peak set $S$ is $L^1$-removable if
$H^d(\cup_{{\cal O}\subset S} {\cal O})=0$.}

\smallskip
Now, we explain the terminology, compare our results to the
codimension one case and give some motivations.

The general feature of our work is that $L^{\rm p}$-removablity
is linked with ${\cal W}$-removability, {\it i.e.}
with envelope of holomorphy results. 

\smallskip
{\it a. About the methods.} 
Of course, different approaches for proving $(L^{\rm
p},\overline{\partial}_b)$-removability are conceivable,
and the distinguished role of holomorphic hulls in our context 
is not clear yet.

For instance, one could consider the problem from the viewpoint of the
general theory of removable singularities for solutions of linear
partial differential operators. Let $\Omega\subset \R^n$ be a domain,
let $K\subset\subset \Omega$ be a compact and let $P=\sum_{|\beta|\leq
e} a_{\beta}(x) \partial_x^{\beta}$ be such an operator. $K$ is called
$(L^{\rm p}, P)$-removable if any $u\in L^{\rm p}(\Omega)$ with
$Pu\equiv 0$ on $\Omega \backslash K$ satisfies $Pu\equiv 0$ on
$\Omega$ (in the distributional sense).  Of course, the notion makes
sense by replacing $L^{\rm p}(\Omega)$ with other differentiability
classes, {\it e.g.}  ${\cal C}^0(\Omega)$, ${\cal C}^k(\Omega)$, or
even ${\cal D}'(\Omega)$.

Harvey and Polking \cite{HAPO} have proved removable singularity
theorems for general $P$, but they give results only in case
${\rm p}>1$. Indeed, their main theorem 4.1 \cite{HAPO}, states that $K$ is
$(L^{\rm p},\overline{\partial}_{\rm b})$-removable if the Hausdorff measure
$H^{n-p'}(K)<\infty$, $p'=p/(p-1)$, $e={\rm deg} \ P$.
The authors further point out that this result
cannot be improved in terms of Hausdorff measures
in the class of all first order differential operators.
Especially information about $L^1$-removability is
never available on this level (more precisely,
in our special setting results for $p<e/(e-1)$ cannot be
derived from the above theorem).

In fact, one of the main argument here ({\it cf.} \cite{HAPO}, \cite{PIN},
\cite{KYREA}, \cite{JO3}, \cite{JO4}) is to write $$\int_{\Omega}
f{\:}{\!}^{\tau}tP\varphi= \int_{\Omega\backslash K}
f{\:}{\!}^{\tau}P((1-\chi_{\varepsilon})\varphi)+ \int_{K^{\varepsilon}}
f{\:}{\!}^{\tau}P(\chi_{\varepsilon}\varphi)= \int_{K^{\varepsilon}}
f{\:}{\!}^{\tau}P(\varphi \chi_{\varepsilon}),$$ where
$\phi\in {\cal C}^{\infty}_c(\R^n)$, $\varepsilon>0$,
$\chi_{\varepsilon}\in {\cal C}^{\infty}_c(\R^n)$,
$|\nabla\chi_{\varepsilon}|\leq C/\varepsilon$,
$\chi_{\varepsilon}\equiv 1$ on $K$,
$\chi_{\varepsilon}\equiv 0$ on $\Omega \backslash K^{\varepsilon}$,
$K^{\varepsilon}= \{p\in \R^n\ \! {\bf :} \ \! {\rm dist} \ (p,K)<\varepsilon\}$,
${\:}{\!}^{\tau}P$ is the transpose of $P$, and to find conditions on $K$,
$P$, $f$ in order that $\int_{K^{\varepsilon}}
f{\:}{\!}^{\tau}P(\varphi \chi_{\varepsilon})$, 
when $\varepsilon \to 0$. For instance, one can
prove along these lines that, 
given $S\subset M$ a ${\cal C}^1$ submanifold with ${\rm
codim}_M S=1$, then ${\cal C}_{CR}^0(M\backslash S) \cap {\cal
C}^0(M)={\cal C}_{CR}^0(M)$ \cite{PIN} and also that $L_{CR}^{\rm
p}(M\backslash S) \cap L_{loc}^{\rm p}(M)= L_{loc,CR}^{\rm p}(M)$
under the additional condition that $S$ is a complex hypersurface in
$M$ \cite{KYREA}.

Actually, these kinds of arguments do not even answer the following
very simple question, which was the model for our work on Theorem 2:
Can one always remove for $L_{CR}^1$ functions a {\it single}
point of a {\it minimal} hypersurface?

\smallskip   {\sc Example 1.} Let $M\subset \C^{m+1}$, $K=\{p\}
\subset M$. First, it is clear that if $p$ belongs to a complex
hypersurface $S\subset M$, then $p$ is $L^1$-removable (\cite{KYREA} again). 
If $M=b\Omega$, $\Omega \subset \subset \C^{m+1}$, there is a trick
(\cite{AC}, p.115).  Let $\varphi$ be a smooth $(m+1,m-1)$-form with
compact support.  Then the form $\psi:=\varphi-\varphi(p)$ satisfies
that ${\rm supp} \ \psi \cap b\Omega$ is still compact, {\it since
$\overline{\Omega}$ is compact}, and $\overline{\partial} \varphi
=\overline{\partial}\psi$, so that, in proving $\int_{b\Omega}
f\overline{\partial} \varphi=0$, {\it one can assume that
$\varphi(p)=0$}. But then $\int_{\{p\}^{\varepsilon}} f
\overline{\partial} (\chi_{\varepsilon} \varphi)=
\int_{\{p\}^{\varepsilon}} f \varphi \overline{\partial}
\chi^{\varepsilon}+ \int_{\{p\}^{\varepsilon}} f \chi^{\varepsilon}
\overline{\partial} \varphi$ tends to zero, because $\varphi(p)=0$,
$|\nabla\chi_{\varepsilon}|\leq C/\varepsilon$ and
$||f||_{L^1(\{p\}^{\varepsilon})}\to 0$.
{\it This argument fails in case $M$ is a local piece of a minimal
hypersurface.} $\square$

\smallskip
{\it b. Obstructions to removability.} Let us consider
the archetypal non-removable singularity: Assume that
a generic CR manifold $M$ intersects transversely
a smooth analytic hypersurface $X$ given as the 
vanishing-locus of a holomorphic function $f$.
Then it is easily seen that the intersection
$M\cap X$ is a smooth submanifold of $M$ of 
codimension $2$ and that $f|(M\backslash (M\cap X))$ is locally
integrable but not CR. We note that $M\cap X$
is itself a CR manifold of CR dimension
$\hbox{CRdim}M-2$ and of course nowhere generic: obstruction
to 4 (iii).

Theorem 4 is a partial generalization to higher dimensions of a theorem
of J\"oricke \cite{JO4} which states that {\sl every
connected, 2-codimensional, smooth subvariety
of a boundary $M\subset\C^n,n\geq 3$ is
removable if it is not maximally complex
in the sense of Harvey and Lawson}.
This paper gave evidence of the fact 
that, excepted for boundaries in $\C^2$, the removal
of singularities in CR-manifolds of CRdimension $1$
should be especially hard.

We explain a second obstruction, du to orbits, after Theorem 6.

\smallskip
{\it c. The strategy.}
Here, a good strategy would be to understand {\it a priori} the {\it a
posteriori} fact that, if $\{p\}$ was $L^1$-removable, if $M$ is
minimal at $p$, then $L_{CR}^1(M\backslash \{p \})$ would extend
holomorphically to one side of $M$ at $p$. In other words, we are
looking for a holomorphic extension to a wedge attached to a full
neighborhood of $p$ in the hope to re-obtain the original function 
around $p$ as a $L_{CR}^1$ boundary value of the extension.

This strategy has been endeavoured by J\"oricke \cite{JO4}
for hypersurfaces or in CR dimension 
$m\geq 2$ by the second author \cite{PO1}. 

Our main objective in this
article is to push forward in any codimension $n\geq 1$
and in any CR dimension $m\geq 1$ the theory of 
$L^{\rm p}$-removability for 
the induced $\overline{\partial}$ in the context of
the {\it extension theory} of Tr\'epreau and Tumanov. 
This theory has acquired recently its final level of
maturity, but there is still no complete survey about its global 
aspects.
%
% La derniere phrase me semble un peu vague. Peut-etre on peut
% l'inserrer dans le paragraphe suivant.

\smallskip

{\it d. CR orbits.} So, let us first recall the fundamental
picture of a CR manifold $M$, its decomposition in the so-called CR
orbits (\cite{SU}, \cite{TR2}, \cite{TU5}, \cite{JO3}, \cite{JO4},
\cite{ME1}).  An immersed connected submanifold $S\subset M$ is called
a {\it CR-integral} submanifold if $T_p S\supset T_p^c M$, for all $p\in S$. 
A fundamental result of Sussmann, applied to our context for the first
time by Treves, shows that $M$ carries a finest partition into
CR-submanifolds - the so called {\it CR-orbits} - which are complete with
respect to their manifold-topology. The description of CR
orbits as subsets of $M$ is very simple. A point $q$ is in the CR
orbit of $p\in M$, {\it i.e} 
$q\in {\cal O}_{CR}(M,p)$, if and only if there
exists a piecewise smooth integral curve of $T^cM$ with origin $p$ and
target $q$. The content of Sussmann's Theorem 4.1 \cite{SU} is that
{\it these sets have a structure of smooth immersed submanifolds of $M$}.

\smallskip
{\sc Analogy.} The CR orbits are kinds of {\it irreducible
components} of CR manifolds, like 
irreducible components of analytic sets, but
in general, of transfinite cardinal.

\smallskip

We emphasize that the behaviour of CR-functions is dominated by
the CR-orbits of the supporting manifold. First, Theorem 1
asserts that the obstructions for the construction of integrable
CR-functions may be expressed in terms of CR-orbits and splits up
into the differential aspect, which is completely localized on
the orbits, and some global part depending on the overall geometry
of the decomposition into orbits.

On a much more advanced level, CR-orbits appear as propagators of 
stuctural properties of CR-objects, like vanishing in a neighborhood 
of a point, or being CR extendable to a manifold with boundary.
For our purposes the theory of analytic extensions into wedges,
as developed by Tr\'epreau and Tumanov (\cite{TR2}, \cite{TU3}) 
is of paramount importance.
Their work culminated in the following theorem, conjectured by Tr\'epreau 
and proved independently by B.~J\"oricke and the first author
(\cite{ME1}, \cite{JO3}):
For each CR orbit ${\cal O}_{CR}$, there
exists an analytic wedge ${\cal W}^{an}$ attached to ${\cal O}_{CR}$,
{\it i.e.} a conic complex manifold with edge ${\cal O}_{CR}$ and with
${\rm dim}_{\C}{\cal W}^{an}= {\rm dim}_{\R} {\cal O}_{CR} -{\rm
dim}_{CR} M$, such that each continuous CR function on ${\cal O}_{CR}$
admits a holomorphic extension to ${\cal W}^{an}$.  This justifies
that we work with these concepts.  

Basically, it is this extension phenomenon which allows us to derive
non-local results. So, in the Theorems 2, 3 and 4, one assumes explicitely 
that $M$ and $M\backslash E$ or $M\backslash N$ are globally minimal.

\smallskip
{\it e. ${\cal W}$-removability.} 
Now, we need to explain the notion of removability which generalizes
the one-sided removability of the hypersurface case and which will
be our major tool in proving our $L^1$-removability results.

In $\C^{m+n}$,
one-sided neighborhoods are replaced by {\it wedges}, {\it i.e.} open
sets in $\C^{m+n}$ of the form ${\cal W}_p=\{z+\eta\ \! {\bf :} \ \! 
z\in U, \  \eta\in
C\}$, for open $p\in U\subset M$ and open truncated cone $C\subset
T_p\C^{m+n} \cap {\cal V}_{\C^n}(p)$.
Let $\Phi \subset M$ be
a closed set, $p\in b\Phi$. Then an open connected set
${\cal W}_0={\cal W}_0(M\backslash \Phi)$ is called a {\it wedge attached to} $M\backslash \Phi$
if there exists a continuous section $\eta: M\to T_M\C^{m+n}:=
T\C^{m+n}|_{M} / TM$ of the normal
bundle to $M$ and ${\cal W}_0$ contains a wedge at $(p,\eta(p))$, 
for all $p\in M\backslash \Phi$.
This is a one-sided neighborhood of $M\backslash \Phi$
in case $n=1$.
Thus a closed set $\Phi\subset M$ is called {\it ${\cal W}$-removable}
if, given a wedge ${\cal W}_0={\cal W}_0(M\backslash \Phi)$ 
attached to $M\backslash \Phi$, 
there exists a wedge ${\cal W}={\cal W}(M)$ attached to $M$ with holomorphic
functions in ${\cal W}_0$ extending holomorphically to ${\cal W}$
[ME2], [MP2]. 

We localize these notions as follows. We say that a point $p\in \Phi$ (here 
$\Phi$ has no interior points) is ${\cal W}$- or $L^1$-removable
if there exists a small neighborhood $V$ of $p$ in $M$ such that
$L_{CR}^1$ or wedge extendable functions over $M\backslash \Phi$
extend to be $L_{CR}^1$ or wedge extendable over $V$.

\smallskip
{\it f. Description of the proofs.}
As the proofs to be explained in the following section are
long and complicated we conclude the introduction by explaining
the underlying scheme of the proof Theorem 4. 

Subsequently,
we will briefly mention the difference to the proof of 
Theorem 2.

\smallskip
 
{\it Step 1: Reduction to a single CR orbit.}
Another version of Theorem 2 would be

\smallskip
{\sc Theorem 2'.} {\it Let 
$M$ be ${\cal C}^{2,\alpha}$, $0<\alpha <1$, ${\rm dim}_{CR} M=m
\geq 1$. Let $E$ be a closed subset of $M$ such that $H_{loc}^{d-3}(E)<\infty$.
If, for almost every CR orbit ${\cal O}_{CR}$, ${\cal O}_{CR}
\backslash ({\cal O}_{CR} \cap E)$ is globally minimal, then
$E$ is $L^1$-removable.}

\smallskip
Thanks to Theorem 1, it is sufficient to remove
$E \cap {\cal O}_{CR}$ for almost
all CR orbits ${\cal O}_{CR}$ (it can be easily shown that
$H^{e-3}_{loc}(E\cap {\cal O}_{CR})<
\infty$ for almost all CR orbits ${\cal O}_{CR}$, 
$e={\rm dim}_{\R} {\cal O}_{CR}$).

Replacing $M$ by any
${\cal O}_{CR}=$ new $M$, and $E$ by 
$E \cap {\cal O}_{CR}$, we can assume that $M$ is globally
minimal in Theorem 2', {\it i.e.} prove only Theorem 2.

\smallskip
 
{\it Step 2: Wedge extension over $M\backslash \Phi$.} 
Here we intend to apply the following
finest possible {\it extension theorem} (solution 
of Tr\'epreau's conjecture) to $M\backslash \Phi$:

\smallskip
 
{\sc Theorem.} (\cite{TR2}, \cite{TU3}, \cite{ME1}, \cite{JO3}.) {\it
If $M$ is a globally minimal locally embeddable generic
${\cal C}^{2,\alpha}$ manifold, there exists a wedge 
${\cal W}_0$ attached to $M$ such that ${\cal C}_{CR}^0(M)$, 
$L_{loc,CR}^1(M)$, ${\cal D}_{CR}'(M)$, extend holomorphically
to ${\cal W}_0$.}

\smallskip
 
The desired application is possible as soon as we know that 
$M\backslash \Phi$ {\it also} is a single orbit.
Even if $\Phi$ is contained in a codimension $2$ submanifold,
it is in general not true that global minimality of $M$
implies global minimality of $M\backslash\Phi$, as shown by the
following example.

\smallskip
 
{\sc Example 2.} A typical obstruction is where $\Phi=bS$ bounds
a proper closed CR manifold $S\subset (M\backslash \Phi)$, with 
${\rm dim}_{CR} S = {\rm dim}_{CR} M$ and
$\overline{S}=S\cup \Phi$, $\Phi$ a smooth submanifold of $M$
with ${\rm dim}_{\R} \Phi= {\rm dim}_{\R} S-1$ and $M\backslash 
(S\cup \Phi)$, and $S$ is a single CR orbit of $M\backslash \Phi$.
For instance, 
in $\C^2_{w,z}$, $(w,z)$, $w=u+iv$, $z=x+iy$, 
the hypersurface $M: \ y=w\bar{w}\varphi(u)+xw\bar{w}$, 
where $\varphi\equiv 0$ on $\{u\leq 0\}$, $\varphi(0)=0$,
$\varphi>0$ on $\{u>0\}$ has $\Phi=\{u=y=x=0\}$, $S=\{y=x=0, u<0\}$,
$M\backslash (S\cup \Phi)$ is a single orbit of $M\backslash \Phi$.

\smallskip
In the situation of Theorem 4 however, we will be able
to realize our strategy by establishing the
following sufficient technical condition in order that 
CR orbits of $M$ are in one-to-one correspondence
with those of $M\backslash N$.

\smallskip
 
{\sc Theorem 6.} {\it Let $M$ be a locally embeddable
${\cal C}^2$-smooth CR manifold with 
${\rm dim}_{CR} M =m \geq 1$, let $N \subset M$ be a
${\cal C}^2$-smooth submanifold with 
$\hbox{codim}_M N \geq 2$ and 
$T_p N \not \supset T_p^c M$, $\forall \ p\in N$, let
$N^c=\{p\in N\ \! {\bf :} \ \! \hbox{dim}_{\R} T_p^c M 
\cap T_p N = 2m-1\}$ and let
$\Upsilon$ be the set of ${\cal C}^1$ 
sections  of $T^cM$ such that $Y|_{N^c}$ 
is tangent to $N$. If for each $p\in N$, 
${\cal O}_{\Upsilon}(M,p)$ is not contained in $N$, then 
every CR orbit of $M\backslash N$ 
is given by ${\cal O}_{CR} \backslash N$, 
for some CR orbit of $M$.}

\smallskip
This condition that $M$ and $M\backslash \Phi$ too 
are globally minimal ({\it cf.} Theorems 2, 3, 4) cannot be
dropped at least for the problem of wedge extension
of ${\cal D}_{CR}'(M\backslash \Phi)$. A modification
of Treves' construction \cite{TREV2} yields counterexamples
as follows (we however do not know how to 
refine the example up to $L_{CR}^1(M\backslash \Phi)$ or
${\cal C}_{CR}^0(M\backslash \Phi)$).

The occurence of non open CR orbits ${\cal O}$ in 
$M\backslash N$ which are closed submanifolds of $M\backslash N$ 
(in a localized situation near $p\in M$ as in Example 2) can give rise
to CR distributions with support on them. We assume $N=b {\cal O}$.

Indeed, since ${\cal O}$ is a closed CR submanifold
of $M\backslash N$ with the same
CR dimension as $M$, the pullbacks to ${\cal O}$ of $dw_1,...
dw_m$, $dz_1,...,dz_d$ span a subbundle $T_{{\cal O}}'$
of $\C T^* {\cal O}$, of rank $r:= m+e$, where $e={\rm dim}_{\R}
{\cal O}-2m$. Let $\omega_1,...,\omega_{m+e}$ be a basis of 
$T_{{\cal O}}'$, 
for some neighborhood
$U$ of $p$ in $M$. Consider the linear functional 
$u$ on ${\cal C}_c^{\infty}(U
\backslash N)$ defined by
$$u(\varphi) = \int_{(U\backslash N)\cap {\cal O}} \varphi
\ \omega_1\wedge \cdots \wedge \omega_r \wedge d\bar{z}_1 
\wedge \cdots \wedge d\bar{z}_n.$$
This is a nonzero CR distribution in $U\backslash N$ with 
support in ${\cal O}$ since the antiholomorphic sections
$L_j$, $j=1,...,m$ of $T^{0,1}M$ which form the
dual basis to $d\bar{z}_1,...,d\bar{z}_n$, {\it i.e.} 
$L_j \bar{z}_k =\delta_{jk}$, $1\leq j,k\leq m$, act as follows
$$ (L_j \varphi) \omega_1 \wedge \cdots \wedge \omega_r \wedge
d\bar{z}_1 \wedge \cdots \wedge d\bar{z}_n=\pm
d(\varphi \omega_1 \wedge \cdots \wedge \omega_r\wedge  d\bar{z}_1 \wedge
\cdots \wedge d\bar{z}_{j-1} \wedge d\bar{z}_{j+1} \wedge \cdots
\wedge d\bar{z}_n),$$
and Stokes' theorem yields
$$ \left< L_j u, \varphi \right> =\pm
\int_{(U\backslash N)\cap {\cal O}}
d(\varphi \omega_1 \wedge \cdots \wedge \omega_r \wedge d\bar{z}_1 \wedge
\cdots \wedge d\bar{z}_{j-1} \wedge d\bar{z}_{j+1} \wedge \cdots
\wedge d\bar{z}_n)=0.$$
Such a nonzero CR distribution $u$ cannot be the boundary value 
of a holomorphic function in a wedge of edge $M$ at $p$, because
of the uniqueness principle and because $u$ vanishes
on the open set $U\backslash (N\cup {\cal O})$. In this example,
$u$ has the manifold $\overline{\cal O}\backslash {\cal O}=N$
as a singular set near $p$ and $TN \cap T^c M|_N$ is a vector bundle
of rank $2m-1$.

\smallskip
 
{\it Step 3: ${\cal W}$-removability of $\Phi$.} 
Here we use massively the theory of Bishop discs.
It will turn out to be possible to construct sufficiently rich families
of analytic discs which sweep out (almost all of) a wedge attached 
to $M$ and whose boundaries are contained in $M\backslash\Phi$.

Next we would like to extend CR-functions along the discs.
Most often this is done by an application of the Baouendi-Treves
approximation theorem. Unfortunately, because of the presence
of the singularity $\Phi$, this tool is not valid. 

To avoid confusion,
let us mention that in the case that $m\geq 2$, there is the following
device used in \cite{JO3} to employ nevertheless the approximation
theorem: Assume that $N\subset M$ is generic ${\cal C}^{2,\alpha}$
of codimension one 
(hence ${\rm dim}_{CR} N =m-1\geq 1$), and that $N$ has already
been minimalized at $p\in b\Phi$,
where $\Phi\subset N$ is closed, $\Phi\neq N$ and $b\Phi$ is 
with respect to the topology of $N$.
Therefore there is a family of
Bishop discs sweeping out an open wedge ${\cal W}_0$ attached to $N$.
By continuity of Bishop's equation we get perturbed families
of discs attached to all nearby manifolds $N_t$, where
$\bigcup_t N_t$ is a local real foliation of $M$ with $N_0=N$.
Now we can apply the approximation theorem to each $N_t,t\not= 0$
to get extensions of continuous CR functions on $M\backslash N$ 
to the wedges ${\cal W}_t,t\not= 0$ glued to $N_t,t\not= 0$. 
As $p\in b\Phi$, all extensions will glue together in the 
wedge ${\cal W}=\bigcup_{t\not=0}$ whose edge is a neighborhood
of $p$ in $M$.

But in case of general CR dimension $m\geq 1$,
$\Phi$ cannot be included in a generic minimal $N$
of positive CR-dimension. So we will follow an alternative
strategy.

First we deform $M$ in a very nearby manifold
$M^d\supset \Phi$ inside the wedge
${\cal W}^{0}(M\backslash \Phi)$. Thereby, ${\cal W}_0$ becomes
a neighborhood $\omega$ of $M^d\backslash \Phi$. Thus
we are in position to use the {\it continuity principle} with discs with 
boundaries in $M^d\backslash \Phi \cup \omega$. Under
the stated conditions on $N, E$, there are enough deformations of discs
to get that ${\cal H}(\omega)$ extends to 
${\cal H}({\cal W}^d)$, where ${\cal W}^d$ is a wedge attached to 
$M^d$. The main tool is Tumanov's theorem on deformations
of analytic discs which has been successfully applied to propagation 
of CR extension in \cite{TU3}.

Finally, we realize that the construction depends smoothly on $d$ 
and not on the size
of $\omega$, such that, letting $d\to 0$, {\it i.e.} $M^d\to M$,
we get a wedge ${\cal W}$ attached to $M$, {\it i.e.} $\Phi=N,E$ is 
${\cal W}$-removable.

\smallskip
 
{\it Step 4. $L^{\rm p}$-removability.} The general principle is that
${\cal W}$-removability ({\it i.e.} a hull result) 
yields $L^{\rm p}$-removability
provided the constructions of discs glued to $M\backslash \Phi$ are
suitable to control the wedge extension in $L^{\rm p}$ norm, up
to the edge $M$ ({\it cf.} \cite{AC}). So that the link between hulls
and $L^{\rm p}$ is clear, now ({\it cf.} \cite{JO4}, \cite{PO2}).

In fact, to prove our results in $L^{\rm p}$ above is equivalent to 
do the constructions which prove ${\cal W}$-removability, 
with $L^{\rm p}$ control (however, 
in case $N$ is generic of codimension $\geq 2$ 
and $m\geq 2$, we give an argument like
estimating $\int_{{N}^{\varepsilon}}f{\:\!}^{\tau}P 
(\varphi\chi_{\varepsilon})\to 0$, 
$\varepsilon\to 0$, to remove $N$ in $L^{\rm p}$,
see Proposition 1.17. This argument fails in case $m=1$).

\smallskip
The paper is organized as follows. 
In Section 1, we will explain
the main geometric approach and prepare
a good deal of technical tools necessary for the sequel.
This section is already
sufficient to understand the proof of the following 
statement, which apparently was not known before:
{\it isolated points in a real ${\cal C}^2$-smooth 
hypersurface in $\C^2$ are $L^1$-removable.}

In section 2, we prove Theorems 2, 4 and 5  
by using the technique
of deformation of discs and the 
continuity principle. 

In Section 3, we prove
Theorem 1, by combining our techniques studied 
before with a result of Shiffman about 
separate meromorphicity.

In Section 4, we derive the sufficient condition
given in Theorem 3 
insuring that CR orbits of $M$ are in one-to-one correspondence
with those of $M\backslash N$.

In Section 5 we prove Theorem 5 and in Section 6, we explain how these 
results extend to hypoanalytic structures
(Treves \cite{TREV2}). 

This paper was written
when Egmont Porten was in a 
postdoctoral position at the University of Paris 6. He
wishes to thank Professors 
Dolbeault, Henkin, Skoda and Tr\'epreau for their hospitality.

The authors are grateful to Burglind J\"oricke
who draw their attention to the 
subject of $L^{\rm p}$-removability. They wish 
to thank the referee for valuable suggestions.

\bigskip
 
{\bf 1. Preliminaries.} {\sc 1.0. Notations.} 
Let $M$ be generic in $\C^{m+n}$, 
$m={\rm dim}_{CR} M$, $n=\hbox{codim}_{\C^n} M$, 
$d=2m+n=\hbox{dim}_{\R} M$,
$N\subset M$, a submanifold, 
$L_{loc,CR}^1=L^1_{loc}\cap {\cal D}_{CR}'$, 
${\cal O}_{CR}$ be a CR-orbit, 
${\cal O}_{CR}(M,p)$ be the CR orbit of $p$ and
${\cal H}(U)$ be the holomorphic 
functions in $U$. We will give a synthetic proof of the following.

\smallskip
 
{\sc Theorem 1.1.} ({\sc J\"oricke, Porten}, 
\cite{JO3}\cite{JO4}\cite{PO2}.)
{\it Let $M$ be a ${\cal C}^3$-smooth locally embeddable
abstract CR manifold. Then a function 
$f\in L_{loc}^{\rm p}(M), {\rm p}\geq 1$,
is in $L_{loc,CR}^{\rm p}(M)$ if and only if 
$f|_{\cal O}$ is in 
$L_{loc,CR}^{\rm p}({\cal O})$ for 
almost all CR orbits ${\cal O}$ in $M$.}

\smallskip
By almost all CR orbits, we mean except on a CR invariant
set of zero $d$-dimensional measure, $d=\hbox{dim}_{\R} M$. 
Since $L_{loc}^{\rm p} (M)$ embeds in $L_{loc}^1(M)$, we set
${\rm p}=1$.

\smallskip
{\it Remark.} 
For integrable functions Theorem 1.1 is proved 
in \cite{PO2} for M of class ${\cal C}^3$. We give here a summary of proof, the
complete proof appears in \cite{PO2}.
The smoothness ${\cal C}^3$ 
is essentially needed to assure that the CR orbits are of 
class ${\cal C}^2$. On the other hand it is semi-known that for an embedded
CR-manifold of class $C^{2,\alpha}$ the local and global
orbits are of class 
${\cal C}^{2,\beta}$ for any $\beta <\alpha$.
Careful examination shows that the arguments carry
over to our situation assuming
only that $M$ is ${\cal C}^{2,\alpha}$. 

A sharper result (${\cal C}^2$ instead
of ${\cal C}^{2,\alpha}$ or ${\cal C}^3$) for hypersurfaces
relies on the fact that in this case orbits are either open
or even analytic hypersurfaces (\cite{JO4}).

\smallskip
First, we recall the following elementary facts about 
CR orbits. Let $T^cM$ be the ${\cal C}^{k-1}$
complex tangent bundle to $M$ of class ${\cal C}^k$, $k\geq 2$.
Let $L\in \Gamma(U,T^cM)$ be a section over an open
set $U\subset M$. Let $p\in U$. Denote by $t\mapsto L_t(p)$ the integral
curve of $L$ starting at the point $p$, $L_0(p)=p$,
$\frac{d}{dt} L_t(p) =L(L_t(p))$. It is known that the
global flow mapping $(p,t)\mapsto L_t(p)$ is of 
class ${\cal C}^{k-1}$ (defined on some domain $\Omega_L \subset
U\times \R$). Two elements $p,q\in M$ are {\it in the same
CR orbit} if there exist $k\in \N_*$, $L^1,...,L^k$, 
$t_1,...,t_k\in \R$ such that $q=L_{t_k}^k\circ\cdots\circ L_{t_1}^1(p)$.
We write $q=L_T(p)$. An {\it immersed
submanifold} $S$ of $M$ is an abstract manifold
$S$ together with a smooth immersion $S\to M$.
$S$ is CR-{\it integral} if $TS\supset T^cM|_S$.

\smallskip

{\sc Theorem.} ({\sc Sussmann}, \cite{SU}.) {\it Let 
$M$ be ${\cal C}^{k,\alpha}$, $k\geq 2$,
$0\leq \alpha <1$. Then every
CR orbit ${\cal O}_{CR}$ of $M$ admits a unique differentiable
structure such that ${\cal O}_{CR}$ is an immersed
${\cal C}^{k-1,\alpha}$ submanifold of $M$.}

\smallskip
{\it Remark.} In our case that $M$ is locally embeddable CR, a better
regularity result holds (Tumanov):
orbits are of class ${\cal C}^{k, \alpha-0}=
\bigcup_{\varepsilon > 0}  {\cal C}^{k,\alpha-\varepsilon}$.
Indeed, for local CR
orbits, this follows from \cite{TU2}. For global CR
orbits, this follows from \cite{TU5}. Furthermore, 
the estimates for Bishop's equation in \cite{TU5} {\it also}
yield that the {\it transversal structure to orbits}
is ${\cal C}^{k, \alpha-0}$, chains of analytic discs 
replacing Sussmann's curves. Therefore, in Proposition 
below, the map $\psi$ can be assumed to be ${\cal C}^{2,\beta}$ if
$\psi$ is constructed from chains of attached analytic discs 
instead of coming from Sussmann's chains. As we would only need
such $\psi$ to be ${\cal C}^2$ in our proof, this entails that 
our Theorem 1.1 is true in ${\cal C}^{2,\alpha}$, $0<\alpha < 1$. $\square$

Let us denote the mapping
$T\mapsto L_T(p)$ by $\gamma_{L,p}$. 

We now formulate
Sussmann's main lemma on orbits.

\smallskip

{\sc Lemma.} ({\sc Sussmann}, \cite{SU}.) {\it
Let ${\cal O}_{CR}$ be an orbit of $M$. Let 
$p\in {\cal O}_{CR}$. Then there exist
$k\in \N_*$, $L^1,...,L^k\in \Gamma(T^cM)$, $T^0=
(t_1^0,...,t_k^0)\in \R^k$ such that $p=L_T(p)$ is defined
and $(d\gamma_{L,p})_{T^0}(\Omega_{L,p})=T_p{\cal O}_{CR}$.}

\smallskip
In other words, the differential
at $T=T^0$ of the multitime parameter variation $T\mapsto L_T(p)$
covers the tangent space to ${\cal O}_{CR}$ at $p$.

\smallskip

{\sc Corollary.} {\it Let $e={\rm dim}_{\R} {\cal O}_{CR}$.  Then
there exists an $e$-dimensional affine subspace $H\subset \R^k$
through $T^0$ and a ball $B$ of center $T^0$ such that $B\cap H \ni T
\mapsto L_T(p)\in {\cal O}_{CR}\cap {\cal V}_M (p)$ is a ${\cal
C}^{k-1,\alpha}$ embedding.}

\smallskip Now, choosing a small manifold $\Lambda$ through $p$
transversal to ${\cal O}_{CR}$ at $p$, {\it i.e.} $T_p \Lambda\oplus
T_p {\cal O}_{CR}= T_p M$, we recall that the local transversal
structure of orbits at $P$ can be described by applied the multitime
map $\gamma_{L,p}$ to points $q\in \Lambda$ close to $p$ to derive :

\smallskip

{\sc Proposition.} {\it Let $p\in M$ and $e={\rm dim}_{\R} {\cal
O}_{CR}(M,p)$. Then there exist ${\cal C}^{k-1,\alpha}$ coordinates
$\psi=(s,t): V\to I^e\times I^{d-e}$ $(I=(-1,1))$ on a neighborhood
$V$ of $p$ in $M$ such that

$1)$ \ $\psi^{-1}(I^e\times \{0\})={\cal O}_{CR}(M,p)\cap V${\rm ;}

$2)$ \ $\forall t\in I^{d-e}$, the set $\psi^{-1} (I^e\times \{t \})$
is contained in a single CR orbit of $M$.}

\smallskip As a corollary, the mapping $q\mapsto {\rm dim}_{\R} {\cal
O}_{CR} (M,q)$ is lower semicontinuous (all of this is classical). We
use this proposition in 2) just below.

\smallskip
Proposition above provides a local slice structure for
global orbits. The main idea in proving Theorem 1.1 is to

\smallskip {\it Proof of Theorem 1.1.} Assume $f$ is $L^1$-locally CR
on $M$. There exists $\bigcup_{j\in J} U_j = {\cal U}$ a countable
open covering of $M$ by relatively compact $U_j\subset M$ with the
properties that

\smallskip

1) there exist embeddable open sets $U_j' \supset U_j$ such that
$L^1_{CR}(U_j')$ is approximable, in the $L^1$ sense, by holomorphic
polynomials on $\overline{U}_j$;

2) each $U_j$ is equipped with ${\cal C}^2$ coordinates $(s,t)\in U_j=
S_j \times T_j$ such that for every $t\in T_j$, $S_j\times \{t\}$ is
contained a single CR orbit of $M$;

3) for every point $p\in M$, there is a $j\in J$ such that
$p=(s_p,t_p) \in U_j$ and such that the leaf $S_j\times \{t_p\}$ is an
open subset of ${\cal O}_{CR}(M,p)$;

\smallskip

In fact, $U_j$ is locally described by full pieces of the piecewise
smooth flow of $L$, each of which is contained in a single CR orbit of
$M$ \cite{SU} \cite{JO3}.  Regularity of the flow is ${\cal C}^2$, so
coordinates $(s,t)$ are ${\cal C}^2$.  According to 1), for each $f\in
L_{loc,CR}^1(M)$, there exists a sequence $(P_{j,l})_{l\in \N}$ of
holomorphic polynomials such that $\lim_{l\to\infty}
|f-P_{j,l}|_{L^1(U_j)}=0$.  Choose coordinates $(s,t)$ on
$U_j=S_j\times T_j$ and let $\hbox{dim}_{\R} S_j=2m+k$, {\em i.e.}
$p_j\in {\cal A}_k$ :=$\{p\in M\ \! {\bf :} \ \! \hbox{dim}_{\R} {\cal
O}_{CR}(M,p)=2m+k\}$.  Since ${\cal A}_k \cap U$ is a measurable
subset of $S_j \times T_j$ of the form $S_j \times T_{j(k)}$, with
$T_{j(k)} = ( \{ 0 \} \times T_j) \cap {\cal A}_k$ closed in $T_j$,
and since $f|_{{\cal A}_k \cap U_j} \in L^1({\cal A}_k \cap U_j)$,
Fubini's theorem implies (passing to a subsequence if necessary) that
$P_{j,l}$ converges in $L^1(S_j\times \{t\})$ to $f|_{S_j \times
\{t\}}$ for almost all $t\in T_{j(k)}$. Therefore, $f|_{S_j\times
\{t\}}$ is CR, for almost all $t\in T_{j(k)}$. Fix $k$. If $U_j\cap
U_i \neq \emptyset$, $f|_{S_{i,j} \times \{ t_{i,j} \} }$ is in
$L_{CR}^1(S_{i,j}\times \{t_{i,j} \})$ for almost all $t_{i,j}\in
T_{i,j(k)}$. Then the set of orbits ${\cal O}_k$ contained in ${\cal
A}_k$ such that $f|_{{\cal O}_k \cap U_j} \in L_{CR}^1({\cal O}_k \cap
U_j)$ for each $j\in J$ consists of almost all of ${\cal A}_k$, since
$J$ is countable. For such orbits, $f\in L_{loc, CR}^1({\cal O}_k)$.

Conversely, assume $f\in L_{loc}^1(M)$ is CR on almost all CR orbits
of $M$.  The sets ${\cal B}_k=\{p\in M\ \! {\bf :} \ \!
\hbox{dim}_{\R} {\cal O}_{CR}(M,p) \geq 2m+k \}$ are open in $M$ ({\em
cf.} 2) above), $0\leq k \leq n$.  Then $f|_{{\cal B}_n}$ is CR on
${\cal B}_n$, since ${\cal B}_n= {\cal A}_n$.  By induction and
localization, it is sufficient to prove that $f\in L_{CR}^1(U_j)$,
given that $f$ is in $L_{CR}^1({\cal B}_{k+1})$, where $U_j=U$ is an
open subset of $M$ as in 2), $p_j\in {\cal A}_k$, $0\leq k \leq
n-1$. Drop $j$, set $F = (\{0\}\times T) \cap {\cal A}_k\subset T$ and
let $(s,t)$ be ${\cal C}^2$ coordinates on $U=S\times T$.

Let $\delta(t)$ be a so-called regularized distance function
associated with $F$, {\em i.e.}  a real function $\delta \in {\cal
C}^0(T) \cap {\cal C}^{2}(T\backslash F)$ satisfying $1/A \
\hbox{dist}(t,F) \leq \delta(t) \leq A \ \hbox{dist}(t,F)$ and
$|\nabla \delta| \leq A$ on $T\backslash F$, for some constant $A>0$
(\cite{STEIN2}).  Let $\psi(x)$ be a ${\cal C}^2$ real-valued function
with $\psi(x) \equiv 1$ on $\{x \leq 1/2\}, \psi(x)\equiv 0$ on
$\{x\geq 1\}$, set $\psi_l(x)=\psi(lx)$ for $l\in \N \backslash \{0\}$
and put $\chi_l(t)=\psi_l(\delta(t))$.  Then $\chi_l \equiv 1$ on $F$,
$\chi_l$ is ${\cal C}^2$-smooth on $T$, $\hbox{supp} \ \chi_l$ is of
order $1/l$ close to $F$ and $|\nabla \chi_l|\leq C l$, for some
constant $C>0$.

Let $L$ be a ${\cal C}^1$ CR vector field on $U$ and let ${\:}^{\tau}
\! L$ denote the transpose of $L$ with respect to coordinates $(s,t)$
on $U$.  Notice that $L$ has ${\cal C}^1$-smooth coefficients, since
$(s,t)$ are ${\cal C}^2$, so ${\:}^{\tau} \! L$ also is ${\cal C}^1$.
Let $\varphi \in {\cal C}^2_c(U)$. We have to prove that $\int_U f
{\:}^{\tau} \! L(\varphi) \ ds dt =0$.  Write ${\:}^{\tau} \! L
\varphi= {\:}^{\tau} \! L((1-\chi_l)\varphi)+ {\:}^{\tau} \! L(\chi_l
\varphi)$.  Since $f$ is CR on the open set $U\backslash ({\cal A}_k
\cap U) \subset {\cal B}_{k+1}$, we have $$\int_U f {\:}^{\tau} \!
L((1-\chi_l)\varphi) \ ds dt=0.$$ Now, write ${\:}^{\tau} \! L=D_t
+D_s$, with $D_t$ involving only derivatives with respect to $t$ and
no zero order terms. Then $$| \int_U f(D_t+D_s)(\chi_l\varphi) \ dsdt|
\leq | \int_U f D_t (\chi_l \varphi) \ ds dt | + |\int_U f \chi_l
D_s(\varphi) \ ds dt|.$$ Letting $l \to \infty$, we get that the last
integral tends to $\int_F dt(\int_{S\times\{t\}} f {\:}^{\tau} \!
L\varphi \ ds)$, which is zero, since $L|_{S\times \{t\}}$ is a CR
vector field tangent to $S\times \{t\}$ for each $t\in F$, and since
we assumed that $f$ is CR on almost all orbits in ${\cal A}_k$.

On the other hand, for each $\varepsilon >0$, there exists $l\in \N$
such that $|\int_UfD_t(\chi_l \varphi) \ ds dt| < \varepsilon$.
Indeed, the ${\cal C}^1$-smooth coefficients of $D_t$ are bounded by
$C_1 \hbox{dist} (., F)$, since $L|_{S\times \{t\}}$ is tangent to
$S\times \{t\}$ for each $t\in F$. So $|D_t(\chi_l\varphi)|\leq C_2$
on $U$.  Moreover, we can write $$|\int_U f D_t(\chi_l \varphi) \ ds
dt| \leq | \int_S ds \int_{F^{(1/2l)}} f D_t(\chi_l \varphi) \ dt|+
|\int_S ds \int_{F^{(1/l)}\backslash F^{(1/2l)}} fD_t(\chi_l \varphi)
\ dt|,$$ where $F^{(r)}=\{t\in T\ \! {\bf :} \ \delta(t) <r\}$.  The
first term tends to $\int_{S\times F} fD_t (\varphi) \ ds dt$, which
is zero, as $D_t=0$ on $S \times F$. As the ${\cal C}^1$-smooth
coefficients of $D_t$ are bounded by $C_1 \hbox{dist} (., S \times
F)$, we may estimate the second term uniformly by $C \int_S ds
\int_{F^{(1/l)}\backslash F^{(1/2l)}} |f| \ ds dt$. For $l\to \infty$,
the volume of the domain of integration goes to zero.

The proof of Theorem 1.1 is complete. $\square$

\smallskip
 
{\sc 1.1.}$'$ {\sc Corollary.}  {\it A function $f\in {\cal C}^0 (M)$
is continuous CR on $M$ if and only if, for each CR orbit ${\cal O}$
of $M$, $f|_{{\cal O}}$ is continuous CR on ${\cal O}$.}

\smallskip In the rest of Section 1, we shall prove Theorem 4 (i) in a
special case which simplifies the geometric treatment of the
singularity, namely ${\rm dim}_{CR} M=1$ and $T_pN \cap T_p^cM =
\{0\}$, $\forall \ p\in N$.  This special case is much easier to
handle and we will see an elementary construction of discs providing
removability.  But we need already in this case to apply the {\it
continuity principle}, so a deformation argument will be necessary.
Therefore, during the course of the proof of Theorem 1.2 below, we will
prove our main deformation Proposition 1.16 which works for Theorems
2, 3 and 4. It states, roughly (see 1.7), that ${\cal
C}^0$-removability or ${\cal W}$-removability implies $L^{\rm
p}$-removability after some supplementary work, {\it i.e.} that 
after it is known that a set $\Phi$ is ${\cal W}$-removable,
one can easily deduce that $\Phi$ is $L^{\rm p}$ removable under
circumstances which are all satisfied in our Theorems 2, 3 (\cite{MP2}) and 4.

\smallskip

{\sc Theorem 1.2.} {\it Let $M$ be a ${\cal C}^{2,\alpha}$-smooth
$(0<\alpha <1)$ generic manifold in $\C^{1+n}$ with ${\rm dim}_{CR} M
= 1$ and let $N\subset M$ be a ${\cal C}^1$ submanifold such that
$N\cap {\cal O}$ is a ${\cal C}^1$ submanifold of ${\cal O}$ with
${\rm codim}_{\cal O} (N\cap {\cal O})= {\rm codim}_M N$ for every CR
orbit ${\cal O}$ and $T_pN\cap T_p^cM=\{0\}, \forall \ p \in N$.  Then
$N$ is $L^1$-removable, i.e.  $L_{loc,CR}^1(M\backslash N) \cap L^1(M)
= L_{loc,CR}^1(M)$.}

\smallskip {\it Example.} We could choose a local piece $M$ of
a generic manifold with ${\rm dim}_{\R} {\cal O}_{CR}^{loc}(M,0)=3$ in
$\C^{1+n}$ passing by $0$ and as $N\ni 0$ a submanifold transversal to
${\cal O}_{CR}^{loc}$ with ${\rm codim}_M N=3$ (so $T_0^cM \cap
T_0N=\{0\}$).

\smallskip {\it Proof.} After replacing $M$ by an orbit ${\cal O}$ of
$M$ (${\cal O} :=M$), we have that ${\cal O}=M$ and ${\cal O}
\backslash (N\cap {\cal O})$ are globally minimal by our lemma 1.3
below, that $N\subset M$ has ${\rm codim}_M N \geq 3$ and that
$T_pN \cap T_p^cM
=\{0\}$ $\forall \ p\in N$. Hence we can assume thanks to Theorem 1.1
that $M$ is a single orbit.  It is clear that we can also include $N$
in a bigger manifold $N_1$ of codimension exactly three and remove
$N_1$ instead of $N$. So we treat the case ${\rm codim}_M N=3$. 

Now, we prove a capital lemma.

\smallskip
 
{\sc Lemma 1.3.} {\it For each $p_0\in M\backslash N$, ${\cal
O}_{CR}(M\backslash N,p_0)= {\cal O}_{CR}(M,p_0) \backslash N$.}

\smallskip {\it Proof.}  Let $p_1\in N$, let $L$ be a ${\cal
C}^1$-smooth section of $T^cM$ in the neighborhood $U$ of $p_1$, and
let $(t,p)\mapsto L_t(p)$ be the flow of $L$. If $\varepsilon >0$ is
small enough, one has $p_{\varepsilon}=L_{\varepsilon}(p_1)\in
U\backslash N$.  Take $\omega_{\varepsilon}$ a ${\cal C}^1$ piece of
${\cal O}_{CR}(U\backslash N, p_{\varepsilon})$ through
$p_{\varepsilon}$. Then
$L_{-\varepsilon}(\omega_{\varepsilon})\backslash \cup_{t\leq 0}
L_t(N)$ is contained in $L_{-\varepsilon}(\omega_{\varepsilon})$
Observe that $\omega_0=L_{-\varepsilon}(\omega_{\varepsilon})$ is
outside of $L_{-\varepsilon}(\omega_{\varepsilon})$ a CR-submanifold
because $\omega_0\backslash L_{-\varepsilon}(\omega_{\varepsilon})$ is
open in ${\cal O}_{CR}(U\backslash N, p_{\varepsilon})$.  By
continuity, $\omega_0$ is a CR-submanifold and contains some
neighborhood $\omega$ of $p_1$ in ${\cal O}_{CR}^{loc}(p_1)$.

So it is enough to show $\omega_0\backslash N\subset{\cal
O}_{CR}(U\backslash N, p_{\varepsilon})$, which is done by introducing
a second CR vector field $R'$ as a small rotation of $R$. We remark
that $R'$ is tangent to $\omega_0$ and that $\cup_{t\leq 0} L_t(N)
\cap \cup_{t\leq 0} L_t'(N) = N$ in a neighborhood of $p_1$ in $M$
by the fact that ${\rm dim}_{\R} (\R L(p_1) + \R L'(p_1))=2$.
$\square$

\smallskip

{\sc 1.4. Special case.}  We will first remove functions in $L^1(M)$
which extend holomorphically into a neighborhood of $M\backslash N$ in
$\C^{1+n}$ and then prove that a deformation argument reduces to that
situation (see 1.7 and 1.16 below). First, we glue families of
analytic discs to $M$ which shrink as a slice pressing $N$.
Considering functions in $L^1_{CR}(M\backslash N)$ which are
holomorphic in a neighborhood of $M\backslash N$, we can apply a
continuity principle to discs with their boundaries in $M\backslash N$
to {\it modify} $M$ in a new manifold where singularities
disappear. Here is the program.

There exist holomorphic coordinates $(z,w)$ in $\C^{n+1}$ such that
$p_0=0$ in these coordinates, $w\in \C$, $z=x+iy\in \C^n$,
$T_0M=\{x=0\}$, $T_0N=\{x=0,w=0,y_1=0\}$ and, if we set
$y'=(y_2,...,y_n)$, $M$, $N$ are given by \begin{equation} M: \ \
x=h(y,w), \ \ \ \ \ \ \ \ \ \ \ \ \ N: \ \ x=h(y,w), \ y_1=k(y'), \
w=g(y'), \end{equation} for a ${\cal C}^{2,\alpha}$ vector valued function
$h$ with $h(0)=0, \nabla h(0)=0$, and ${\cal C}^1$ functions $g,k$
with $g(0)=0, \nabla g(0)=0$, $k(0)=0, \nabla k(0)=0$. First, we
construct a family of analytic discs such that their boundaries
parametrize a neighborhood of $0$ in $M$. Let $A^{y}_{\rho}(\zeta)=
(Z^y_{\rho}(\zeta),W^y_{\rho}(\zeta))$ be the family of analytic discs
parametrized by $\rho, y$ with $w$-component
$W^y_{\rho}(\zeta)=\rho\zeta+g(y')$ and let $Y^y_{\rho}$ be the
harmonic extension of the solution on $b\Delta$ of Bishop's equation
(see after Lemma 1.10 the theorem of Tumanov for existence)
                    \begin{equation}
                                  Y^y_{\rho} = T(h(Y^y_{\rho}, \rho
\zeta+ g(y')))+ (k(y')+y_1,y').
                    \end{equation} ($T$ denotes the usual Hilbert
transform on the unit circle, {\em i.e.}  $T: L^2(b\Delta) \to
L^2(b\Delta)$ is the unique operator such that $u+iTu$ extends
holomorphically into $\Delta$ and the 
harmonic extension $PH(Tu)(0)=0$ \cite{KOO}.)  Because
of the differentiability of the solution, there exists an open
neighborhood ${\cal Y}$ of $0$ in the $y$-space and $\rho_0 >0$,
$I_{\rho_0}=[0,\rho_0)$, such that the mapping ${\cal Y} \times
I_{\rho_0} \times (\R/2\pi\Z) \ni (y,\rho, \theta) \mapsto
A^y_{\rho}(e^{i\theta})\in M$ is a smooth embedding.  This gives polar
coordinates on $M$ 
and the volume form $\rho d\rho d\theta dy$ is proportional to
$\hbox{dVol}_M$. Fix small $\eta >0$, write $|y|= \max_{1\leq j \leq
n-1} |y_j|$ and set, for $0< \varepsilon << \eta$, $$M_{\eta,
\varepsilon} = \{ p\in M\ \! {\bf :} \ \! |y(p)| < \eta, \rho(p) <
\varepsilon\}.$$ Define also
                    \begin{equation}
                                  M_{\varepsilon}= (M_{\eta,\eta}  \backslash
M_{\eta,\varepsilon}) \cup \bigcup_{|y|<\eta}
A_{\varepsilon}^y(\overline{\Delta})= (M_{\eta,\eta} \backslash
M_{\eta,\varepsilon}) \cup \widetilde{M}_{\varepsilon}.
                                                            \end{equation}
$M_{\varepsilon}$ is a piecewise smooth Lipschitz CR graph, a weakly
differentiable manifold on which the existence of
CR functions still makes sense
({\it cf.} \cite{CHI}, \cite{CS}).  
${\cal V}(E)$ will denote a small open neighborhood of a set
$E\subset \C^{1+n}$.
 
\smallskip
 
{\sc Proposition 1.5.}  {\it If $f\in L_{loc,CR}^{1}(M) \cap {\cal
H}({\cal V}(M\backslash N))$, there exists $f_{\varepsilon}\in
{\cal C}^0_{CR}(M_{\varepsilon})$ such that $f_{\varepsilon}\equiv f$ on
$M_{\varepsilon}\backslash \widetilde{M}_{\varepsilon} (\subset
M)$. Moreover, there exists a sequence $\varepsilon_j \to 0$ such that
$\int_{\widetilde{M}_{\varepsilon_j}} |f_{\varepsilon_j}| \to 0$ as
$j\to \infty$.}

\smallskip We claim that 
Proposition 1.5 implies the theorem in the special case where our function
$f\in L_{loc,CR}^1(M\backslash N)$ is holomorphic in a neighborhood of
$M\backslash N$ in $\C^{1+n}$.  Indeed, take an arbitrary ${\cal
C}^2$-smooth $(1+n,0)$-form $\varphi$ with support in a small ball
centred at $0$.  Then $f_{\varepsilon_j}\to f$ in $L^1$ norm
according to Proposition 1.5 and $$\int_{M_{\eta,\eta}}
f\overline{\partial}\varphi =\lim_{j\to \infty} \int_{M_{\eta,\eta}\backslash
M_{\eta,\varepsilon_j}} f \overline{\partial}\varphi= \lim_{j\to \infty} \int_{M_{\eta,\eta}\backslash
M_{\eta,\varepsilon_j}} f_{\varepsilon_j} \overline{\partial}\varphi=
-\lim_{j\to \infty} \int_{\widetilde{M}_{\varepsilon_j}}
f_{\varepsilon_j} \overline{\partial} \varphi
\stackrel{1.5}{=}0,$$ since, because $f$ is $L^1$
on $M_{\eta,\eta}$, we have $\lim_{j\infty}
\int_{M_{\eta,\varepsilon_j}} f\overline{\partial}\varphi=0$, since
$f\equiv f_{\varepsilon_j}$ on $M_{\eta,\eta} -
M_{\eta,\varepsilon_j}$ and since $f_{\varepsilon_j}\in
{\cal C}^0_{CR}(M_{\varepsilon_j})$. $\square$

\smallskip {\it Proof of Proposition 1.5.} 
We recall facts about analytic isotopies of
discs \cite{ME2}.  Let $M$ be generic. An embedded analytic disc $A$
attached to $M$ is said to be {\it analytically isotopic} to a point
in $M$ if there exists a ${\cal C}^1$-smooth mapping $(s,\zeta)\mapsto
A_s(\zeta), 0\leq s \leq 1, \zeta \in \overline{\Delta}$, such that
$A_0=A$, each $A_s$ is an embedded analytic disc attached to $M$ for
$0\leq s <1$ and $A_1$ is a constant mapping $\overline{\Delta} \to
\{pt\} \in M$.  We shall use the following continuity principle:

\smallskip
 
{\sc Proposition 1.6.} (\cite{ME2})  
{\it Let $M$ be generic, ${\cal C}^{2,\alpha}$, let
$\Phi$ be a proper closed subset of $M$ and let $\omega$ be a
neighborhood of $M\backslash \Phi$ in $\C^{m+n}$. If an embedded disc
 $A$ attached to $M\backslash \Phi$ is analytically isotopic to a
point in $M\backslash \Phi$, then, there exists a neighborhood ${\cal
V}(A(\overline{\Delta}))$ in $\C^{m+n}$ such that, for each function
$f\in {\cal H}(\omega)$, there exists a function $F\in {\cal H}({\cal
V}(A(\overline{\Delta})))$ such that $F=f$ in a neighborhood of
$A(b\Delta)$.} $\square$

\smallskip We then have to isotope the analytic discs to points with
boundaries in $\omega={\cal V}(M\backslash N)$ 
to show that $\zeta \mapsto f \circ
A_{\varepsilon}^y(\zeta)$ extends holomorphically into $\Delta$ for
each $\varepsilon,y$.  One explicit mean to get the isotopy is the
following: we extend the family $A^{y}_{\rho}$ for the 
parameter $\eta \leq y_1 \leq
2\eta$ (and similarly also for $-2\eta \leq y_1 \leq -\eta$) by setting
$A^y_{\rho}(\zeta)=(Z^y_{\rho}(\zeta), W^y_{\rho}(\zeta))$, where
$W^y_{\rho}(\zeta)=\rho/\eta |2\eta-y_1| \zeta + g(y')$ and
$Y^y_{\rho}$ is the solution of Bishop's equation
                    \begin{equation}
                                  Y^y_{\rho}= T (h(Y^y_{\rho},
\rho/\eta | 2\eta -y_1| \zeta + g(y'))) + (k(y')+y_1,y').
                    \end{equation}
Varying $y_1$ in $(0,2\eta)$, we see that
each disc $A_{\varepsilon}^y$ is
analytically isotopic in $\omega$ to
disc $A_{\varepsilon}^{\widetilde{y}}$ with 
$\widetilde{y}_1=2\eta$ and
$\widetilde{y}'=y'$, which is
a point with $(y,w)$-coordinates on $M$ equal to $(k(y')+2\eta,
y', g(y'))$ (due to the fact that
the solution of (5) is a constant, because
the term behind $T$ in independent of $\zeta$).
Therefore, in view of Proposition
1.6, each $f\in {\cal H}(\omega)$ extends to be a {\it continuous} CR
function $f_{\varepsilon}$ on $\widetilde{M}_{\varepsilon}$ as
$f_{\varepsilon}\equiv f$ on $M_{\varepsilon} \backslash
\widetilde{M}_{\varepsilon}$ and $$f_{\varepsilon}(p)=\frac{1}{2i\pi}
\int_{b\Delta} \frac{f\circ A_{\varepsilon}^y(\xi)}{\xi-\zeta_p}
d\xi,$$ for each point $p=A_{\varepsilon}^y(\zeta_p)$. 
Put $f\circ A_{\rho}^y(e^{i\theta})=f(\rho, \theta, y)$. 
For each small
$\varepsilon >0$, there exists $\varepsilon'<\varepsilon$ such that
                    \begin{equation} \varepsilon' \int_{|\theta|\leq
\pi, |y|<\eta} |f(\varepsilon',\theta,y)| \ d\theta dy \leq
\frac{1}{\varepsilon} \int_{\rho < \varepsilon, |\theta| \leq \pi, |y|
<\eta} |f(\rho,\theta,y)| \rho d\rho dy = I_{\varepsilon}
/\varepsilon.
                    \end{equation} On the other hand, there exists a
constant $C>0$, depending only on the ${\cal C}^1$-norm of the embedding
${\cal Y} \times I_{\rho_0} \times I_{2\pi} \ni (y,\rho, \theta)
\mapsto A^y_{\rho}(e^{i\theta})\in M$ with $$J_{\varepsilon'} =
\int_{\widetilde{M}_{\varepsilon'}} |f_{\varepsilon'}| \ d
\hbox{Vol}_{\widetilde{M}_{\varepsilon'}} \leq C (\varepsilon')^2
\int_{\rho< 1,|\theta| \leq \pi, |y|<\eta} |f\circ A_{\varepsilon'}^y
(\rho e^{i\theta})| \ \rho d\rho d\theta dy.$$ Recall that, for 
$u\in {\cal C}^0(\overline{\Delta}) \cap {\cal H}(\Delta)$,
the subharmonicity of $|u|$ implies $\int_0^{2\pi} |u(re^{i\theta})| \
d\theta \leq \int_0^{2\pi} |u(e^{i\theta})| \ d\theta$, for every
$0\leq r \leq 1$.  Then $$J_{\varepsilon'} \leq C(\varepsilon')^2
\int_{|y|<\eta} dy \int_{\rho < 1} \rho d\rho \int_{|\theta| \leq \pi}
|f\circ A_{\varepsilon'}^y(\rho e^{i\theta})| \ d\theta \leq \ \ \ \ \ \ \ \ \ \ \ \ \ \ \ \ $$
$$ \ \ \ \ \ \ \ \ \ \ \ \ \ \ \ \ \leq C
(\varepsilon')^2 \int_0^{1} \rho d\rho \int_{|\theta| \leq \pi, |y|<
\eta} |f\circ A_{\varepsilon'}^y(e^{i\theta})| \ d\theta dy$$ and this
is $\leq C I_{\varepsilon} \varepsilon'/2\varepsilon\leq C 
I_{\varepsilon} /2$, in view of
$(6)$. But $I_{\varepsilon} \to 0$, since $f$ is $L^1$.

The proof of Proposition 1.5 is complete. $\square$

\smallskip
 
{\sc 1.7. General case.} Let $M$ be generic, let $p\in M$ and suppose
that $\Phi\subset M$ is a given closed singularity set for locally
integrable CR functions which are $L_{loc,CR}^1$ on $M\backslash
\Phi$.  Our aim is to {\it reduce the problem of removing $\Phi$ to the
problem of removing $\Phi$ for the space of functions holomorphic in
some neighborhood of $M\backslash \Phi$ and locally integrable on $M$},
${\cal H}({\cal V}(M\backslash \Phi))\cap L_{loc}^1(M)$. Informally, this
can be performed by suitable very small deformations $M^d$ of
$M$ over a wedge attached to $M\backslash \Phi$, so small 
that we only very slightly change a fixed function $f\in L^1_{loc}(M)\cap
L_{loc,CR}^1(M\backslash \Phi)$. 

To be precise, we localize
first $M$ near one of its points $p\in b\Phi$.  Suppose that we can
prove (see Proposition 1.16) that, given $f\in L^1(M)\cap
L_{loc,CR}^1(M\backslash \Phi)$, for each $\varepsilon >0$, there
exists a ${\cal C}^{2,\alpha}$-smooth compactly supported deformation
$M^d$ of $M$, $M^d\supset \Phi$, with $||M^d-M||_{{\cal C}^{2,\alpha}}
< \varepsilon$ 
(the ${\cal C}^{2,\alpha}$ distance $||M^d-M||_{{\cal C}^{2,\alpha}}$
between $M$ and $M^d$ can be measured by $||h^d-h||_{{\cal C}^{2,\alpha}}$
for graphing functions $h$ and $h^d$ like in (2)) 
 such that there exists $f^d\in L^1(M^d)$ which is
holomorphic in a neighborhood of $(M^d\backslash \Phi)$ in $\C^{m+n}$
and such that $|f^d-f|_{L^1(M)} < \varepsilon$. If, for each $M^d$,
$$L^1(M^d)\cap {\cal H}({\cal V}(M^d \backslash
\Phi))=L_{loc,CR}^1(M^d),$$ then $f\in L_{loc,CR}^1(M)$ near $p$.
Indeed, fix $\varepsilon >0$ and choose $d$, $M^d$, $f^d\in L^1(M^d)$
as above.  Let $\varphi$ be a ${\cal C}^2$-smooth $(m+n,m-1)$-form,
compactly supported near $p$.  Then $|\int_M
(f^d-f)\overline{\partial} \varphi | \leq C_{\varphi} \varepsilon$,
since $|f^d-f|_{L^1(M)}< \varepsilon$ and we can identify abusively
the volume form on $M$ and on $M^d$.  Therefore, $\int_{M^d} f^d
\overline{\partial} \varphi =0$ implies that $\int_M f
\overline{\partial} \varphi =0$, since $\varepsilon$ was arbitrary.

\smallskip To fulfill this program, it is convenient to introduce
families of discs which enjoy better regularity properties than the
general discs of minimal defect \cite{TU1}.
These discs are suitable to apply measure theoretic arguments
in $L^{\rm p}$, {\it e.g.} Fubini's theorem or convergence almost
everywhere and they in a certain sense reduce the $L^{\rm p}$ analysis
on wedges to the $L^{\rm p}$ analysis on bundles of
analytic discs.

\smallskip
 
{\sc Definition 1.8.} By a {\it regular family of analytic discs attached
to $M$ at $p$}, we mean a ${\cal C}^{2,\beta}$-smooth mapping $A: {\cal
S}\times {\cal V}\times \overline{\Delta} \to \C^{m+n}$,
$(s,v,\zeta)\mapsto A_{s,v}(\zeta)$, $A_{0,v}(1)=p$, holomorphic with
respect to $\zeta$, the parameters $s,v$ run over $0$-neighborhoods
${\cal S}\subset \R^{2m+n-1}$, ${\cal V}\subset \R^{n-1}$, respectively,
such that

1) The mapping ${\cal S}\times b\Delta \to M$, $(s,\zeta) \mapsto
A_{s,v}(\zeta)$ is an embedding, $\forall \ v\in {\cal V}$ uniformly;

2) The vector $Y=-\frac{\partial }{\partial \zeta} A_{0,0}(1) \not\in
T_pM$ and the rank of the mapping

$v\mapsto \hbox{pr}_{T_p \C^{m+n}/(T_pM \oplus \R Y)} \left
( -\frac{\partial }{\partial \zeta} A_{0,v}(1)\right)$ is equal to
$n-1$;

3) There exists a neighborhood $V$ of $p$ such that $\{A_{s,v}(\zeta)\
\! {\bf :} \ \! \zeta \in b\Delta\} \subset V$ and $L_{loc,CR}^1$
functions are approximable by holomorphic polynomials in $L^1(V)$.

\smallskip Here, we require that one set of parameters ${\cal V}$ describes
the total amount of necessary outer directions 
$-\frac{\partial }{\partial \zeta} A_{0,v}(1)$
to cover a cone in $T_p\C^{m+n}/T_pM$ of dimension
$n$ and that the set of parameters ${\cal S}$ describes, together
with the arc-length on $b\Delta$ near 1, a set which gives coordinates
on $M$ near $p$. This is crucial to use Fubini's theorem in order
to exhibit $L^1$ traces of $f$ on almost every circle
$A_{s,v}(b\Delta)$ if $f\in L^1(M)$ ({\it cf.} \cite{KYREA},
\cite{JO4}).

\smallskip {\it Remark 1.} We stress that existence of
regular families of analytic discs does not 
hold automatically at minimal points, because a disc
of minimal defect need not be an embedding \cite{TU1}, \cite{BRT} and
that on a globally minimal $M$ which is nowhere locally minimal,
condition 2) above cannot at all be satisfied. Fortunately, the
existence of regular families of analytic discs attached to {\it
deformations of $M$} follows from 
the work of Tumanov on propagation \cite{TU3},
\cite{TU5} (but is nontrivial) and is
{\it equivalent to global minimality}. 

Furthermore, it follows from \cite{TU3} in the 
construction of such $A_{s,v}$ that one deforms
a small round disc whose projection on the
$T_p^cM$ space is a trivial disc, like $\zeta \mapsto 
(c(1-\zeta),0)\in \C^m$, implying that each disc $A_{s,v}$ is an embedding 
({\it cf.} also Lemma 2.4).

We shall use $L^{\rm p}$ wedge extension
of $L^{\rm p}_{CR}$ functions on a globally minimal
$M$, which is proved in \cite{PO1}, \cite{PO2}
using the minimalization theorem in \cite{JO3}.

\smallskip {\it Remark 2.} According to conditions 1) and 2), there
exist ${\cal S}_1\subset {\cal S}$ a subneighborhood of $0$, ${\cal
V}_1 \subset {\cal V}$ a subneighborhood of $0$, $\Delta_1\subset
\overline{\Delta}$ a neighborhood of $1$, such that the mapping
$${\cal S}_1\times {\cal V}_1 \times \stackrel{\circ}{\Delta}_1 \ni
(s,v,\zeta) \mapsto A_{s,v}(\zeta)\in \C^{m+n}\backslash M$$ is an
embedding, whose image $${\cal W}=\{A_{s,v}(\zeta) \in \C^{m+n}; \
(s,v,\zeta) \in {\cal S}_1\times {\cal V}_1 \times
\stackrel{\circ}{\Delta}_1\}$$ constitutes a wedge open set with edge
$M$ at $p$. $\square$

\smallskip {\it Remark 3.}  For sufficiently small families of discs,
the approximation property is automatic, according to a version of a
theorem due to Baouendi and Treves proved by J\"oricke \cite{JO4}.

\smallskip
 
{\sc Theorem 1.9.} \cite{BT} \cite{JO4}.  {\it Let $M$ be a ${\cal
C}^2$-smooth locally embeddable CR manifold, let ${\rm p}$, $1\leq
{\rm p} < \infty$. Then, for every neighborhood $U$ of $p$ in $M$,
there exists another neighborhood $V \subset \subset U$ of $p$ such
that, for each $f\in L^{\rm p}_{loc,CR}(U)$, each $\varepsilon >0$,
there exists a polynomial $P$ in $z_1,...,z_{m+n}$ with $$\int_V |f-P|
\ d\hbox{Vol}_M < \varepsilon. \ \square $$} 

{\it Remark 4.}  In general, we shall deal below with regular families
such that ${\cal S}\times b\Delta \to M$, $(s,\zeta) \mapsto
A_{s,v}(\zeta)$ is an embedding, $\forall \ v\in {\cal V}$.  But the weakened
condition 1) that the disc mapping with parameters should only be an immersion
along the boundary can be reached directly by suitable modifications
of a given disc as follows.  If $E$ is any set in a metric space $B$, we shall
denote by ${\cal V}(E,r)$ the set $\{p\in E\ \! {\bf :} \ \!
\hbox{dist} (E,p) < r\}$ and ${\cal V}(E)$ any small open neighborhood of
$E$ in $B$. For instance the space $B$ of ${\cal C}^{2,\alpha}$
discs is a Banach space and we write
$A'\in {\cal V}( A,\delta)$.

\smallskip
 
{\sc Lemma 1.10.} \cite{KYREA}. {\it Let $M$ be generic, ${\cal
C}^{2,\alpha}$-smooth, and let $A$ be a sufficiently small ${\cal
C}^{2,\beta}$-smooth analytic disc attached to $M$, $A(1)=p_0$. Then,
for each $\delta >0$, there exists a family of analytic discs $A_s$
attached to $M$ with $A_0(1)=p_0$, the parameter $s$ runs over a
$0$-neighborhood ${\cal S} \subset \R^{2m+n-1}$, such that $A_s \in
{\cal V}(A, \delta)$, for $s\in {\cal S}$ and the rank of the mapping
${\cal S}\times \R \ni (s,\theta) \mapsto A_s(e^{i\theta}) \in M$ is
equal to $2m+n={\rm dim}_{\R} M$.}

\smallskip To begin with, let $T$ denote the Hilbert transform {\it i.e.} the
harmonic conjugation operator) on the unit circle $b\Delta$.  Recall
that for a function $\phi$ on $b\Delta$, $T\phi$ is the unique
function on $b\Delta$ such that $\phi + i T\phi$ extends
holomorphically into $\Delta$ and $PH(T\phi)(0)=0$, where $PH$ denotes
the harmonic extension operator. We also denote by $T_1$ the Hilbert
transform vanishing at $1$, $T_1\phi=T\phi-T\phi(1)$.  It is known
that $T$ and $T_1$ are bounded operators ${\cal C}^{k,\alpha} \to
{\cal C}^{k,\alpha}$ $(k\geq 0$, $0<\alpha <1$) and in $L^{\rm p}$
($1< {\rm p} \leq \infty$.

The following equation, called Bishop's equation with parameters,
arises in constructing analytic discs with boundaries in a given
generic manifold
                    \begin{equation}
                     y=T_1H(y,.,t)+y_0,
                    \end{equation} where $H$ is a given $\R^n$-valued
function depending on $y\in \R^n$, $\zeta \in b\Delta$ and a parameter
$t\in \R^l$. The solution $y=y(\zeta,t,y_0)$ is a function of $\zeta
\in b\Delta$, the variables $t\in \R^l$ and $y_0\in \R^n$ being
parameters. We shall repeatedly use the best result below on
solvability and regularity of Bishop's equation. Let $B_r^n$ denote
the ball of radius $r$ centered at the origin.

\smallskip
 
{\sc Theorem.} ({\sc Tumanov}, \cite{TU5}). {\it Let $H\in {\cal
C}^{k,\alpha}(B_1^n\times b\Delta \times B_1^l, \R^n)$, $k\geq 1$,
$0<\alpha < 1$. For every constant $C>0$, there exists $c>0$ such
that, if $||H||_{{\cal C}^{k,\alpha}} <C$, $||H||_{{\cal C}^0}< c$,
$||H_y||_{{\cal C}^0}<c$ and $||H_{\zeta}||_{{\cal C}^0} <c$, then

{\rm (i)} \ $(6)$ has a unique solution $\zeta\mapsto y(\zeta,t,y_0)$
in $L^2(b\Delta)${\rm ;}

{\rm (ii)} \ $y\in {\cal C}^{k,\alpha} (b\Delta)$ and $||y||_{{\cal
C}^{k,\alpha}}$ is uniformly bounded with respect to parameters{\rm ;}

{\rm (iii)} \ $y\in {\cal C}^{k,\beta} (b\Delta\times B_1^l \times
B_c^n,\R^n)$ for all $0<\beta <\alpha$.} $\square$

\smallskip {\it Proof of Lemma 1.10.} Assume that $M$ is given by (2)
and write $A(\zeta)=(Z(\zeta),W(\zeta))$. We can slightly perturb the
component $W_m( \zeta)$ to insure that $W_m'(\zeta)\neq 0$ for every
$\zeta\in b\Delta$ and that this property holds for every disc in
${\cal V}(A,\delta')$, $\delta' >0$, $\delta' <<\delta$ small
enough. Still denote by $A$ this disc, where $y(\zeta)$ is obtained by
taking the harmonic extension of the solution of Bishop's equation
$Y=T_1h(Y,W)$ on $b\Delta$. Using the notation
$W^*(\zeta)=(W_1(\zeta),..., W_{m-1}(\zeta))$, we set, for small
$|t|<\varepsilon$, $\varepsilon << \delta$, $\varepsilon>0$,
$$W^t(\zeta)=(W^*(\zeta), W_m(\zeta)+t\zeta W_m'(\lambda \zeta)),$$
where $\lambda <1$ is close enough to 1 to insure that the jacobian of
the map $(t,\theta)\mapsto W_m^t(e^{i\theta})$ is nowhere vanishing on
$\{|t| < \varepsilon\} \times \R$. Then we consider the family of
analytic discs $$A_{t,w^{*0},y^0}(\zeta)= (Z_{t,w^{*0},y^0}(\zeta),
W^*(\zeta)+w^{*0},W_m^t(\zeta)-W_m^t(1)),$$ where $Y_{t,w^{*0},y^0}$
satisfies $$Y_{t,w^{*0},y^0}=T_1h(Y_{t,w^{*0},y^0},
W^*(.)+w^{*0},W_m^t(.)-W_m^t(1))+y^0$$ on $b\Delta$. 
It is ${\cal C}^{2,\beta}$-smooth with
respect to all the variables by the Theorem of Tumanov. 
Then the jacobian matrix $E=D_{y^0}
Y_{t,w^{*0},y^{0}}$ has non vanishing determinant, if the disc
$A_{t,w^{*0},y^0}$ is small enough, since it satisfies
$E=T_1h(E,W_{t,w^{*0},y^0})+I$, and $dh(0)=0$. Since the jacobian
matrix $D_{(w^{*0},t,\theta)} W_{t,w^{*0},y^0}$ has non vanishing
determinant too, we get the lemma by letting $s=(t,w^{*0},y^0)$ vary
in a sufficiently small
 neighborhood ${\cal S}$ of $0$ in $\R^{2m+n-1}$.

The proof of Lemma 1.10 is complete. $\square$

\smallskip From now on, since it will be more convenient and not
restrictive for our purpose, we shall assume that regular families of
discs {\it embed} ${\cal S} \times b\Delta$ in $M$ as
in Definition 1.18. We check:

\smallskip
 
{\sc Lemma 1.11.} {\it Let
 $f\in L_{loc,CR}^1(M)$ and let $v\in {\cal V}$ be given. Then, for
almost all $s\in {\cal S}$, the mapping $b\Delta \ni \zeta \mapsto
f\circ A_{s,v}(\zeta)\in \C$ is well-defined and belongs to
$L^1(b\Delta)$.}

\smallskip {\it Proof.} The mapping ${\cal S} \times b\Delta \ni
(s,\theta) \mapsto A_{s,v}(e^{i\theta})$ is a smooth embedding of
${\cal S} \times b\Delta$ onto a tubular open connected neighborhood
${\cal T}$ of $A(b\Delta)$ in $M$ and this gives coordinates $(s,\theta)$ on
${\cal T}$ together with a volume element $ds d\theta$ on ${\cal T}$ which is
proportional to $d\hbox{Vol}_M$. Since $f\in L^1({\cal T})$, Fubini's theorem
written in the form $$\int_{{\cal S}\times b\Delta} f\circ
A_{s,v}(e^{i\theta}) \ ds d\theta = \int_{{\cal S}} ds \left
( \int_0^{2\pi} f\circ A_{s,v}(e^{i\theta}) \ d\theta \right) <
\infty,$$ implies that $b\Delta \ni \zeta \mapsto f\circ
A_{s,v}(\zeta)$ belongs to $L^1(b\Delta)$ for almost all $s\in {\cal
S}$, which completes the proof. $\square$

\smallskip A $L_{loc,CR}^1$ function $f$ on $M$ is said to be
holomorphically extendable to a wedge ${\cal W}={\cal W}(U,C)$ if
there exists $F\in {\cal H}({\cal W})$ such that $F|_{U_{\eta}} \to f$
in the $L^1$ sense, as $\eta \to 0$, $U_{\eta}=U+\eta$, uniformly for
$\eta \in C$. Recall that this implies that, given a subcone $C_1
\subset \subset C$, given $U_1\subset \subset U$, for every one-parameter
${\cal C}^2$-smooth family of open sets
$U_{1\varepsilon}$ with  $U_{1\varepsilon}\to U_1$ in ${\cal C}^2$
norm  and such that $U_{1\varepsilon} \subset {\cal W}(U_1,C_1)$, then
$F|_{U_{1\varepsilon}} \to f|_{U_1}$ in $L^1$ norm, and that, conversely,
if $F|_{U_{1\varepsilon}} \to f|_{U_1}$ in $L^1$, then
$F|_{{\cal W}(U_2,C_2)}$ is a holomorphic extension of the
$L_{loc,CR}^1$ function $f$ in the first sense, for
possibly smaller $U_2 \subset \subset U_1$, $C_2 \subset
\subset C_1$. Therefore, the sizes of $U$ andXS $C$ are not essential.

Let $H^1_{\rm a}(\Delta)$ denote Hardy space
in $\Delta$, $H_{\rm a}^1(\Delta)=\{u\in {\cal H}(\Delta)\ \! {\bf :} \ \!
\sup_{r<1} \int_0^{2\pi} |u(re^{i\theta})| d\theta < \infty\}$. 
Accordingly, we denote by $H_{\rm a}^1(\Delta)$ the Hardy space
of holomorphic functions in ${\cal W}$ with $L^1$ boundary values on 
$M$, {\it i.e.} $F|_{U_{1\varepsilon}} \to f|_{U_1}$ in $L^1$ norm.

The existence of regular families of discs is well suited to get
holomorphic extension of measurable CR functions. This is the content
of the next proposition.

\smallskip
 
{\sc Proposition 1.12.}  {\it Let $M$ be generic, ${\cal
C}^{2,\alpha}$, let $p\in M$ and assume that there exists a regular
family $A: {\cal S}\times {\cal V} \times \stackrel{\circ}{\Delta}_1
\to \C^{m+n}$ of analytic discs attached to $M$ at $p$.  Then there
exist open neighborhoods $0\in {\cal S}_1\subset {\cal S}$, $0\in
{\cal V}_1 \subset {\cal V}$, $1 \in \Delta_1 \subset
\overline{\Delta}$, such that every $L_{loc,CR}^1$ function on $M$
extends holomorphically as $F\in H_{\rm a}^1({\cal W})$
into the wedge open set $${\cal
W}=\{A_{s,v}(\zeta)\in \C^{m+n}\ \! {\bf :} \ \!  (s,v,\zeta) \in
{\cal S}_1 \times {\cal V}_1 \times \stackrel{\circ}{\Delta}_1\}.$$}

{\it Remark.} As ${\cal W}$ is foliated by pieces of discs $A_{s,v}$,
the overall scheme of our proof will be a reduction to Hardy
spaces of discs $A_{s,v}$. Hence
$H_{\rm a}^1({\cal W})$ appears as a {\it bundle of Hardy spaces on
$\stackrel{\circ}{\Delta}_1$}. Our proof is a repetition of 
\cite{JO4} in higher codimension.

\smallskip
{\it Proof.} Choose ${\cal S}_1 \times {\cal V}_1 \times \Delta_1
\subset {\cal S} \times {\cal V} \times \overline{\Delta}$ such that
$\Delta_1$ is in the form $\Delta_1= \{\zeta=re^{i\theta} \in
\overline{\Delta}\ \! {\bf :} \ \!  r_1 < r \leq 1, |\theta| <
\theta_1\}$, $\theta_1 > 0$ and ${\cal S}_1 \times {\cal V}_1 \times
\stackrel{\circ}{\Delta}_1 \ni (s,v,\zeta)\mapsto A_{s,v}(\zeta) \in
\C^{m+n}$ is an embedding.  Set $I_1=(-\theta_1,\theta_1)$.  For $r_1
< r \leq 1$, $v\in {\cal V}_1$, define a partial copy of $M$ contained
in ${\cal W}$ $$M_r(v)=\{A_{s,v}(re^{i\theta}); \ \theta \in I_1, s\in
{\cal S}_1\}.$$ There exists a uniform constant $C>0$ such that $(1/C)
d\hbox{Vol}_M \leq d\hbox{Vol}_{M_r(v)} \leq C d\hbox{Vol}_M$, where
$d\hbox{Vol}_M(p)=g_M(s,\theta) ds d\theta$,
$d\hbox{Vol}_{M_r(v)}(p)=g_{M_r(v)}(s,\theta) ds d\theta$ are computed
with respect to the ${\cal C}^2$-smooth parametrization
$(s,\theta)\mapsto A_{s,v}(re^{i\theta})$ with respect
to which both $d\hbox{Vol}_M$ and 
$d\hbox{Vol}_{M_r(v)}$ can be expressed, since $M_r(v)$
is a small deformation of $M$.  In the rest of the paper,
we shall denote by $C$ an unspecified constant $>0$, depending on the
context.

\smallskip
 
{\sc Lemma 1.13.} {\it There exists a constant $C>0$ such that the
following estimate holds $$\int_{M_r(v)} |P| \ d\hbox{Vol}_{M_r(v)}
\leq C \int_V |P| \ d\hbox{Vol}_M,$$ for every holomorphic polynomial
$P$, every $r\in (r_1,1]$ and every $v\in {\cal V}_1$.}

{\it Proof.}  Recall $g_{M_r(v)}(p)/g_M(p) \leq C$.  By
plurisubharmonicity, for every polynomial $P$, $$\int_{{\cal
S}_1}\int_{-\theta_1}^{\theta_1} |P\circ A_{s,v}(re^{i\theta})|
g_{M_r(v)}(s,\theta) \ ds d\theta \leq C_1 C \int_{{\cal
S}_1}\int_{-\pi}^{\pi} |P\circ A_{s,v}(e^{i\theta})| g_{M}(s,\theta) \
ds d\theta,$$ which yields, since the disc mapping $(s,\theta)\mapsto
A_{s,v}(e^{i\theta})$ is an embedding (uniform in $v$), 
$$\int_{{\cal
S}_1}\int_{-\theta_1}^{\theta_1} |P\circ A_{s,v}(re^{i\theta})|
g_{M_r(v)}(s,\theta) \ ds d\theta \leq C|P|_{L^1(V)}.$$

The proof of Lemma 1.13 is complete. $\square$

\smallskip Integrate now 1.13 with respect to $v\in {\cal V}_1$ and to
$r\in (r_1,1]$. We get
                    \begin{equation} \int_{{\cal W}} |P| \
d\hbox{Vol}_{\C^{m+n}} \leq C \int_V |P| \ d\hbox{Vol}_M.
                    \end{equation} Let $f$ be a CR function on $M$ of
class $L^1_{loc}$. According to the approximation theorem, there
exists a sequence $\{P_{\mu}\}_{\mu=1}^{\infty}$ of polynomials in
$z_1,...,z_{m+n}$ such that $f=\lim_{\mu\infty} P_{\mu}$ in
$L^1(V)$. Then, for $P_{\mu}-P_{\nu}$, we have the estimate (8), so
$|P_{\mu}-P_{\nu}|_{L^1({\cal W})}$ tends to $0$ as $\mu, \nu$ go to
infinity. Hence $P_{\mu}$ converges to a holomorphic function $F$ in
${\cal H}({\cal W})$ (Cauchy's Theorem).

To exhibit $f$ as the $L^1$ boundary value of $F$, take $M'$ as a
neighborhood of $p$ contained in $$\bigcap_{v \in {\cal V}}
\{A_{s,v}(e^{i\theta})\ \! {\bf :} \ \! s\in {\cal S}_1, \theta\in
I_1\}.$$ Define $M_r'(v)=\{A_{s,v}(r e^{i\theta}); s\in {\cal S}_1,
\theta\in I_1, A_{s,v}(e^{i\theta})\in M'\}$. We have to show the uniform
$L^1$ convergence $F|_{M_r'(v)} \to f|_{M'}$, as $r$ tends to $1$. For
$\varepsilon >0$, 1.9 shows that there is a polynomial $P$ with $|P|_V
-f|_V|_{L^1(V)} < \varepsilon$. 
We notice that $f=\lim_{\mu\infty} P_{\mu}$ in $L^1(U)$ also
implies that $f\circ A_{s,v}(e^{i\theta})$ extends
holomorphically to $\Delta$ for almost all $s\in {\cal S}$, {\it i.e.}
$f\circ A_{s,v}\in H_{\rm a}^1(\Delta)$.
By 1.13, applied to $F-P$, which
is possible, because
$f\circ A_{s,v}\in H_{\rm a}^1(\Delta)$, we get that 
for all $r,v$, the estimate $|P|_{M_r'(v)}-
F|_{M_r'(v)}|_{L^1(M_r'(v))} < C \varepsilon$ holds true. By continuity
of $P$, we can now choose $r_0>0$ such that, for $0<r<r_0$,
$|P|_{M_r'(v)}-P|_{M'}|_{L^1(M')} < \varepsilon$. This implies that
$|F|_{M_r'(v)}-f|_{M'}|_{L^1(M')}$ can be made arbitrarily small
(with the identifications of volume forms on $M'$ and on 
$M_{r}'(v)$).

The proof of Proposition 1.12 is complete. $\square$

\smallskip

{\sc Proposition 1.14.} {\it Under the hypotheses of 1.12, let
$W_2=\{A_{s,0}(e^{i\theta})\in M; \ 2s \in {\cal S}_1, 2\theta \in
I_1\}$.  Then, for each $\varepsilon >0$,
 there exists a ${\cal C}^{2,\beta}$-smooth deformation $M^d\subset M
\cup {\cal W}$ of $M$ with $W_2 \subset \hbox{supp} \ d \subset
\subset V$ such that $||M^d - M ||_{{\cal C}^{2,\beta}}
 < \varepsilon$ and there exists a function $f^d\in
L_{loc,CR}^1(M^d)\cap {\cal H}({\cal V}(W_2^d))$ such that
$|f^d-f|_{L^1(V)}<\varepsilon$ and $f^d\equiv f$ in $M\backslash
\hbox{supp} \ d$.}

\smallskip {\it Proof.} Fix $\varepsilon >0$ and define
$W_1=\{A_{s,0}(e^{i\theta})\in M; \ s\in {\cal S}_1, \theta \in
I_1\}\subset V$, which is a neighborhood of $0$ in $M$. Let $r: {\cal
S}_1 \times I_1 \to (r_1, 1]$ be a ${\cal C}^{2,\beta}$-smooth
function with ${\cal S}_1 \times I_1 \supset \hbox{supp} \ (1-r)
\supset \frac{1}{2} {\cal S}_1 \times \frac{1}{2} I_1$ and consider
the deformation $M^d$ of $M$ with support in $W_1$ for which
$$W_1^d=\{A_{s,0}(r(s,\theta)e^{i\theta}); \ s\in {\cal S}_1, \theta
\in I_1\} \subset {\cal W} \cup M.$$ Let $f\in L_{loc,CR}^1(M)$.
Since $f$ extends holomorphically into the wedge ${\cal W}$ in 1.12,
one can define $f^d$ to be equal to $f$ in $M\backslash \hbox{supp} \
d$ and to be the restriction 
$F|_{W_1^d}$ 
of the holomorphic extension $F$
of $f$ to ${\cal W}$.  By a similar reasoning as in 1.12, assume
$f=\lim_{\nu \infty} P_{\nu}$ in $L^1(V)$ as in 1.9.  Then
$F|_{W_1^d}=f^d=\lim_{\nu\infty} P_{\nu}|_{W_1^d}$ in $L^1(W_1^d)$,
since $\int_{W_1^d} |P_{\nu}-P_{\mu}| \leq C \int_{W_1}
|P_{\nu}-P_{\mu}|$.  Now, we have $$|f^d-f|_{L^1(W_1)}<
|f-P_{\nu}|_{L^1(W_1)}+ |P_{\nu}|_{W_1^d}-P_{\nu}|_{W_1}|_{L^1}+
|P_{\nu}-f^d|_{L^1(W_1^d)}.$$ The first term is estimated in terms of
the approximation theorem. The third term is estimated by Carleson's
embedding theorem.
According to Carleson's embedding theorem (see \cite{KOO}, p. 181),
given a function $r(\theta)\in {\cal C}^1$ such that $r\equiv r_1$ in
$(-\theta_1,\theta_1)$, $0<r_1<1$ and $\hbox{supp} \ (1-r)\subset
(-\theta_0,\theta_0)$, $0<\theta_1 <\theta_0 <\pi$, there exists a
constant $C_1>0$ such that $$\forall \ f \in H^1_{\rm a}(\Delta), \ \ \ \ \
\int_{-\theta_0}^{\theta_0} f(re^{i\theta}) d\theta \leq C_1
|f|_{L^1(b\Delta)}$$ and $C_1$ depends only on $r$. 
Applying therefore this inequality to each disc and
choosing
$d$ sufficiently small in order that
$|P_{\nu}|_{W_1^d}-P_{\nu}|_{W_1}|<\varepsilon$, we get
$|f^d-f|_{L^1(W_1)}$ arbitrarily small.

The proof of Proposition 1.14 is complete. $\square$

\smallskip
{\it Remark.} 
In 1.12, we had $M_r'(v)\subset {\cal W}$ being small
approach manifolds {\it not glued} to $M$ and we had that
$f|_{M_r'(v)}\to f|_{M'}$ in $L^1$ followed by a 
trivial subharmonic inequality,
whereas Carleson's estimate was needed in 1.14 for the smoothing
and gluing of $M_r'(v)$ to $M$ $\square$

\smallskip The following is a consequence of the minimalization
processus given by successive deformations of discs along the orbits
as in \cite{TU3} or \cite{JO3}.

\smallskip
 
{\sc Proposition 1.15.} \cite{PO1}.  {\it Let $M$ be generic, let
$p\in M$ and assume that $M$ is globally minimal at $p$. Then, for
each $\varepsilon >0$, each $0< \beta < \alpha$, there exists a
compactly supported ${\cal C}^{2,\beta}$ deformation $M^d$ of $M$ with
$M^d=M$ in a neighborhood of $p$ in $M$, $||M^d-M||_{{\cal
C}^{2,\beta}} < \varepsilon$, such that

\hspace{0.25cm} $1)$ There exists a ${\cal C}^{2,\beta}$ regular
family of analytic discs $A_{s,v}$ attached to $M^d$ at $p${\rm ;}

\hspace{0.25cm} $2)$ There exists a bounded linear extension operator
$$L_{loc,CR}^1(M) \to L_{loc,CR}^1(M^d), \ \ \ \ \ f \mapsto f^d, \
f^d=f \ on \ M\cap M^d{\rm ;}$$

\hspace{0.25cm} $3)$ Moreover, given $f\in L_{loc,CR}^1(M)$, there
exists such a deformation with $|f^d-f|_{L^1} < \varepsilon$.}

\smallskip Then Propositions 1.14 and 1.12 show that the regular
family $A_{s,v}$ given in 1.15 enables one to extend $f$ with $L^1$
control $2\varepsilon$ at $p$ into the wedge of edge $M$ at $p$
defined by $A_{s,v}$.  By means of the above result on the existence
of regular families of discs, and Proposition 1.14, we now can prove
the main result of Section 1.  We recall that Proposition 1.16
below will be suitable to apply the continuity principle and that, 
as explained in 1.7, it reduces the problem to considering $L^1{M}_{loc}\cap
{\cal H} ({\cal V}(M\backslash N))$ instead of simply 
$L^1_{loc}(M)\cap L^1_{loc,CR}(M\backslash N)$. The proof is technical and
consists in a great number of deformations
of $M$ into wedges defined by regular families of analytic discs attached
to $M$ and attached to its subsequent 
deformations and in applying the estimates
given by Propositions 1.14 and 1.15.

\smallskip
 
{\sc Proposition 1.16.} {\it Let $M$ be generic, ${\cal
C}^{2,\alpha}$-smooth, let $\Phi\subset M$ be a proper closed subset
of $M$, assume that $M\backslash \Phi$ is globally minimal, let $f\in
L_{loc,CR}^1(M\backslash \Phi)\cap L_{loc}^1(M)$ and let $U$ be an
open subset of $M\backslash \Phi$ whose closure is compact in $M$.
Then, for each $\varepsilon >0$, there exists a ${\cal C}^{2,\beta}$
deformation $d$, $M^d \supset \Phi$, of $M$ with $\hbox{supp} \ d
=\overline{U}$, $||M^d-M||_{{\cal C}^{2,\beta}} < \varepsilon$, such
that there exists a function $$f^d \in L_{loc}^1(M^d) \cap
L_{loc,CR}^1 (M^d \backslash \Phi )\cap {\cal H} ({\cal V}(U^d)),$$
coinciding with $f$ on $M\backslash (U\cup \Phi)$ and such that
$|f^d-f|_{L^1(U)} < \varepsilon.$}

\smallskip {\it Proof.} Consider an exhaustion of $U$ by compact sets
$$U=\bigcup_{\nu=0}^{\infty} K_{\nu}, \ \ \ \ \ K_{\nu} \subset
\stackrel{\circ}{K}_{\nu+1}.$$ Since global minimality is a stable
 property under smooth deformations vanishing sufficiently fast at
infinity, there exists a decreasing sequence $(\delta_{\nu})_{\nu\in
\N}$ of positive numbers, $0<\delta_{\nu+1} \leq \delta_{\nu}$, such
that for every ${\cal C}^{2,\beta}$ deformation $d$ of $M$ with
$||M^d-M||_{{\cal C}^{2,\beta}(K_{\nu} \backslash
\stackrel{\circ}{K}_{\nu-1})} < \delta_{\nu}$, $\nu=0,1,...$,
$K_{-1}=\emptyset$, $M^d \backslash \Phi$ is globally minimal.

We shall assume by induction on $\nu$ that, for each sequence
$0<\varepsilon_{\nu+1} \leq \varepsilon_{\nu} \leq \cdots \leq
\varepsilon_1 \leq \varepsilon_0 \leq \varepsilon$, there exists a
deformation $d_{\nu}$ of $M$, $M^{d_{\nu}}$, such that $K_{\nu}
\subset \hbox{supp} \ d_{\nu} \subset \stackrel{\circ}{K}_{\nu+1}$,
$||M^{d_{\nu}}-M||_{{\cal C}^{2,\beta} (K_j\backslash
\stackrel{\circ}{K}_{j-1})}< \varepsilon_j$, $j=0,...,\nu+1$ and there
exists a function $$f^{d_{\nu}} \in L_{loc,CR}^1(M^{d_{\nu}}) \cap
{\cal H}({\cal V}(K_{\nu}^{d_{\nu}}))$$ such that
$|f^{d_{\nu}}-f|_{L^1(U)} < \varepsilon$.

For $\nu=0$, this is a consequence of 1.14, if one chooses
$K_0=\overline{W}_2^0$ and $K_1=\overline{W}_1^0$, changing only the
first two compact sets in the exhaustion if necessary.

Fix $\varepsilon>0$ and $0< \varepsilon_{\nu+2} \leq
\varepsilon_{\nu+1} \leq \varepsilon_{\nu} \leq \cdots \leq
\varepsilon_1 \leq \varepsilon_0 \leq \varepsilon/2$ arbitrary, with
$\varepsilon_j \leq \delta_j$, $j=0,...,\nu+2$.  Let $d_{\nu}$,
$M^{d_{\nu}}$, $f^{d_{\nu}}$ be as in the induction hypothesis, with
$\varepsilon /2$ in place of $\varepsilon$ and
$$||M^{d_{\nu}}-M||_{{\cal C}^{2,\beta} (K_{\nu+1}\backslash
\stackrel{\circ}{K}_{\nu})} < \varepsilon_{\nu+1}/2, \ \ \ \ \
||M^{d_{\nu}}-M||_{{\cal C}^{2,\beta} (K_{\nu}\backslash
\stackrel{\circ}{K}_{\nu-1})} < \varepsilon_{\nu}/2,$$ and
$||M^{d_{\nu}}-M||_{{\cal C}^{2,\beta} (K_{j}\backslash
\stackrel{\circ}{K}_{j-1})} < \varepsilon_j$, for
$j=\nu-1,...,0$. Then $|f^{d_{\nu}}-f|_{L^1(U)} \leq \varepsilon /2$.

$M^{d_{\nu}}\backslash \Phi$ is globally minimal.  By compactness of
$K_{\nu+1}^{d_{\nu}} \backslash
(\stackrel{\circ}{K}_{\nu})^{d_{\nu}}$, there exists a finite number
$\mu=\mu(\nu)$ of points $p_1,...,p_{\mu} \in K_{\nu+1}^{d_{\nu}}
\backslash (\stackrel{\circ}{K}_{\nu})^{d_{\nu}}$ and there exist
regular families $A_{s_1,v_1}: {\cal S}^1 \times {\cal V}^1 \times
\overline{\Delta} \to \C^{m+n},...,A_{s_{\mu},v_{\mu}}:
 {\cal S}^{\mu} \times {\cal V}^{\mu} \times \overline{\Delta}\to
\C^{m+n}$ of analytic discs attached to suited deformations of
$M^{d_{\nu}}$ as in 1.15, with $A_{s_j,v_j}(1)=p_j$, such that
$A_{s_j,0}(e^{i\theta})$ embeds ${\cal S}_1^j \times I_1^j$ onto
$W_1^j=\{A_{s_j,0}(e^{i\theta})\in M^{d_{\nu}}; \ s_j\in {\cal S}_1^j,
\theta\in I_1^j\}$ and such that $$K_{\nu+2}^{d_{\nu}} \backslash
(\stackrel{\circ}{K}_{\nu-1})^{d_{\nu}}\supset W_3^1 \cup \cdots \cup
W_3^{\mu} \supset K_{\nu+1}^{d_{\nu}} \backslash
(\stackrel{\circ}{K}_{\nu})^{d_{\nu}},$$ for some $S_1^j \subset
\subset S^j$, $I_1^j \subset \subset I^j$,
$W_2^j=\{A_{s_j,0}(e^{i\theta})\in M^{d_{\nu}}; \ 2s_j\in {\cal
S}_1^j, 2\theta\in I_1^j\}$, $W_3^j=\{A_{s_j,0}(e^{i\theta})\in
M^{d_{\nu}}; \ 3s_j\in {\cal S}_1^j, 3\theta\in I_1^j\}$.  Set also
${\cal S}_2^j=\frac{1}{2} {\cal S}_1^j$, $I_2^j=\frac{1}{2} I_1^j$,
${\cal S}_3^j=\frac{1}{3} {\cal S}_1^j$, $I_3^j=\frac{1}{3} I_1^j$.

Moreover, we can assume that there exists $\gamma_{\nu}>0$ such that
for every ${\cal C}^{2,\beta}$ deformation $(M^{d_{\nu}})^d$ of
$M^{d_{\nu}}$ with $\hbox{supp} \ d \subset K_{\nu+2}^{d_{\nu}}
\backslash K_{\nu-1}^{d_{\nu}}$, $||(M^{d_{\nu}})^d-
M^{d_{\nu}}||_{{\cal C}^{2,\beta}(K_{\nu+2}^{d_{\nu}} \backslash
(\stackrel{\circ}{K}_{\nu-1})^{d_{\nu}})} < \gamma_{\nu}$, and the
support of $d$ not meeting the part where $A_{s_j,v_j}^d$ is not
attached to $M^{d_{\nu}}$ (as in 1.15), then the deformed discs
$A_{s_j,v_j}^d$ exist and give regular families.  Indeed, Bishop's
equation allows small perturbations.

Let $\varphi_3^1$ be a smooth function on ${\cal S}_1^1 \times I_1^1$
such that $\hbox{supp} (1-\varphi_3^1) \subset \subset {\cal S}_1^2
\times I_1^2$, $\varphi_3^1 \equiv r_1$ near ${\cal S}_3^1 \times
I_3^1$, with $r_1 < 1$ very close to $1$ and $\varphi_3^1$ very close
to the constant $r_1$ in ${\cal C}^{\infty}$ norm. We set a
deformation of $M^{d_{\nu}}$ to be $(M^{d_{\nu}})^{d_{\nu}^1}=
M^{d_{\nu}}$ outside $W_1^1$ and $$(M^{d_{\nu}})^{d_{\nu}^1}=
\{A_{s_1,0}(\varphi_3^1(s_1,\theta_1)e^{i\theta_1})\ \! {\bf :} \ \!
s_1\in {\cal S}_1^1, \theta_1 \in I_1^1\},$$ We can assume
$\varepsilon_{\nu+2} \leq \gamma_{\nu}$.  Proposition 1.14 with
$\varepsilon=\varepsilon_{\nu+2}/4$, implies that we can choose $r_1$
and $\varphi_3^1$ close to $1$ in order that
$$||(M^{d_{\nu}})^{d_{\nu}^1}-M||_{{\cal C}^{2,\beta}
(K_{\nu+1}\backslash \stackrel{\circ}{K}_{\nu})} \leq
\varepsilon_{\nu+1}/2+ \varepsilon_{\nu+2}/4\leq
\varepsilon_{\nu+1}/2+ \varepsilon_{\nu+1}/4,$$
$$||(M^{d_{\nu}})^{d_{\nu}^1}-M||_{{\cal C}^{2,\beta}
(K_{\nu}\backslash \stackrel{\circ}{K}_{\nu-1})} \leq
\varepsilon_{\nu}/2+ \varepsilon_{\nu+2}/4\leq \varepsilon_{\nu}/2+
\varepsilon_{\nu}/4,$$ and such that there exists
$(f^{d_{\nu}})^{d_{\nu}^1} \in
L_{loc,CR}^1((M^{d_{\nu}})^{d_{\nu}^1})$, such that
$$|f-(f^{d_{\nu}})^{d_{\nu}^1}|_{L^1(U)} < |f-f^{d_{\nu}}|_{L^1(U)}+
|f^{d_{\nu}}-(f^{d_{\nu}})^{d_{\nu}^1}|_{L^1(U)}< \varepsilon /2 +
\varepsilon_{\nu+2} /4 < \varepsilon /2+ \varepsilon/4.$$ To make a
second (and decisive) step, let $A_{s_2,v_2}$ and $W_1^2$, $W_2^2$,
$W_3^2$ be as above. We denote by $\psi_1: (s_1,\theta_1) \mapsto
A_{s_1,0} (e^{i\theta_1})$, $\psi_2: (s_2,\theta_2) \mapsto
A_{s_2,0}(e^{i\theta_2})$.  Choose a smooth function $\chi$ on ${\cal
S}_1^1 \times I_1^1$ such that $\hbox{supp} (1-\chi)\subset {\cal
S}_1^1 \times I_1^1$, $\chi\equiv \varphi_3^1$ near $\psi_1^{-1}(W_2^1
\cap W_2^2)$ and $\chi \equiv 1$ outside $\psi_1^{-1}(W_1^1\cap
W_1^2)$. Define $$(M^{d_{\nu}})^{d_{\nu}^{1/2}}= \{A_{s_1,0}(\chi
(s_1,\theta_1) e^{i\theta_1})\ \! {\bf :} \ \!  (s_1, \theta_1) \in
{\cal S}_1^1 \times I_1^1\},$$ and
$(M^{d_{\nu}})^{d_{\nu}^{1/2}}=M^{d_{\nu}}$ outside $W_1^1$. By
Proposition 1.14, we can choose $\chi$ in order that
$$||(M^{d_{\nu}})^{d_{\nu}^{1/2}}-M||_{{\cal C}^{2,\beta}
(K_{\nu+1}\backslash \stackrel{\circ}{K}_{\nu})} \leq
\varepsilon_{\nu+1}/2+ \varepsilon_{\nu+2}/4\leq
\varepsilon_{\nu+1}/2+ \varepsilon_{\nu+1}/4,$$
$$||(M^{d_{\nu}})^{d_{\nu}^{1/2}}-M||_{{\cal C}^{2,\beta}
(K_{\nu}\backslash \stackrel{\circ}{K}_{\nu-1})} \leq
\varepsilon_{\nu}/2+ \varepsilon_{\nu+2}/4\leq \varepsilon_{\nu}/2+
\varepsilon_{\nu}/4,$$ and such that there exists
$(f^{d_{\nu}})^{d_{\nu}^{1/2}} \in
L_{loc,CR}^1((M^{d_{\nu}})^{d_{\nu}^{1/2}})$ such that
$$|f-(f^{d_{\nu}})^{d_{\nu}^{1/2}}|_{L^1(U)} \leq \varepsilon/2 +
\varepsilon /4.$$ Now, we consider the deformed regular family of
discs $(A_{s_2,v_2})^{d_{\nu}^{1/2}}$ attached to
$(M^{d_{\nu}})^{d_{\nu}^{1/2}}$. The crucial fact is that
$(M^{d_{\nu}})^{d_{\nu}^{1/2}}$ is not deformed where the discs
$A_{s_2,v_2}$ is not attached to $M^{d_{\nu}}$, as stated in 1.15.
Moreover, since $\varepsilon_{\nu+2}/4 \leq \gamma_{\nu}$,
$(A_{s_2,v_2})^{d_{\nu}^{1/2}}$ exists and gives a regular family.
Choose a smooth function $\varphi_3^2$ such that $\hbox{supp}
(1-\varphi_3^2) \subset {\cal S}_2^2 \times I_2^2$, $\varphi_3^2\equiv
r_2 <1$ near ${\cal S}_3^2\times I_3^2$, and set
$$(M^{d_{\nu}})^{d_{\nu}^2}= \{(A_{s_2,v_2})^{d_{\nu}^{1/2}}
( \varphi_3^2(s_2,\theta_2)e^{i\theta_2})\ \! {\bf :} \ \!
(s_2,\theta_2)\in {\cal S}_2^2 \times I_2^2\},$$ and
$(M^{d_{\nu}})^{d_{\nu}^2}=M^{d_1}$ outside
$(W_2^2)^{d_{\nu}^{1/2}}$. If $\varphi_3^2$ is sufficiently close to
1, we have $$||(M^{d_{\nu}})^{d_{\nu}^{2}}-
(M^{d_{\nu}})^{d_{\nu}^{1}}||_{{\cal C}^{2,\beta} (K_{\nu+1}\backslash
\stackrel{\circ}{K}_{\nu})} \leq \varepsilon_{\nu+2}/8 \leq
\varepsilon_{\nu+1}/8,$$ $$||(M^{d_{\nu}})^{d_{\nu}^{2}}-
(M^{d_{\nu}})^{d_{\nu}^{1}}||_{{\cal C}^{2,\beta} (K_{\nu}\backslash
\stackrel{\circ}{K}_{\nu-1})} \leq \varepsilon_{\nu+2}/8 \leq
\varepsilon_{\nu}/8,$$ and there exists $(f^{d_{\nu}})^{d_{\nu}^{2}}
\in L_{loc,CR}^1((M^{d_{\nu}})^{d_{\nu}^{2}})$ such that
$$|(f^{d_{\nu}})^{d_{\nu}^{2}}- (f^{d_{\nu}})^{d_{\nu}^{1}}|_{L^1(U)}
\leq \varepsilon/8.$$ It suffices to repeat the deformations with
$j=3,...,\mu$.

The proof of Proposition 1.16 is complete. $\square$

We end this section with an elementary proof of the following.

\smallskip
 
{\sc Proposition 1.17.} {\it Let $M$ be generic, ${\cal C}^2$-smooth,
in $\C^{m+n}$, ${\rm dim}_{CR} M = m$, let $N \subset M$ be a generic
two-codimensional ${\cal C}^2$-smooth submanifold.  If $m\geq 2$, $N$
is $L^1$-removable.}

\smallskip {\it Proof.} Let $p\in N$, let $U\ni p$ be a neighborhood
of $p$ in $M$ and choose two one-codimensional ${\cal C}^2$-smooth
submanifolds $M_1,M_2$ of $M$ containing $N$ near $p$ such that
$T_pM_1+T_pM_2=T_pM$. Let $\Upsilon_j$, $j=1,2$, denote the set of 
${\cal C}^1$ sections $L_j\in
\Gamma(T^{0,1} U)$ such that $L_j|_{M_j} \in \Gamma(T^{0,1}M_j)$,
$j=1,2$. Any $L\in \Gamma(T^{0,1}U)$ is a linear combination of a
section $L_1$ of $\Upsilon_1$ and a section $L_2$ of $\Upsilon_2$.
For this, we use $m\geq 2$.
 Therefore, it
suffices to show that, given $f\in L_{CR}^1(U\backslash N)$, then
$L_j(f)=0$ in the weak sense, for any $L_j\in \Upsilon_j$, $j=1,2$.  Let
$\varphi\in {\cal C}^{\infty}_c(U)$, let $\chi_{j,l}\in {\cal
C}_c^{\infty}(U)$, $l\in \N$, $|\nabla\chi_{j,l}|\leq Cl$, with
$\chi_{j,l}\equiv 1$ in a neighborhood of $M_j\cap \hbox{supp} \
\varphi$, and $\lim \chi_{j,l}=0$ in $L^1(U)$, as $l\to \infty$.  Then
$$\int_U {\:}^{\tau}L_j(\varphi) f= \int_{U\backslash M_j} {\:}^{\tau}
L_j((1-\chi_{j,l}) \varphi)f+ \int_{{\cal V}(M_j\cap U)}
{\:}^{\tau}L_j(\chi_{j,l}\varphi) f= \int_{{\cal V}(M_j\cap U)}
{\:}^{\tau}L_j(\chi_{j,l}\varphi) f,$$ since $f$ is CR outside $M_j$,
where the transposition is relative to a fixed measure on $U$. This
last term tends to zero as $l$ tends to infinity, since $L_j|_{M_j}$
is tangent to $M_j$, exactly as in the proof of Theorem 1.1.
For this we need $M_1, M_2$ ${\cal C}^2$.

The proof of Proposition 1.17 is complete. $\square$

\smallskip

{\bf 2. Proofs of Theorem 4 (i), (ii) and (iii).}  Now, we come to the the
general case.  We plan to show first that these three theorems reduce
to a single statement, namely Proposition 2.2 below. This is done as
follows. Recall that a point $p\in \Phi$ is called {\it removable
in} $M\backslash \Phi$ ({\it removable}, for short)
if there exists a small neighborhood
$V$ of $p$ in $M$ such that $L_{CR}^1$ or ${\cal W}$ functions
over $M\backslash \Phi$ extend in the $L^1_{CR}$ or ${\cal W}$ sense over
$V$.

{\sc Lemma 2.1.} {\it  
(i) and (iii) in Theorem 4 reduce to (ii).}

\smallskip
{\it Proof.}
We have $m\geq 2$.  Since
$M\backslash N$ is globally minimal, Proposition 1.16 with 1.7 show
that it suffices to prove that $L_{loc,CR}^1(M\backslash N) \cap {\cal
H} ({\cal V}(M\backslash N))\subset L_{loc,CR}^1(M)$.  Let $p\in M$,
${\rm dim}_{CR} T_p N =m-2$.  By Proposition 1.17, $p$ is
$L^1$-removable.  An alternative argument is as follows.  By the
remark after Lemma 2.10 or by Theorem 5.A.1 in \cite{ME2}, $p$ is
${\cal W}$-removable.  By Proposition 2.11, this implies that $p$ is
$L^1$-removable.  This proves that it suffices to show that the proper
closed subset $\Phi=\{p\in N\ \! {\bf :} \ \! {\rm dim}_{CR} T_p N
\geq m-1\}$ is $L^1$-removable. Clearly also,
Theorem 4 (i)
is implied by (ii).

The proof that (i) and (iii) reduce to (ii) in Theorem 4 is complete. $\square$

\smallskip
To prove Theorem 4 (ii), proceed as follows. Set $${\cal A}=
\{\Psi \subset \Phi \ \hbox{closed} \ \! {\bf :} \ \!  M\backslash
\Psi \ \hbox{is globally minimal} \ \ \ \ \ \ \ \ \ $$
$$\ \ \ \ \ \ \ \ \ \ \ \ \ \ 
\ \ \ \ \ \ \ \ \ \ \hbox{and} \ L_{loc,CR}^1(M\backslash \Phi)
\cap L_{loc}^1(M)=
L^1_{loc,CR}(M\backslash \Psi)\cap L_{loc}^1(M)\},$$ and 
$${\cal B}=\{\Psi \subset
\Phi \ \hbox{closed} \ ; M\backslash \Psi \ \hbox{globally minimal
and} \ \Phi\backslash \Psi \ {\rm is} \ {\cal W}{\rm -removable} \
\hbox{in} \ M\backslash \Phi \}.$$ As
${\cal A}$ ({\em resp.} ${\cal B}$) is closed under arbitrary
intersections (obvious), it contains especially $$\Psi_{{\rm
nr}}=\bigcap_{\Psi \in {\cal A}} \Psi \ \in {\cal A}, \ \ \ \ \
(resp. \ \ \Psi_{{\rm nr}}=\bigcap_{\Psi \in {\cal B}} \Psi \ \in
{\cal B}).$$ 
Assume that $\Psi_{{\rm nr}} \neq \emptyset$.  We shall
reach to a contradiction by proving that there is
some point $p_1$ in $\Psi_{{\rm
nr}}$ such that 

1. There exists a section $L\in \Gamma(T^cM)$ near $p_1$
with $L_t(p_1) \in M\backslash \Psi_{{\rm nr}}$, $0< t < \delta$, 
for some
$\delta >0$ and

2. The point $p_1$ is $L^1$- ({\em resp.} ${\cal W}$-) removable
in $M\backslash \Psi_{\rm nr}$.

\noindent
Then by 1., if $V$ denotes a small neighborhood
of $p_1$ in $M$ which is removable as in 2., we have that
$(M\backslash \Psi_{\rm nr})\cup V$ is 
globally minimal and that $\Phi \backslash (V \cup
(M\backslash \Psi_{\rm nr}))$ is removable  in
$M\backslash \Phi$, contradicting the
definition of ${\cal A}$ or ${\cal B}$.
 
Denote from now on $\Psi_{{\rm nr}}$ again by $\Phi\neq\emptyset$.

Since $M\backslash \Phi$ is globally minimal,
$L_{loc,CR}^1(M\backslash \Phi)$ extends holomorphically into a wedge
at every point of $M\backslash \Phi$.  Notice that we have
derived in Proposition 1.16 (see also 1.7) that we can slightly deform $M$ over
$M\backslash \Phi$ in a manifold $M^d$ in order that we are given the
space $L_{loc,CR}^1(M^d\backslash \Phi) \cap {\cal H} ({\cal
V}(M^d\backslash \Phi))$ for $L^1$-removability, with $L^1$-control, and
${\cal H} ({\cal V}(M^d\backslash \Phi))$ for ${\cal W}$-removability,
which is a crucial reduction to apply the continuity principle.
Since we can let $d$ tend to $0$ and since the construction of any
wedge ${\cal W}^d$ attached to $M^d$ by means of discs will depend
smoothly on $d$ (see also Section 5 in \cite{ME2}), it is sufficient
to prove the following statement, a single proposition which implies
therefore ${\cal W}$-removability in Theorem 4 (ii).  
$L^1$-removability in Theorem 4 (ii) will be derived thereafter in
Proposition 2.11 from Proposition 2.2. We rename
such $M^d$ as $M$.

\smallskip
  
{\sc Proposition 2.2.}  {\it Let $M$ be generic in $\C^{m+n}$, ${\cal
C}^{2,\alpha}$, ${\rm dim}_{CR} M= m \geq 1$, let $N \subset M$ be a
connected ${\cal C}^2$-smooth submanifold with $\hbox{codim}_M N =2$,
let $\Phi \subset N$ be a proper closed subset, $\Phi \neq \emptyset$,
assume that $M$ and $M\backslash \Phi$ are globally minimal and let
$\omega$ be a neighborhood of $M\backslash \Phi$ in $\C^{m+n}$.  Then
there exists a point $p_1\in b\Phi$ such that there exists a section
$L\in \Gamma(T^cM)$ near $p_1$ with $L_t(p_1) \in M\backslash \Phi$,
$0< t < \delta$, for some $\delta >0$, and there exists a wedge ${\cal
W}_{p_1}$ at $p_1$ such that, for every function $f\in {\cal
H}(\omega)$, there exists a function $F\in {\cal H}({\cal W}_{p_1})$
 with $F=f$ in the intersection of ${\cal W}_{p_1}$ with a
neighborhood of $M\backslash \Phi$ in $\C^{m+n}$.}

\smallskip {\it Proof.}  The proof will be given in four steps.  Let
us begin with the following. We consider the interior and the boundary 
of $\Phi$ as a subset of the manifold $N$.

\smallskip
 
{\sc Lemma 2.3.} {\it There exists a point $p_1\in b\Phi$ such that
there exists a section $L\in \Gamma(T^cM)$ near $p_1$ with $L_t(p_1)
\in M\backslash \Phi$, $0< t < \delta$, some $\delta >0$, and such
that there exists a generic one codimensional ${\cal C}^2$-smooth
manifold $M_1$ through $p_1$ with $\Phi$ contained in one closed 
side $M_1^{-}$ of $M_1$ at $p_1$, {\it i.e.} a side 
$V^-$ of a small neighborhood $V$ of $p_1$ in $M$ divided
by $M_1$ in two closed parts $V^+$ and $V^-$ with 
$V^+\cap V^-=V\cap M$.}

\smallskip {\it Proof.}  Assume first that there exists a point $p_0
\in b\Phi$ such that $T_{p_0} N \not\supset T_{p_0}^cM$. 
Then we easily construct a generic hypersurface $M_1\subset M$ 
touching $N$ exactly in $p_0$ and therefore containing 
$N\supset\Phi$ in one of its sides.
Existence of $L$ with $L_t(p_1)\in M\backslash \Phi$, 
$0< t < \delta$ is trivial.

Assume now that $T_q N \supset T_q^cM$, for each point $q$ in the
closed subset $b\Phi\subset N$. By the 
{\it local CR orbit} ${\cal O}_{CR}^{loc}(p)$ of a
point $p\in M$, we mean a representative of
$\lim_{U\ni p} {\cal O}_{CR}(U,p)$ (\cite{ME1}).
Recall 
$N$ is ${\cal C}^2$-smooth.  Choose a point
$p_0\in b\Phi$ such that ${\cal O}^{loc}_{CR}(p_0)\not\subset b\Phi$
and let $U$ be a neighborhood of $p_0$ in $M$.  (Otherwise, if ${\cal
O}^{loc}_{CR}(q)\subset b\Phi$, $\forall \ q\in b\Phi$, $b\Phi$
contains a CR orbit, which contradicts the assumption 
that $M$ and $M\backslash N$ are globally minimal). We have $T_{p_0} N \supset
T_{p_0}^cM$. Let $S$ be a two-dimensional submanifold of $M$ through
$p_0$ with $T_{p_0} S + T_{p_0} N = T_{p_0} M$.  Let $\pi: M \to N$ be
a submersion parallel to $S$ in smooth (linear) coordinates on $M$
such that $S$ and $N$ are coordinate spaces.  For any vector field $L
\in \Gamma(U, T^cM)$, one can define a vector field $L_N$ on $N\cap U$
by taking $L_N(q)=\pi_*(L(q))$, $q\in N\cap U$. This defines a vector
bundle $K\subset TN$ of rank $2m$ (this holds since $T_{p_0} N \supset
T_{p_0}^cM$).  Assume that ${\cal O}^{loc}_K(p_0)\subset b \Phi$.
Since, by assumption, $K(q)=T_q^cM$, for each $q\in b\Phi$, this
implies that ${\cal O}^{loc}_{CR}(p_0)\subset b\Phi$, which is not
true. Therefore, ${\cal O}^{loc}_K(p_0)\not\subset b\Phi$ near $p_0$.

Following Bony \cite{BONY}, we shall say that a vector field $L_N\in
\Gamma(K)$ is tangent to the closed set $b\Phi$ near $p_0$ if, for any
point $q\in U\cap N$, any open ball $B$ with center $q$ such that $B
\subset N\backslash b\Phi$, then, for every point $p_1 \in bB \cap
b\Phi$, $L_N(p_1)$ is tangent to $bB$ at $p_1$.  We shall use the
following theorem of Bony \cite{BONY}. {\it If a Lipschitz real vector
field $X$ on $\R^{k}$ is tangent in the above sense to a closed subset
$F\subset \R^k$, then every integral curve meeting $F$ stays in $F$}.
By this theorem, the condition that ${\cal O}^{loc}_K(p_0)\not \subset
b\Phi$ implies that there exists an open ball $B\subset N\backslash
b\Phi$ such that $bB \cap b\Phi=\{p_1\}$ (is a single point, which
holds true after a homothety with center $p_1$) and there exists a
section $L\in \Gamma(U,K)$ such that $L(p_1)\not \in T_{p_1}
bB$. Choose ${\cal C}^1$ coordinates $s_1,...,s_{2m+n}$ on $M$ such
that lines $(s_2,...,s_{2m+n})=const$ correspond to integral curves of
$L$. Since $p_1\in b\Phi$, we can choose a point $q\in
N\backslash\Phi$ close to $p_1$. The set $\{L_{N,t}(q)\ \! {\bf :} \
\! t\in I\}$, $I$ an open interval with origin $0$ in $\R$, uniform in
$q$, {\em i.e.} the integral curve of $L_N$ with origin $q$, meets
$\Phi$ in $U\cap N$ or does not meet $\Phi$.  If it does meet $\Phi$,
it meets $b\Phi$ for a smallest $|t|$, say $t_0<0$ (changing $L_N$ in
$-L_N$ if necessary), so we get the existence of a new
$p_0=L_{N,t_0}(q) \in b\Phi$ such that $L_{N,t}(p_0) \in N\backslash
\Phi$, $0<t<\delta$, some $\delta>0$. Since $\pi$ is parallel to $S$
and since $T_{p_0}^cM\subset T_{p_0}N$, this implies that $L_t(p_0)\in
M\backslash \Phi$, $0< t <\delta$. Assume on the contrary that, for
each $q\in N\backslash \Phi$, the set $\{L_{N,t}(q) \ \! {\bf :} \ \!
t\in I\}$ does not meet $\Phi$ in $U\cap N$. As a consequence,
$\{L_{N,t}(p_1)\ \! {\bf :} \ \! t\in I\}\subset \Phi$. But since
$p_1\in b\Phi$, there exist points $q\in N \backslash \Phi$
arbitrarily close to $p_1$. Therefore, in fact, $\{L_{N,t}(p_1)\ \!
{\bf :} \ \! t\in I\} \subset b\Phi$, which is contrary to the fact
that $B\subset N\backslash b\Phi$, $p_1\in b \Phi$ and
$L_N(p_1)\not\in T_{p_1} bB$.

Fix a point $p_0\in b \Phi$ as above with $L_{N,t}(p_0)\in N\backslash
\Phi$, $0<t<\delta$. Let $q=L_{N,t}(p_0)\in N\backslash \Phi$, some
$t>0$.  We can assume that $s_1,s_2,...,s_{2m+n-2}$ give coordinates
which are ${\cal C}^1$ on $N$ (via a graphing function for $N$) and that
the lines $s_2,...,s_{2m+n-2}=const$ correspond
to integral curves of $L_N$.  Consider, for $\varepsilon\geq 1$, the
increasing family of domains with ${\cal C}^1$-smooth boundary
contained in $N$ $$D_N(q,\varepsilon)=\{(s_1-s_1(q))^2/ \varepsilon^2+
s_2^2+\cdots+s_{2m+n-2}^2<r^2\},$$ where $r>0$ is small enough to
 insure that $D_N(q,1)\subset \subset N\backslash \Phi$.  Then, for
each $\varepsilon \geq 1$, a point $r\in b D_N(q,\varepsilon)$
satisfies $L_N(r)$ tangent to $ D_N(q,\varepsilon)$ if and only if
$s_1(r)=s_1(q)$, and this implies that $r\in N\backslash \Phi$.
Therefore, there exists the unique least $\varepsilon_0>1$ such that
$b D_N(q, \varepsilon_0)\cap \Phi =b D_N(q, \varepsilon_0)\cap
b\Phi\neq 0$ and for each $p_1$ (new $p_1$) in the intersection, one
has $L(p_1)$ nontangential to $b D_N(q,\varepsilon_0)$.  For ${\cal
C}^2$ coordinates close in ${\cal C}^1$ norm to the coordinates
$(s_1,...,s_{2m+n})$, in this new geometric situation, we can repeat
the argument of increasing balls and get the same result with $b
D_N(q, \varepsilon_0)$ ${\cal C}^2$.  Now a part of the ${\cal C}^2$-smooth
boundary of the ball $bD_N(q, \varepsilon_0)$ can be included in a
${\cal C}^2$-smooth manifold $M_1\subset M$ of codimension one in $M$
which is generic in $\C^{m+n}$ and such that there exists $p_1\in
b\Phi\cap M_1$ and $\Phi \subset M_1^-$.

The proof of Lemma 2.3 is complete. $\square$

\smallskip {\it Step one: existence of a disc.}  Choose holomorphic
coordinates $(w,z)$ such that $p_1=0$, $M_1$ is a germ through $p_1$
of a generic one codimensional ${\cal C}^2$-smooth submanifold of $M$
with $T_0 M_1 =\{x=0, u_1=0\}$, $M=\{x=h(y,w)\}$ and $M_1$ is given in
$M$ by the supplementary equation $u_1=k(v_1,w_2,...,w_m,y)$, 
$w_1=u_1+iv_1$, for a
${\cal C}^2$-smooth function $k$ with $k(0)=0$ and $dk(0)=0$. We note
$M_1^-=\{u_1\leq k(v_1,w_2,...,w_m,y)\}$.

\smallskip
 
{\sc Lemma 2.4.} {\it There exists an embedded analytic disc $A\in
C^{2,\beta}(\overline{\Delta})$ with $A(1)=p_1$, $A(b\Delta)\backslash
\{1\}\subset M\backslash M_1^{-}$ and $\frac{d}{d\theta}|_{\theta=0}
A(e^{i\theta})= v_0\in T_{p_1}M_1$.}

\smallskip {\it Proof.}  For small $c>0$, take
$W_{c}(\zeta)=(c(1-\zeta),0,...,0)$ and consider the analytic disc
$A_c(\zeta)=(Z_c(\zeta),W_c(\zeta))$, where $Z_c$ is the ${\cal C}^{2,
\beta}$ solution of Bishop's equation $Y_c=T_1h(Y_c, W_c)$ on
$b\Delta$. Then, for $\zeta\in b\Delta$, $$U_1(\zeta)=\frac{c}{2}
|1-\zeta|^2, \ \ \ \ \ |Z_c(\zeta)|\leq O(c|1-
\zeta|(c^{\beta}+|1-\zeta|^{\beta})),$$ and there exists $c_0>0$ and a
constant $C>0$, depending on the second derivatives of $k$ in a ball
containing $0$ in $\R\times\C^{m-1}\times\R^n$, such that, for each
$c\leq c_0$, $|k\circ A(\zeta)| < C(c^2|1-\zeta|^2)$.  Recall
$M_1^+=\{u_1\geq k(v_1,w_2,...,w_m,y)\}$.  Choose now $c$ with
$cC<\frac{1}{2}$.

The proof of Lemma 2.4 is complete. $\square$

\smallskip Therefore, ${\cal W}$-removability of $p_1$ 
in $M\backslash \Phi$ is a consequence of the
following.

\smallskip
 
{\sc Proposition 2.5.} {\it Let $M$ be generic, ${\cal C}^{2,\alpha}$-smooth,
let $p_1\in M$, let $N\ni p_1$ be a ${\cal C}^2$ submanifold with $codim_M N
=2$, let $\Phi \subset N$ be closed, let $p_1 \in b\Phi$ and assume
that there exists a one codimensional generic ${\cal C}^2$-smooth
manifold $M_1 \subset M$ such that $\Phi \subset M_1^-$ and let
$\omega$ be a neighborhood of $M \backslash N$ in $\C^{m+n}$.  Let, as
in 2.4, a sufficiently small embedded analytic disc $A\in
C^{2,\beta}(\overline{\Delta})$ be attached to $M$, $A(1)=p_1$,
$\frac{d}{d\theta}|_{\theta=0} A(e^{i\theta})= v_0\in T_{p_1}M_1$,
with $A(b \Delta \backslash \{1\}) \subset M \backslash M_1^{-}$.
Then for each $\varepsilon >0$, there exist $v_{00}\in T_{p_1}M_1$
with $|v_{00}-v_0|<\varepsilon$, $v_{00}\not\in T_{p_1}N$,
$v_{00}\not\in T_{p_1}^cM$ and a wedge ${\cal W}$ of edge $M$ at
$(p_1, Jv_{00})$ such that for every holomorphic function $f \in {\cal
H}(\omega)$ there exists a function $F \in {\cal H}({\cal W})$ with
$F=f$ in the intersection of ${\cal W}$ with a neighborhood of
$M\backslash N$ in $\C^{m+n}$.}

\smallskip {\it Proof.}  Fix a function $f \in {\cal H}(\omega)$.  The goal
will be to construct deformations of our given original disc $A$ as in
\cite{ME2} with boundaries in $M \cup \omega$ to show that the
envelope of holomorphy of $\omega$ contains a (very thin) wedge of
edge $M$ at $A(1)$ (instead of appealing to a Baouendi-Treves
approximation theorem, a version of which {\it by no means} being 
{\it a priori} valuable here, which is explained
in the remarks before step five below), the natural tool being the
so-called {\it continuity principle} 1.6.

We can assume that $A(1)=0$ and that $M$ is given in a coordinate
system as in Lemma 1.12.  Set $A(\zeta)=(X(\zeta)+ i Y(\zeta),
W(\zeta))$.  Since $v_0\in T_0M_1$,
$\frac{d}{d\theta}|_{\theta=0}w_1(e^{i\theta})$ is purely imaginary.

The proof of Proposition 2.5 will be divided in four more
steps. During the second one, we will introduce a large family of
normal deformations of analytic discs, and during the third one, we
will check and use the isotopy properties of this family.
This is the most important step in our approch.

\smallskip {\it Step two: normal deformations.}  The following result
shows that any disc can be included in a regular family.

\smallskip
 
{\sc Proposition 2.6.} {\it Let $M$ be generic, ${\cal
C}^{2,\alpha}$-smooth, let $A\in {\cal
C}^{2,\beta}(\overline{\Delta})$ be a sufficiently small analytic disc
attached to $M$, $ A(1)=p_0$, with $v_0=\frac{d}{d\theta}|_{\theta=0}
A(e^{i\theta})\neq 0$ and let $\omega$ be a neighborhood of $A(-1)$ in
$\C^{m+n}$. Then there exists a ${\cal C}^{2,\beta}$-smooth family of
analytic discs attached to $M\cup \omega$, $A_{t,\tau,a,p}(\zeta)$,
with $t$ in a neighborhood ${\cal T}$ of $0$ in $\R^n$, $\tau \in
I_{\tau_0}=(-\tau_0,\tau_0)$, $\tau_0>0$, $a\in {\cal A}\subset
\C^{m-1}$, ${\cal A}$ a neighborhood of $0$ in $\C^{m-1}$, with $p$ in
a neighborhood ${\cal M}$ of $p_0$ in $M$, $A_{t,\tau,a,p}(1)=p$, such
that the rank of the mapping $t\mapsto -\partial A_{t,0,0,p_0}/
\partial \zeta (1)$ is equal to $n$ and such that the set
${\Gamma}_0=\{s \frac{d}{d\theta}|_{\theta=0}
A_{t,\tau,a,p_0}(e^{i\theta})\ \! {\bf :} \ \! s>0, t\in {\cal T},
\tau \in I_{\tau_0}, a\in {\cal A}\}$ is a $(2m+n)$-dimensional cone
with vertex $p_0$ in $T_{p_0}M$.}

\smallskip {\it Proof.} We include a proof of this result which is
crucial for the proof of Theorem 1. 
The main argument Lemma 2.7 below will be recalled here for
completeness (\cite{TU3}, \cite{ME2}): it is this availability
of variation of the outer direction of discs as to describe
a {\it cone} in $T_p\C^{m+n}/T_pM$ which underlies our propagation
of removability process.

We can assume that $A(1)=0$ in a
coordinate system as (1), and that the projection of $v_0=
\frac{d}{d\theta}|_{\theta=0} A(e^{i\theta})$ on the $v_1$-axis is non
zero.

Let $\mu=\mu(y,w)$ be a ${\cal C}^{\infty },$ $\R$-valued function with
support near the point $(y(-1), w(-1))$ that equals 1 there and let
$\kappa: \R^{n} \to \R^{n}$ be a ${\cal C}^{\infty }$ function with
$\kappa(0)=0$ and $\kappa'(0)=Id$.  We can assume that the supports of
$\mu$ and $\kappa$ are sufficiently concentrated in order that every
manifold $M_t$ with equation
              $$
              x=H(y,w,t)=h(y,w)+\kappa(t)\mu(y,w)
              $$ is contained in $\omega$ and the
deformation is localized in a neighborhood of $A(-1)$ in $\C^{m+n}$.
Let $\chi= \chi(\zeta)$ be a smooth function on the unit circle
supported
 in a small neighborhood of $\zeta=-1$ that will be chosen later.  Set
$A(\zeta)=(Z(\zeta),W(\zeta))$.  For every small $t$, $t\in {\cal T}$,
every small $\tau$, $\tau \in I_{\tau_0}$, every $a\in {\cal A}$,
every $p\in {\cal M}$, after denoting by $(y^0,w_1^0,...,w^0_m)$ the
coordinates of $p$ on $M$, we consider the disc
$$A_{t,\tau,a,p}(\zeta)=
(X_{t,\tau,a,p}(\zeta)+iY_{t,\tau,a,p}(\zeta),
e^{i\tau}W_1(\zeta)+w_1^0, W^*(\zeta)+a(\zeta-1)+w^{*0}),$$ where
$Y_{t,\tau,a,p}$ is the solution of Bishop's equation with parameters
$$Y_{t,\tau,a,p}(\zeta)=T_1H(Y_{t,\tau,a,p},e^{i\tau}W_1+w_1^0,
W^*+a(.-1)+w^{*0},t\chi)+y^0,$$ which exists and depends in a ${\cal
C}^{2,\beta}$-smooth fashion on $(t,\tau,a,p,\zeta)$. Then
$A_{t,\tau,a,p}(1)=p$. When $\tau=0$, $a=0$ and $p=p_0=0$, simply
denote $A_{t,0,0,0}$ by $A_t$.  We prove the crucial lemma 2.7 below,
due to Tumanov, on normal deformations of discs: by pushing $A$ into
$\omega$ near $A(-1)$ along every direction given by $t$ in the normal
bundle to $M$ at $A(-1)$, the inner tangential direction
$-\frac{\partial A_t}{\partial \zeta} (1)$ will describe a {\it whole}
open cone in the normal bundle to $M$ at $A(1)$.

let $\Pi$ denote the canonical bundle epimorphism $\Pi: T\C^{m+n} |_M
\to T\C^{m+n} |_M / TM$ and consider the ${\cal C}^{1,\beta}$ mapping
                    \begin{equation}
              D: \ \ \ \ \ \R^{n} \ni t \ \longmapsto \ \Pi \left( -
\frac{\partial A_t}{\partial \zeta} (1) \right) \in T_0 \C^{m+n} / T_0
M \simeq \R^{n}.
                    \end{equation}

\smallskip
 
{\sc Lemma 2.7.} ({\sc Tumanov} \cite{TU3}). {\it $\chi$ can be chosen
in order that $ \ rk \ D'(0)= n$.}

\smallskip {\it Proof.} Set $t= (t_1,...,t_n)$.  Differentiating the
equation $X_t(\zeta) = H(Y_t(\zeta),W(\zeta), t\chi(\zeta))$, $\zeta
\in b \Delta$, with respect to $t_j, j=1,...,n,$ we obtain that the
holomorphic disc $\frac{\partial}{\partial t_j} |_{t=0} A_t(\zeta) =
 \dot{A}(\zeta) = (\dot{X}(\zeta)+i\dot{Y}(\zeta),0)$ satisfies the
following equation on the unit circle
                    \begin{equation}
         \dot{X}= H_y\circ A \dot{Y} + \chi H_{t_j} \circ A.
                    \end{equation} We also introduce some notations.
For a ${\cal C}^{1, \beta}$-smooth function $g(\zeta)$ on the unit circle
with $g(1)=0$, we write $${\cal J}(g)= \frac{1}{\pi} \int_0^{2\pi}
\frac{g(e^{i\theta})}{|e^{i\theta} - 1|^2} d \theta,$$ where the
integral is understood in the sense of principal value.  Then, if $g
\in C^{1, \beta} (\overline{\Delta})$ is holomorphic in $\Delta$ and
vanishes at $1$, we have
                    \begin{equation}
         {\cal J}(g) = - \frac{ \partial g}{\partial \zeta}(1) = i
       \frac{d}{d \theta} {\mid}_{\theta = 0} g(e^{i \theta}).
                    \end{equation} Notice also that for ${\cal C}^{1, \beta}$
real-valued functions $g, g'$ with $g(1)=g'(1)=0$, applying $(11)$ to
the holomorphic function $(g+iT_1g)(g'+iT_1g')$ vanishing to second
order at $1$, we obtain
                    \begin{equation}
                                  {\cal J}(gg' - T_1g T_1g')=0.
                    \end{equation}

Associate with $A$ and $H$ a $n \times n$ matrix-valued function
$G(\zeta)$ on the unit circle as a solution to the equation
                    \begin{equation}
              G= I - T_1 (GH_y \circ A).
                    \end{equation} The definition of $G$ implies that
$G(1)=I$ and $T_1 G = GH_y \circ A - H_y \circ A(1)= GH_y \circ A$,
since $A(1)=0$ and $h(0)=0, dh(0)=0$.  Using $(14)$ and $G$, we can
write on the unit circle

\smallskip

\hspace{4cm} $G \chi H_{t_j} \circ A = G(\dot{X} - H_y \circ A
\dot{Y})$

\hspace{6cm} $=G\dot{X} - (T_1G)(T_1\dot{X})$

\hspace{6cm} $=\dot{X}+ (G-I)\dot{X} -T_1(G-I) T_1 \dot{X}.$

\smallskip
 
By virtue of $(12)$,
                    \begin{equation}
               {\cal J}(G\chi H_{t_j} \circ A ) =
               {\cal J}(\dot{X}).
                    \end{equation} On the other hand, according to
$(11)$ and the fact that $\dot{Y}=T_1\dot{X}$,
                    \begin{equation}
             -i {\cal J}(\dot{X})+{\cal J}(T_1\dot{X})=
           i\frac{\partial \dot{Z}}{\partial \zeta} (1)=
         \frac{d}{d\theta}|_{\theta=0} (\dot{X}+i\dot{Y}).
                    \end{equation} Identifying the real part of the
two extreme terms and taking $(10)$ into account, we have
                    \begin{equation}
                {\cal J}(T_1\dot{X})=
                   \frac{d}{d\theta}|_{\theta=0} \dot{X} =
                \frac{d}{d\theta}|_{\theta=0}
               (H_y\circ A \dot{Y}+\chi H_{t_j} \circ A)=0,
                    \end{equation} if we choose $\chi$ in order that
$\chi$ is equal to zero near $\zeta=1$ and since $\dot{Y}(1)=0$,
$dH(0)=0$.  $(14)$, $(15)$ and $(16)$ therefore yield
$$-\frac{\partial \dot{Z}}{\partial \zeta}(1) ={\cal J}(G\chi H_{t_j}
\circ A).$$ Natural coordinates on $T_0\C^{m+n} / T_0 M $ being given
by $x_1,...,x_n$, we obtain in these coordinates
                    \begin{equation}
                               \frac{\partial D}{\partial t_j} (0) =
                                \Pi \left(- \frac{\partial
                                \dot{Z}}{\partial \zeta}(1) \right) =
                                {\cal J}(G \chi H_{t_j} \circ A).
                    \end{equation} Furthermore, choose $\chi$ in order
that ${\cal J}(\chi)=1$ and the support of $\chi$ is concentrated near
$\zeta= -1$ so that the vectors ${\cal J}(G \chi H_{t_j}\circ A)$ are
close to the vectors $G(-1)H_{t_j}\circ A(-1)$ and linearly
independent, for $j=1,...,n$ respectively.  This is possible, since
$G$ is non singular at every point on the unit circle and the
$H_{t_j}\circ A(-1)$, $j=1,...,n$ are linearly independent by the
choice of $\kappa$.

The proof of Lemma 2.7 is complete. $\square$

\smallskip When $p=0$, simply denote $A_{t,\tau,a,0}$ by
$A_{t,\tau,a}$.

\smallskip
 
{\sc Lemma 2.8.} {\it $\chi$ can be chosen in order that the following
holds: there exist $\tau_0 >0$, ${\cal T}$ a neighborhood of $0$ in
$\R^q$ and ${\cal A}$ a neighborhood of $0$ in $\C^{p-1}$ such that
the set
                    \begin{equation}
                  \Gamma_0=
                 \{s \frac{dA_{t,\tau,a}}{d\theta} (1);
                 \ s>0, t \in {\cal T},
                  \tau\in I_{\tau_0}, a\in {\cal A}\}
                    \end{equation} is a $(2m+n)$-dimensional open
connected cone with vertex $0$ in $T_0M$.}

\smallskip {\it Proof.} Indeed, the smooth mapping $$E: {\cal T}
\times I_{\tau_0} \times {\cal A} \ni (t,\tau,a) \mapsto
\frac{d}{d\theta}|_{\theta=0} A_{t,\tau,a}(e^{i\theta}) \in T_0M$$
satisfies $rk \ E'(0)=2m+n-1$, for a choice of the function $\chi$ as
in Lemma 1.11, since we defined $W_{t,\tau,a}(\zeta)$ not depending on
$t$, and the partial rank of $E$ with respect to $(\tau,a)$ is equal to
$(2m-1)$. Then the rank of $(s,t,\tau,a)\mapsto sE(t,\tau,a)$ is equal
to $2m+n$.

The proof of Proposition 2.6 is complete. $\square$

\smallskip It easy to see that if one chooses the disc $A$ as in
Proposition 2.6, then a subfamily of a slight modification of
$A_{t,\tau, a, p}$ gives a regular family for $(0,M)$ in the sense of
Definition 1.8.

\smallskip {\it Step three: isotopies.}  Choose a disc
$A_{t_1,\tau_1,a_1}$ such that $v_{00}=\frac{d}{d\theta}|_{\theta=0}
A_{t_1,\tau_1,a_1} (e^{i\theta}) \in T_{0}M_1$, $v_{00}\not\in T_0N$,
$v_{00} \not\in T_0^cM$ and $|v_{00}-v_{0}| <\varepsilon$. This is
possible, since the conification of the set of tangential directions
to discs in the family $A_{t,\tau,a}$ covers an open cone ${\Gamma}_0$
in $T_0M$.  We now prove that isotopy properties are satisfied.

\smallskip
 
{\sc Lemma 2.9.} {\it There exist an open cone ${\Gamma}_1 \subset
\subset {\Gamma}_0$, $v_{00} \in {\Gamma}_1$, a generic ${\cal
C}^2$-smooth one-codimensional submanifold $K$ of $M$ with $K\supset
N$ near $0$ and ${\Gamma}_1\cap T_0 K=\emptyset$ and there exists
${\cal K}$ a neighborhood of $0$ in $K$ such that, if ${\cal
P}_1=\{(t, \tau, a) \in {\cal T} \times I_{\tau_0} \times {\cal A}\ \!
{\bf :} \ \!  \frac{d}{d\theta}A_{t, \tau, a} (1) \in {\Gamma}_1\}$,

{\rm (i)} for each $p \in {\cal K}$ and $(t,\tau,a)\in {\cal P}_1$,
$A_{t, \tau, a, p}(b\Delta)\cap \Phi \neq \emptyset$ if and only if $p
\in {\cal K} \cap \Phi$,

{\rm (ii)} each $A_{t,\tau,a,p}$ for $p\in {\cal K} \backslash \Phi$
and $(t,\tau,a)\in {\cal P}_1$ is analytically isotopic to a point in
$\omega$.}

\smallskip {\it Proof.}  By shrinking all the open sets in the
parameter space, (i) follows, since $v_{00}\not \in T_0K$ and the
embedded disc $A_{t_1,\tau_1,a_1}$ satisfies
$A_{t_1,\tau_1,a_1}(b\Delta\backslash \{1\}) \subset M\backslash
M_1^{-}$, and $\Phi \subset M_1^{-}$.

The second part is as follows.

$(a)$ Each disc $A_{t,\tau,a,p}$ with $p\in {\cal K}\backslash \Phi$
is analytically isotopic in $\omega$ to
$A_{t_1,\tau_1,a_1,p}$. Indeed, the discs $A_{t,\tau,a,p}$ can meet
$\Phi$ only if $p\in \Phi$, by $(i)$.

$(b)$ The discs $A_{t_1,\tau_1,a_1,p}$ for different $p\in {\cal
K}\backslash \Phi$ are analytically isotopic between each other in
$\omega$.  Indeed, $\Phi$ does not divide ${\cal K}$ near $0$, since
$\hbox{codim}_K N =1$ and $0 \in b\Phi$. Then one can take a curve
$p_s, 0\leq s \leq 1,$ in ${\cal K}\backslash \Phi$ between two $p_1,
p_2$ and all the discs $A_{t_1,\tau_1,a_1,p_s}$ do not meet $\Phi$
along their boundaries.

$(c)$ For $u_1^{0}>0$, set $M_{1, u_1^{0}}=
\{u_1=u_1^0+k(v_1,w_2,...,w_p,y)\}$ viewed in $M$.  For small enough
$u_1^0>0$, take $p_2 \in {\cal K} \cap M_{1, u_1^0}$. By $(a)$, we can
assume that after an isotopy, we have a disc $A_{t_2,\tau_2,a_2,p_2}$
such that $\frac{d}{d\theta}A_{t_2,\tau_2,a_2,p_2}(1) \in T_{p_2}
M_{1,u_1^0}$.  Now, push the disc away from $M_{1,u_1^0}$ in the
direction of the positive $u_1$-axis by letting it be attached to $M$
in the region $\{u_1 \geq u_1^0 + k(v_1,w_2,...,w_p,y)\} \cap M$. This
can be done by using a single parameter family through Bishop's
equation. When the disc is far enough from $N$, it is analytically
isotopic to a point in $\omega$.

The proof of Lemma 2.9 is complete. $\square$

\smallskip

{\it Step four: holomorphic extension.}  Let $v_{00} \in C$ be a
$n$-dimensional proper linear cone in the $(2m+n)$-dimensional space
$T_0M$ and contained in ${\Gamma}_1$ such that the projection $T_0 C
\to T_0 M / T_0^cM$ is surjective and $\overline{C} \cap T_0^cM
=\{0\}$.  Let ${\cal P}$ denote the set of parameters $${\cal
P}=\{(t,\tau,a)\in {\cal T} \times I_{\tau_0}\times {\cal V}; \
\frac{d}{d\theta}A_{t,\tau,a}(1) \in C\},$$ which is a ${\cal
C}^1$-smooth $(q-1)$-dimensional submanifold of ${\cal T} \times
I_{\tau_0} \times {\cal V}$.  We choose a similar germ of a manifold,
still denoted by ${\cal P}$, with same tangent space at 0 which is
${\cal C}^2$.  As in \cite{ME2}, one observes that a consequence of
the isotopy property 1.16 and of the fact that the mapping $${\cal
P}\times {\cal K} \times \stackrel{\circ}{\Delta}_1 \ni
(t,\tau,a,p,\zeta) \mapsto A_{t,\tau,a,p}(\zeta) \in
\C^{m+n}\backslash M$$ is a smooth embedding is that ${\cal
H}(\omega)$ extends holomorphically into the open wedge set $${\cal
W}=\{A_{t,\tau,a,p}(\zeta); \ (t,\tau,a)\in {\cal P}, \in {\cal K},
\zeta \in \stackrel{\circ}{\Delta}_1\}$$ minus the set $$\Phi_{\cal
P}=\{A_{t,\tau,a,p}(\zeta)\ \! {\bf :} \ \! (t,\tau,a)\in {\cal P},
p\in {\cal K}\cap \Phi, \zeta \in \stackrel{\circ}{\Delta}_1 \}.$$
Indeed, since the mapping remains injective on ${\cal P} \times {\cal
K}\backslash \Phi \times \stackrel{\circ}{\Delta}_1$, we can set
unambiguously $$F(z):= \frac{1}{2i\pi} \int_{b\Delta} \frac{f\circ
A_{t,\tau,a,p}(\eta)}{\eta - \zeta} d\eta$$ as a value at points
$z=A_{t,\tau,a,p}(\zeta)$ for an extension of $f|_{M\backslash \Phi}$,
$p \in {\cal K} \backslash \Phi, (t,\tau,a)\in {\cal P}, \zeta \in
\stackrel{\circ}{\Delta}_1$.  Since $f$ extends holomorphically to the
interior of these discs, we get a continuous extension $F$ on each
$A_{t,\tau,a,p}(\Delta_1)$, $p\in {\cal K}\backslash \Phi$.  Thus, the
extension $F$ of $f|_{M\backslash \Phi}$ also becomes continuous on
$$({\cal W}\backslash \Phi_{\cal P}) \cup (M\backslash \Phi),$$ where
$\Phi_{\cal P}$ is a proper closed subset of the one generic closed
one codimensional submanifold of ${\cal W}$ $$N_{\cal
P}=\{A_{t,\tau,a,p}(\zeta)\ \! {\bf :} \ \! (t,\tau,a)\in {\cal P},
p\in {\cal K}\cap N, \zeta \in \stackrel{\circ}{\Delta}_1 \}.$$ Since
$f|_{M\backslash \Phi}$ extends analytically to a neighborhood of
$A_{t, \tau, a, p}(\overline{\Delta})$, $F$ is holomorphic in ${\cal
W}\backslash \Phi_{\cal P}$. Indeed, fix a point $(\tilde{t},
\tilde{\tau}, \tilde{a}, \tilde{p}_0) \in {\cal P} \times ({\cal K}
\backslash \Phi)$ and let $\tilde{\cal P} \times \tilde{\cal K}$ be a
neighborhood of $(\tilde{t}, \tilde{\tau}, \tilde{a}, \tilde{p}_0)$ in
${\cal P} \times ({\cal K}\backslash \Phi)$ such that for each $(t,
\tau, a, p)\in \tilde{\cal P}\times \tilde{\cal K}$,
$A_{t,\tau,a,p}(\overline{\Delta})$ is contained in some neighborhood
$\tilde{\omega}$ of $A_{\tilde{t}, \tilde{\tau}, \tilde{a},
\tilde{p}_0}(\overline{\Delta})$ in $\C^{m+n}$ such that there exists
a holomorphic function $\tilde{f} \in {\cal H}(\tilde{\omega})$ with
$\tilde{f}$ equal to $f$ near $A_{\tilde{t}, \tilde{\tau}, \tilde{a},
\tilde{p}_0}(b\Delta)$.  Let $\tilde{\zeta}\in
\stackrel{\circ}{\Delta}_1$ and $\tilde{z}=A_{\tilde{t}, \tilde{\tau},
\tilde{a}, \tilde{p}_0}(\tilde{\zeta})$. To check that the previously
defined function $F$ is holomorphic in a neighborhood of $\tilde{z}$,
we note that for $z=A_{t, \tau, a, p}(\zeta)$, $(t, \tau, a, p) \in
\tilde{\cal P}\times \tilde{\cal K}$, $\zeta$ in some neighborhood
$\tilde{\Delta}_1$ of $\tilde{\zeta}$ in $\stackrel{\circ}{\Delta}_1$,
$\tilde{f}(z)$ is given by the Cauchy integral formula
$$\tilde{f}(z)=\frac{1}{2i\pi} \int_{b \Delta} \frac{\tilde{f} \circ
A_{t, \tau, a ,p}(\eta)}{\eta-\zeta} d\eta= \frac{1}{2i\pi} \int_{b
\Delta} \frac{f\circ A_{t, \tau, a ,p}(\eta)}{\eta-\zeta} d\eta=
F(z).$$ As a consequence, $\tilde{f}(z)=F(z)$ for $z$ in a small
neighborhood of $\tilde{z}$ in $\C^{m+n}$, since the mapping $(t,
\tau, a, p,\zeta) \mapsto A_{t, \tau, a, p}(\zeta)$ from $\tilde{\cal
P}\times \tilde{\cal K} \times \tilde{\Delta}_1$ to $\C^{m+n}$ has
rank $2n$ at $(\tilde{t}, \tilde{\tau}, \tilde{a}, \tilde{p}_0,
\tilde{\zeta})$.

This proves that $F$ is holomorphic into ${\cal W}\backslash
\Phi_{\cal P}$.
 
By shrinking $\omega$ near $0$, which does not modify the possible
disc deformations, we can insure that $\omega \cap {\cal W}$ is
connected, since $\overline{C}\cap T_0^cM=\{0\}$ and then also $\omega
\cap ({\cal W}\backslash \Phi_{\cal P}),$ since $N_{\cal P}$ is a
closed one-codimensional submanifold of ${\cal W}$. Therefore $f\in
{\cal H}(\omega)$ and $F\in {\cal H}({\cal W}\backslash \Phi_{\cal
P})$ stick together in a single holomorphic function in $\omega \cup
({\cal W} \backslash \Phi_{\cal P})$, since both are continuous up to
$M\backslash \Phi$, which is a uniqueness set, and coincide there.

\smallskip Recall the following result of J\"oricke \cite{JO4}.

\smallskip
 
{\sc Lemma 2.10.} {\it Let $U\subset \C^{m+n}$ be an open set,
$m+n\geq 2$, let $M \subset U$ be a connected ${\cal C}^2$-smooth
hypersurface and let $\Phi \subset M$ be a proper closed subset. If
$\Phi$ does not contain any CR orbit of $M$, then, for every function
$f\in {\cal H} (U\backslash \Phi)$, there exists a function $F\in
{\cal H}(U)$ with $F|_{U\backslash \Phi} =f$.}

\smallskip {\it Proof.} Let $\Phi_r$ denote the set of removable
points of $\Phi$ and assume that $\Phi_{{\rm nr}}= \Phi\backslash
\Phi_r$ is nonempty.  Replace $\Phi_{{\rm nr}}$ by
$\Phi\neq\emptyset$. By Lemma 2.3, which applies also in that
situation, there exists a point $p_1\in M\backslash \Phi$ such that
half of $M$ near $p_1$ is contained in $U\backslash \Phi$ and the half
is divided by a complex tangential direction to $M$ at $p_1$. If $M$
is not minimal at $p_1$, $M$ contains a germ through $p_1$ of a
complex hypersurface $H\subset M$, with a half $H^+$ of $H$ contained
in $U\backslash \Phi$. Slightly deform $M$ on that side to drop $H$,
in order that $M$ becomes minimal at $p_1$. We can therefore assume
that $M$ is minimal at $p_1$.  By Tr\'epreau's extension theorem, a
side, say $M^-$, of $M$, has the holomorphic extension property at
$p_1$ to the other side $M^+$.  This proves that $f$ extends through
$\Phi$ at $p_1$. Indeed, the univalency is easy to see.

The proof of Lemma 2.10 is complete. $\square$

\smallskip {\it End of proof of Proposition 2.2.}  Recall that the set
$$N_{\cal P}=\{A_{t,\tau,a,p}(\zeta)\ \! {\bf :} \ \! (t,\tau,a)\in
{\cal P}, p\in {\cal K}, \zeta \in \stackrel{\circ}{\Delta}_1\}$$
%
% K au lieu de K\cap N
% 
is a closed ${\cal C}^2$-smooth one-codimensional submanifold of
${\cal W}$ and that $\Phi_{{\cal P}}$ is a proper closed subset of
$N_{{\cal P}}$. Furthermore, the closed subset $N_{{\cal P}}
\backslash (N_{{\cal P}} \cap \omega)$ of the generic ${\cal
C}^2$-smooth manifold $N_{{\cal P}}$ cannot contain any CR invariant
subset. Indeed, each piece of complex curve $A_{t,\tau,a,p}
(\stackrel{\circ}{\Delta}_1)$ meets $\omega \backslash N_{{\cal P}}$,
since $\omega$ is a neighborhood of $M\backslash N$ in $\C^{m+n}$ and
$\frac{d}{d\theta}|_{\theta=0} A(e^{i\theta})\not\in T_{p_1} N$.
Therefore, Lemma 2.10 yields that ${\cal H}(({\cal W} \backslash
N_{{\cal P}})\cup \omega)$ extends holomorphically into ${\cal W}$.

The proof of Proposition 2.2 is complete. $\square$

\smallskip {\it Remark.} Assume that $T_{p_1} N$ is generic in
$T_{p_1}\C^{m+n}$ and that $m\geq 2$.  By a geometric adaptation of
the approximation theorem given in \cite{BT}, it is easy to prove that
there exist two neighborhoods $U, V \subset \subset U$ of $p_1$ in $M$
such that each function $f\in {\cal C}_{CR}^0(U\backslash N)$ can be uniformly
approximated by holomorphic polynomials on compact subsets of
$V\backslash N$ (\cite{ME2}, Proposition 5.B;
the $L_{CR}^1$ case as in \cite{JO4}).  This shows that the
continuity principle is not needed.  Therefore, Proposition 2.2 holds,
with $\Phi=N$ near a point where $N$ is generic, if ${\rm dim}_{CR} M
\geq 2$.

\smallskip {\it Remark.} An adaptation of Baouendi-Treves approximation
fails necessarily in many cases, {\it e.g.} when 
$m=1$, $\Phi=\{p\}$. Indeed, let $U\ni p$ be a small neighborhood, let
$V\subset \subset U$, let $q_1, q_2\in V\backslash \{p\}$. One should
find maximally real manifolds $L_{q_1} \subset U$, $L_{q_2}\subset U$
with $L_{q_1} \cap (U\backslash V)\equiv 
L_{q_2} \cap (U\backslash V)$  and
a manifold $\Sigma\subset V$ with 
$b\Sigma = L_{q_1}-L_{q_2}$ {\it such that} 
$\Sigma\subset V\backslash \{p\}$ and this
is clearly impossible for arbitrary choices of $q_1$, $q_2$, since
${\rm dim}_{\R} L_{q_1} ={\rm dim}_{\R} M -1=1+n$, 
${\rm dim}_{\R} \Sigma={\rm dim}_{\R} M=2+n$. $\square$

\smallskip

The proof of ${\cal W}$-removability in Theorem 4 is complete now.

\smallskip {\it Step five: End of proof of 
Theorem 4, $L^1$-removability.}  
We got analytic extension of $L_{CR,loc}^1$ functions
on $M\backslash N$ over $M$ and 
to get $L_{CR}^1$ extension over $M$ we need Hardy space
like estimates for the extension as in the nonsingular case Proposition 1.12.

We assume
the situation given as in the end of the proof of Proposition 2.2.  We
first include the disc given in 2.4 in a partial regular family by
choosing parameters $s=(\rho, w^{*0},y^0)$ running in a neighborhood
${\cal S}$ of $0$ in $\R^{2m+n-1}$, $w^{*0}=(w_2^0,...,w_m^0)$,
$y^{0}=(y_1^0,...,y_n^0)$, and by setting
$$A_s(\zeta)=(X_s(\zeta)+iY_s(\zeta), c(1-(1+\rho)\zeta), w^{*0}),$$
where $Y_s$ is the solution of the equation $$Y_s=T_1h(Y_s,
c(1-(1+\rho)\zeta), w^{*0})+y^0.$$ This family satisfies the
requirements of the following Proposition. 

\smallskip
 
{\sc Proposition 2.11.} {\it Let $M$ be generic, ${\cal
C}^{2,\alpha}$-smooth, let $p_1\in M$, let $N$ be a two codimensional
submanifold of $M$ with $p_1 \in N$ and let $\omega$ be a neighborhood
of $M\backslash N$ in $\C^{m+n}$. Assume that there exists a family
$A_s$, $s\in {\cal S}$, ${\cal S}$ a neighborhood of $0$ in
$\R^{2m+n-1}$, of discs such that the mapping ${\cal S} \times b\Delta
\ni (s,\zeta) \mapsto A_s(\zeta) \in M$ is an embedding and
$A(b\Delta) \not\subset N$, $A_0=A$, $A(1)=p_1$ and each $A_s$ is
analytically isotopic to a point.  If ${\cal H}(\omega)$ extends
holomorphically into a wedge ${\cal W}$ of edge a small ball in $M$
around $p_1$ containing $A(b\Delta)$, then $p_1$ is $L^1$ removable.}

\smallskip {\it Proof.}  Extend first $A_s$ in a regular family
$A_{s,v}$ by deforming $A_s$ near a point $A(\zeta)\not\in N$
(Proposition 2.6).  Let $${\cal W}_A=\{A_{s,v}(\zeta)\in \C^{m+n}\ \!
{\bf :} \ \!  (s,v,\zeta) \in {\cal S}_1 \times {\cal V}_1 \times
\stackrel{\circ}{\Delta}_1\}$$ denote a wedge defined by $A_{s,v}$ at
$p_1$ (Proposition 1.12).  Introduce a one parameter family $M^d$,
$d\geq 0$, of small smooth deformations of $M$ into ${\cal W}$,
construct the deformed family of analytic discs $A_{s,v}^d$ and
isotope each of these discs to a point on $M^d$, which is possible,
since the wedge has a big simply connected and round edge containing
$A_{s,v}(b\Delta)$.

This family gives holomorphic extension of ${\cal H}(\omega \cup {\cal
W})$ to ${\cal W}_{A^d}$ for each $d$.  By letting $d$ tend to zero,
the deformed wedge tends in a smooth fashion to ${\cal W}_A$, so we
can control the connectedness of the intersections of the sets where
the extensions are defined.  By taking a smaller set of $v$
parameters, we see that ${\cal H}(\omega \cup {\cal W})$ extends
holomorphically into ${\cal W}_A$. (Notice that two wedges of edge $M$
that are in general position have {\it empty} intersection.)

By Lemma A.6.4 in \cite{CHI}, the closed set ${\cal S}_N=\{s\in {\cal
S}\ \! {\bf :} \ \!  A_s(b\Delta) \cap N \neq \emptyset\}$ is of
Hausdorff dimension $2m+n-2$.

Take a sequence of open neighborhoods $U_n \supset {\cal S}_N$ with
$\hbox{Vol}(U_n) \to 0$ as $n\to \infty$ and a sequence of small
nonnegative functions $\chi_n \in {\cal C}^2({\cal S}_N)$,
$\hbox{supp} \ \chi_n\subset U_n$, positive on ${\cal S}_N$, such that
$\chi_n \to 0$ in ${\cal C}^2({\cal S})$. Define $$M_n= \bigcup_{\zeta
\in b\Delta, s\in {\cal S}} A_s((1-\chi_n(s))\zeta) \subset {\cal W}_A
\cup \omega$$ and define $f_n \in CR(M_n)$ by restriction of the
extension $f\in {\cal H}( {\cal W}_A \cup \omega)$. By subharmonicity
of almost all functions $\zeta\mapsto A_s(\zeta)$ and integrating over
$s\in U_n$ one gets $$\int_{s\in U_n} \int^{\pi}_{-\pi} |f_n
(A_s(1-\chi_n(s)) \zeta)| \ ds d\theta
 < \int_{s\in U_n} \int^{\pi}_{-\pi} |f\circ A_s(\zeta)| \ ds d\theta
\longrightarrow_{n\infty} 0.$$ 
So we conclude that $\lim_{n\infty} f_n
=f$ in $L^1$, hence that $f$ is CR near $p_1$, since each $f_n\in
{\cal H}({\cal V}(M_n))$.

The proof of Proposition 2.11 is complete. $\square$

\smallskip

The proof of Theorem 4 is complete now.

\smallskip
 
{\sc Corollary 2.12.} {\it Let $M$ be a {\rm real analytic} generic
manifold in $\C^{m+n}$ of finite type at every point with ${\rm
dim}_{CR} M =m \geq 2$.  Then every connected real analytic
submanifold $N\subset M$ with $\hbox{codim}_M N =2$ is 
${\cal W}$- and $L^1$-removable if it
does not consist of a CR manifold with ${\dim}_{CR} N=(m-1)$.}

\smallskip {\it Proof.} Since $M$ is everywhere minimal, $N$ cannot
contain any (open) CR orbit of $M$, so the hypotheses of Theorem 4 (iii) are
satisfied. $\square$

\bigskip
 
{\bf 3. Metrically thin singularities.}  This paragraph is devoted to
prove Theorem 2.

\smallskip
 
{\sc Theorem 3.1.} {\it Let $M$ be a locally embeddable ${\cal
C}^{2,\alpha}$-smooth CR manifold, of dimension $d=2m+n$, ${\rm
dim}_{CR} M=m\geq 1$, let $1\leq {\rm p} \leq \infty$, and let $E$ be
a closed subset of $M$ such that the Hausdorff measure $H_{d-3}(K) <
\infty$ for each compact set $K\subset E$. Assume that for almost all
CR orbits, ${\cal O}_{CR}\backslash E$ is globally minimal.  Then $E$
is $L^{\rm p}$-removable.}

\smallskip {\it Proof.}  Fix a function $f\in L_{loc}^1(M)\cap
L_{loc,CR}^1(M\backslash E)$. Call $M_f$ the union of all CR orbits
${\cal O}_{CR}$ of $M$ such that $f|_{{\cal O}_{CR}}$ is locally
integrable on ${\cal O}_{CR}$ and CR on ${\cal O}_{CR}\backslash
({\cal O}_{CR} \cap E)$, such that ${\cal O}_{CR}\backslash ({\cal
O}_{CR} \cap E)$ is globally minimal and such that ${\cal O}_{CR} \cap
E$ is locally of Hausdorff codimension at least three in ${\cal
O}_{CR}$. By a measure-theoretic lemma \cite{CHI}, combined with
arguments in the proof of Theorem 1.1, $M_f$ is of full measure. Therefore it
is now sufficient to prove the theorem in the case that $M$ and
$M\backslash E$ are globally minimal.

Define $${\cal A}= \{\Psi \subset E \ \hbox{closed} ; M\backslash \Psi
\ \hbox{is globally minimal and} \ f\in L_{loc,CR}^1(M\backslash \Psi)
\cap L_{loc}^1(M)\}$$ and define $E_{{\rm nr}}=\cap_{\Psi\in {\cal A}}
\Psi$.  Then $M\backslash E_{{\rm nr}}$ is globally minimal too.
Use for $E_{{\rm nr}}$ the previous notation 
$E$.  According to Proposition 1.16, we can
assume that $f\in L_{loc}^1(M)\cap {\cal H}({\cal V}(M\backslash E))$.
Assume $E\neq\emptyset$. We shall reach a contradiction.

According to the idea of the proof of lemma 2.3, there exists a
generic manifold $M_1$ of codimension one in $M$ through a point
$p_1\in E$ such that $T_{p_1} M_1 \not\supset T_{p_1}^cM$ and $E
\subset M_1^{-}$ locally.  By the definition of ${\cal A}$ and
genericity of $M_1$, it suffices to show that $p_1$ is
$L^1$-removable.  Lemma 2.4 then provides an embedded analytic disc
$A\in {\cal C}^{2,\alpha} (\overline{\Delta})$ with $A(1)=p_1$,
$A(b\Delta\backslash \{1\}) \subset M\backslash M_1^{-}$ and
$\frac{d}{d\theta}|_{\theta=0} A(e^{i\theta}) \in T_{p_1} M_1$.  In a
coordinate system as in $(2)$ and as in 2.4, we can assume that
$A(\zeta)=(Z(\zeta), W(\zeta))$, with
$W(\zeta)=(\rho_1(\zeta-\rho_1),0,...,0)$, $\rho_1 >0$.  For the rest
of the proof below, we only need a disc $A$ with $A(1)\in E$, but
$A(b\Delta) \not\subset E$.

We shall first develope $A$ in a partial regular family $A_{\rho, s'}$
of analytic discs, $0 \leq \rho < \rho_2$, $\rho_2 >\rho_1$,
$I_{\rho_2}=(0,\rho_2)$, $s'=(a_2,...,a_m,y_1^0,...,y_n^0)$ running in
a neighborhood of $0$ in $\C^{m-1}\times \R^m$, as follows:
$W_{\rho,a}(\zeta)=(\rho(\zeta-\rho_1),a_2,...,a_{m-1})$ and
$A_{\rho,s'}(\zeta)=X_{\rho,s'}(\zeta)+
iY_{\rho,s'}(\zeta),\rho(\zeta-\rho_1),a_2,...,a_{m-1})$, where
$Y_{\rho,s'}=T_1h(Y_{\rho,s'},W_{\rho,s'})+y^0$.  Then there exists
${\cal V}$, a neighborhood of $0$ in $\C^{m-1}$, ${\cal Y}$, a
neighborhood of $0$ in $\R^n$, such that the mapping $$I_{\rho_2}
\times {\cal A} \times {\cal Y} \times b\Delta \ni (\rho,a,y^0,\zeta)
\mapsto A_{\rho,a,y^0}(\zeta) \in M$$ is an embedding. This shows that
 a neighborhood in $M$ of $A_{0,0,0}(\overline{\Delta})=p_1$ is
foliated by ${\cal C}^2$-smooth real discs $D_{a,y^0}=D_{s'}=
\{A_{\rho,a,y^0}(\zeta)\in M\ \! {\bf :} \ \! 0\leq \rho < \rho_2,
\zeta \in b \Delta\}$.  Moreover, since $H_{d-3}(E) < \infty$, the set
${\cal S}_E'=\{s'\in {\cal S}'\ \! {\bf :} \ \! D_{s'} \cap E \neq
\emptyset\}$ is a closed subset of ${\cal S}'={\cal A} \times {\cal
Y}$ of Hausdorff codimension $\geq 1$. By construction, each disc
$A_{\rho,s'}$ with $s'\not\in {\cal S}_{E'}$ is therefore analytically
isotopic to a point in $M\backslash E$.

Furthermore, by means of normal deformations of the family near
$A(-1)$ as in Proposition 2.6, we can develope $A$ in a regular family
$A_{\rho,s',v}$ which has the property that, for each $v\in {\cal V}$,
the set ${\cal S}_{E,v}'= \{s'\in {\cal S}'\ \! {\bf :} \ \! D_{s',v}
\cap E\neq \emptyset\}$ is a closed subset of ${\cal S}'$ of Hausdorff
codimension $\geq 1$.  Therefore, each disc $A_{\rho,s',v}$ with
$A_{\rho,s',v}(b\Delta)\cap E = \emptyset$ is analytically isotopic to
a point in $M\backslash E$, since ${\cal S}'\backslash {\cal
S}_{E,v}'$ is dense and open in ${\cal S}'$.

Then the isotopy property and Proposition 1.6 imply that ${\cal
H}({\cal V}(M\backslash E))$ extends holomorphically into $${\cal
W}=\{A_{\rho,a,y^0,v}(\zeta)\in \C^n\ \! {\bf :} \ \! \rho \in
I_{\rho_1}, a\in {\cal A}_1, y^0\in {\cal Y}_1, v\in {\cal V}_1, \zeta
\in \Delta_1\}$$ minus the set $$E_{{\cal W}}=\{A_{\rho,s',v}(\zeta)
\in \C^n\ \! {\bf :} \ \! A_{\rho,s',v}(b\Delta) \cap E \neq
\emptyset\}$$ for which $H_{2m+2n-2}(E_{{\cal W}}) < \infty$.  If even
$H_{2m+2n-3}(E_{{\cal W}}) <\infty$, then $E_{{\cal W}}$ is removable
and we are done.  Let $f\in {\cal H}(\omega)$ and let $F$ denote its
extension to ${\cal W}\backslash E_{{\cal W}}$.

By using a subharmonic estimate as in 1.12 and 1.13, we get $F\in
L^1({\cal W})$.  If $P$ is a polydisc into some $\C^n$, we call $D$ a
{\it coordinate disc} if $D$ is the intersection of $P$ with a
coordinate line in $\C^n$.

\smallskip
 
{\sc Lemma 3.2.} {\it Let $p_0\in {\cal W}$ and let $p_0 \in P \subset
\subset {\cal W}$ be a polydisc. Then, for almost every coordinate
disc $D \subset P$, $D \cap E_{{\cal W}}$ is a union of isolated
points and $F|_D$ is meromorphic with poles of order at most one.}

\smallskip {\it Proof.} By Lemma A.6.4 from \cite{CHI}, for almost
every coordinate disc $D$, $D \cap E_{{\cal W}}$ is of Hausdorff
dimension zero, {\em i.e.} is a union of isolated points.  Fubini's
theorem implies that $F|_D\in L^1(D)$ for almost all $D$. To show
meromorphicity, if $g\in {\cal O}(\Delta \backslash \{0\})\cap
L^1(\Delta)$, then $\overline{\partial}(zg)=0$ in the distribution
sense, which completes the proof.

\smallskip We can then apply the following theorem of Shiffmann on
separate meromorphicity to get that $F\in L^1({\cal W}) \cap {\cal
H}({\cal W}\backslash E_{{\cal W}})$ is {\it meromorphic} on ${\cal
W}$. A subset $Q$ of a polydisc $P$ is said to be a {\it full subset}
of $P$ if $Q \cap D$ is a set of full measure in $D$ for almost every
coordinate disc $D\subset P$.

\smallskip
 
{\sc Theorem 3.3.} ({\sc Shiffman}, \cite{SH}). {\it Let $P \subset
\subset \C$ be a polydisc and let $Q \subset P$ be a full subset of
$P$.  Then a function $F: Q \to \P^1(\C)$ has a meromorphic extension
in $P$ if and only if, for almost every coordinate disc $D \subset P$,
$F|_{D\cap Q}$ extends meromorphically in $D$.}

\smallskip Assume that $F$ is not holomorphic in ${\cal W}$, {\em
i.e.} the polar variety $P_F=\{p\in {\cal W}; F(p) \not\in \C\}$, of
pure dimension $(n-1)$, is nonempty.  We shall prove that either
$A_{s,v}(\stackrel{\circ}{\Delta}_1)$ does not intersect $P_F$ or is
contained in it.

\smallskip
 
{\sc Lemma 3.4.} {\it If $A_{s,v}(\stackrel{\circ}{\Delta}_1)\cap
P_F\neq \emptyset$, then $A_{s,v}(\stackrel{\circ}{\Delta}_1)\subset
P_F$.}

\smallskip {\it Proof.} Indeed, if $A_{s,v}(\zeta_1) \in P_F$ and the
intersection number of $P_F$ with
$A_{s,v}(\stackrel{\circ}{\Delta}_1)$ at $A_{s,v}(\zeta_1)$ is a
finite positive number, $P_F$ will intersect every nearby disc
$A_{\tilde{s},\tilde{v}}(\stackrel{\circ}{\Delta}_1)$, contradicting
the fact that $F \in {\cal H}({\cal
V}(A_{\tilde{s},\tilde{v}}(\stackrel{\circ}{\Delta}_1)))$ for almost
all $(\tilde{s},\tilde{v})\in {\cal S} \times {\cal V}$.

\smallskip {\it End of proof of Theorem 3.1.} Consider the manifolds
with boundary near $p_0$ which foliate ${\cal W}$,
$M_v=\{A_{s,v}(\zeta)\in \C^n\ \! {\bf :} \ \! s\in {\cal S}_1, \zeta
\in \stackrel{\circ}{\Delta}_1\},$ for $v\in {\cal V}_1$. Since the
$(2m+2n-2)$-dimensional measure of $P_F$ is positive, for almost all
$v\in {\cal V}_1$, the $(2m+n-1)$-dimensional measure of $P_F \cap
M_v$ is strictly positive.  Fix $v_1 \in {\cal V}_1$ with that
property.  Then the set ${\cal S}_2$ of $s\in {\cal S}_1$ such that
$A_{s,v}(\stackrel{\circ}{\Delta}_1) \subset P_F \cap M_{v_1}$ is of
positive $(2m+n-3)$-dimensional Hausdorff measure, because of Lemma
3.4. Since $E$ is of Hausdorff codimension at least 3 in $M$ and
$\{A_{s_1,v_1}(b\Delta\cap \Delta_1)\in M; \ s_2 \in {\cal S}_2\}$ is
of positive $(2m+n-2)$-dimensional measure, there exists $s_2 \in
{\cal S}_2$ such that $A_{s_2,v_1}(b\Delta \cap \Delta_1)\not\subset
E$, say $A_{s_2,v_1}(\zeta_0) \not\in E$.  But
$A_{s_2,v_1}(\stackrel{\circ}{\Delta}_1) \subset P_F$ and $F$
holomorphic in a neighborhood of $A_{s_2,v_1}(\zeta_0)$ in $\C^n$
gives the desired contradiction.  Finally, the proof of Proposition
2.11, applies verbally to that situation and yields the desired
$L^1$-removability of $E$ near $p_1$.

The proof of Theorem 3.1 is complete. $\square$
\bigskip
 
{\bf 4. Orbit decomposition.} In this section, we shall prove the
technical Theorem 6 about orbits. Since the result is a real differential
geometric lemma, let us start as follows.

Let $M$ be a real ${\cal C}^2$-smooth manifold of dimension $m+n$, $n\geq 1$,
$m\geq 2$ and let $K \subset TM$ be a ${\cal C}^1$-smooth real subbundle of
rank $m$.  Given a section $L$ of $K$, we denote by $\R \ni t \mapsto
L_t(p)$ an integral curve of $L$ with origin $p$, and call it a
$K$-curve.  According to the analysis of Sussmann \cite{SU}, orbits
${\cal O}_K(M, p)$ under $K$ of points $p\in M$ can be naturally
equipped with a structure of a ${\cal C}^1$-smooth manifold making the
canonical inclusion $i : {\cal O}_{K}(M,p) \to M $ an injective
immersion of class ${\cal C}^1$.

By a {\it leaf $\omega_p$ of ${\cal O}_K(M,p)$ through $p$}, we shall mean
an open small neighborhood of $p$ in ${\cal O}_K(M,p)$ for the
topology of the orbit, such that $\omega_p$ is a ${\cal C}^1$-smooth
submanifold of $M$ too.  By a {\it $K$-integral manifold $S$}, we shall mean
a submanifold with the property that $T_q S \supset K(q)$, $\forall
\ q\in S$. Orbits are immersed $K$-integral submanifolds.

{\it Local $K$-orbits} are defined to be the inductive limit, as $U$
ranges through the set of open neighborhoods of $p$ in $M$, of the
$K$-orbit of $p$ in $U$. Since the dimension of orbits is a
well-defined integer $\geq 2$, these stabilize and define the unique
germ through $p$ of a $K$-integral submanifold of $M$ with minimal
possible dimension. We call $M$ {\it $K$-minimal at $p$} if ${\cal
O}_K^{loc} (p)$ is an open neighborhood of $p$ in $M$, {\em i.e.} has
maximal possible dimension.  We shall prove the following.

\smallskip
 
{\sc Theorem 4.1.} {\it Let $M$ be a real ${\cal C}^2$-smooth
manifold, $\hbox{dim}_{\R} M =m+n$, $m\geq 2$, $n\geq 1$, let
$K\subset TM$ be a ${\cal C}^1$ subbundle of rank $m$, let $N\subset
M$ be a ${\cal C}^2$-smooth submanifold, $\hbox{codim}_M N \geq 2$,
$T_p N \not\supset K(p)$, $\forall \ p\in N$, let $N^c=\{p\in N\ \! {\bf
:} \ \! \hbox{dim}_{\R} K(p) \cap T_pN=m-1\}$ and let $\Upsilon$
denote the set of ${\cal C}^1$ sections of $K$ such that $Y|_{N^c}$ is
tangent to $N$. If, for each $p\in N$, ${\cal O}_{\Upsilon}(M,p)$ is
not contained in $N$, then every K-orbit of $M\backslash N$ is given
by ${\cal O}_{K} \backslash N$, for some K-orbit of $M$.}

\smallskip {\it Proof.} Theorem 4.1 can be reduced to two lemmas, 4.2
and 4.4 below.

\smallskip
 
{\sc Lemma 4.2.} {\it Let $M$, $N$
 be as in Theorem 4.1. Then, for each $q\in N$, there exists
$\omega_q$ a $K$-integral manifold through $q$ such that $T_q \omega_q
\not\subset T_q N$, $\omega_q \cap N \subset H_q$, $H_q$ a one
codimensional closed submanifold of $\omega_q$, and $\omega_q
\backslash (\omega_q \cap N)$ is contained in a single $K$-orbit of
$M\backslash N$.}

\smallskip {\it Proof.} Remark that $\omega_q\cap N$ is thin in
$\omega_q$. Remark also ${\cal O}_K(M\backslash N, p) \subset {\cal
O}_K(M,p)\backslash N$, if $p\in M\backslash N$ (obvious). Without
loss of generality, we shall denote by $[0,t_1] \ni t\mapsto L_t(p)$
any piecewise smooth $K$-integral curve, as if there were a single
smooth piece.

Since ${\cal O}_{\Upsilon}(M,q) \not\subset N$, there exists a point
$q_0 \in {\cal O}_{\Upsilon}(M,q)$, $q_0=L_{s_0}(q)$, such that
$L_s(q) \in N$, $\forall \ s, 0\leq s \leq s_0$ and ${\cal
O}_{\Upsilon}^{loc}(M,q_0) \not\subset N$. Choose a germ of a ${\cal
C}^2$-smooth generic one codimensional manifold $M_1 \subset M$
through $q_0$ such that $M_1$ contains $N$. Then $K|_{M_1} \cap TM_1$
defines a ${\cal C}^1$ vector bundle $K_1$ of rank $m-1$ and we have
that ${\cal O}^{loc}_{K_1}(M_1, q_0 ) \not \subset N$.  Indeed,
otherwise, ${\cal O}^{loc}_{K_1}(M,q_0)=S$ satisfies $S\subset N$,
$\hbox{dim}_{\R} (T_q S \cap K(q))=m-1$, $\forall \ q \in S$,
therefore ${\cal O}_{\Upsilon}^{loc}(M,q_0) \subset S \subset N$.

Choose a section $L'$ of $K$ near $q_0$ with $L'(q_0) \not\in
T_{q_0}N$. Let $q_1=L_{t_1}(q_0)\in M_1\backslash N$ be the endpoint
of a piecewise smooth integral curve of $K_1$-tangent vector fields
and extend $L$ in a neighborhood of $M_1$ in $M$. Choose a leaf
$\omega_1$ of ${\cal O}_K(M\backslash N,q_1)$ through $q_1$. Set
$\omega_t=L_{t-t_1}(\omega_1)$, $0\leq t \leq t_1$,
$\omega_{t_1}=\omega_1$. Then we have the
following. Notice that since $T_{q_1} M_1 \not\supset K(q_1)$ and
since $L'(q_1)\in T_{q_1} \omega_1$, then $T_{q_1} \omega_1
\not\subset T_{q_1} M_1$.

\smallskip
 
{\sc Lemma 4.3.} {\it For every $t$, $0 \leq t \leq t_1$, $T_{q_t}
\omega_t \not\subset T_{q_t} N$, $\omega_t$ is a $K$-integral manifold
and $\omega_t \backslash (\omega_t \cap N)$ is contained in a single
$K$-orbit of $M \backslash N$.}

\smallskip {\it Proof.} Since $L|_{M_1}$ is tangent to $M_1$, the flow
of $L$ stabilizes $M_1$. Therefore $T_{q_t} M_1 + T_{q_t}
\omega_t=T_{q_t}M$, for every $0 \leq t \leq t_1$. This implies also that
for any point $r\in \omega_1 \backslash N$, $L_{t-t_1}(r)$ does not
meet $N$ also. So $\omega_t \backslash M_1$ is a $K$-integral
manifold, since everything flows in $M\backslash N$
and since we choosed $\omega_1$ as a leaf of ${\cal O}_K(M\backslash N, q_1)$
through $q_1$. Moreover, since
$H_t=\omega_t\cap M_1$ is a one codimensional submanifold, the closure
$\omega_t$ of $\omega_t \backslash H_t$ still is a $K$ integral
manifold of $M$. Since $\omega_t \backslash M_1$ is contained in a single
$K$-orbit of $M\backslash N$, all points of $M_1 \cap \omega_t$ not in
$N \cap \omega_t$ can be reached by means of $L'$, which completes the
proof. $\square$

We have $q_0=L_{s_0}(q)$.  A repetition of 4.3 along $L_{s-s_0}(q_0)$,
$0\leq s \leq s_0$, gives 4.2.
 
The proof of Lemma 4.2 is complete. $\square$

Let $p\in M\backslash N$ and let $p_1 \in {\cal O}_{K}(M,p)\backslash
N$.  The following lemma shows $p_1\in {\cal O}_K(M\backslash N,p)$.

\smallskip {\sc Lemma 4.4.} {\it Let $p\in M\backslash N$, let $t
\mapsto L_t(p)$ be piecewise smooth integral curve of $K$, set
$p_t=L_t(p)$, $t\in I=[0,t_1]$. Then, for every $t\in I$, there exists
a $K$-integral submanifold $\chi_t$ through $p_t$ such that $\chi_t
\backslash (\chi_t \cap N) \subset {\cal O}_K(M\backslash N,p)$.}

\smallskip {\it Proof.} Let $E$ denote the set of $t\in I$ such that,
for every $s\leq t$, there exists a $K$-integral manifold $\chi_s$
through $p_s$ with $\chi_s \backslash (\chi_s \cap N) \subset {\cal
O}_K(M\backslash N,p)$. Then $E$ contains $[0,\delta)$, for small
$\delta >0$, since $p\in M\backslash N$.  $E$ is open, by definition.
Indeed, let $t \in E$ with $0< t < t_1$, let $\chi_t$ be a
$K$-integral manifold through $p_t$. For every small $\delta>0$,
$L_{\delta}(p_t)\in \chi_t$, so a neighborhood $\chi_{t+\delta}$ of
$p_t$ in $\chi_t$ satisfies $\chi_{t+\delta} \backslash
(\chi_{t+\delta} \cap N) \subset {\cal O}_K(M\backslash N,p)$.

To prove closedness, let $t \leq t_1$ such that $s\in E$, for every $s
<t$. Let $J=\{t\in I\ \! {\bf :} \ \! p_t \in N\}$ and let
$\stackrel{\circ}{J}$ denote the interior of $J$. Assume first that
$(t-\delta, t)$ is contained in $\stackrel{\circ}{J}$, for small
$\delta >0$.  Let $\omega_s$ be as in Lemma 4.2.  We can replace
$\chi_s$ by $\omega_s$, for each $s\in (t-\delta,t)$. Indeed, notice
that ${\cal O}_K^{loc}(p_s)\subset \chi_s \cap \omega_s$, since
$\omega_s$ and $\chi_s$ are $K$-integral manifolds through $p_s$. But
${\cal O}_K^{loc}(p_s)\not\subset N$, since $T_{p_s} N \not\supset
T_{p_s}^cM$.  Since both $\chi_s \backslash (\chi_s \cap N)$ and
$\omega_s \backslash (\omega_s \cap N)$ are contained in a single $K$
orbit of $M\backslash N$, this yields that $\omega_s\backslash
(\omega_s \cap N) \subset {\cal O}_K(M\backslash N,p)$. Then $\omega_t
\backslash (\omega_t \cap N) \subset {\cal O}_K(M\backslash N,p)$
also, since $\omega_t=L_{\delta'}(\omega_{t-\delta'})$, $0< \delta'
<\delta$.  Assume now that there exists $\delta >0$ such that
$L_{-\delta}(p_t) =p_s \not\in N$. Since $\omega_t$ is a $K$-integral
manifold, $L_{-\delta}(p_t)\in \omega_t$. But $\chi_s \subset {\cal
O}_K(M\backslash N,p)$. So a neighborhood of $p_s$ in $\omega_t$ is
also contained in a single $K$ orbit of $M\backslash N$, namely ${\cal
O}_K (M\backslash N, p)$.  It suffices to take $\chi_t=\omega_t$.

The proof of Theorem 4.1 is complete. $\square$

\medskip \noindent {\bf 5. ${\cal C}^{\lambda}$ peak sets.} Here, we
the constructions are simpler than those about global minimality and
continuity principle, since they do not involve considerations of
envelope of holomorphy. Let $S\subset M$ be a ${\cal C}^{\lambda}$
peak set, $0< \lambda< 1$, {\it i.e.} there exists a {\it nonconstant}
function $\varpi\in {\cal C}^{\lambda}_{CR}(M)$ such that $S=\{\varpi
=1\}$ and $|\varpi| \leq 1$ on $M$.  We use three lemmas from
\cite{KYREA}, after which Theorem 5 is an 
easy consequence. Recall that a regular family of analytic discs $A_{s,v}$
generates a wedge ${\cal W}={\cal W}_{A,p}$.

\smallskip

{\sc Lemma 5.1.} (\cite{KYREA}.) {\it Let $M$ be ${\cal C}^{2,\alpha}$
globally minimal and let $p\in M$. If $A_{s,v}$ is a regular family of
analytic discs at $p$, then for each ${\cal C}^{\lambda}$ peak
function $\varpi: M\to \C$:

$1)$ \ $|\varpi|< 1$ in ${\cal W}_{A,p}$ and $\varpi\circ A_{s,v} < 1$
on $\Delta${\rm ;} and{\rm :}

$2)$ \ $\varpi\circ A_{s,v} < 1$ almost everywhere on $b\Delta$.}

\smallskip {\it Proof.} Of course, we can assume that $S\neq
\emptyset$. By the {\it maximum principle on discs}, $|\varpi \circ
A_{s,v}(\zeta) | \leq 1$ on $\overline{\Delta}$, since $|\varpi \circ
A_{s,v}|_{b\Delta}|\leq |\varpi|_M| \leq 1$.  Hence $|\varpi|\leq 1$
in ${\cal W}_{A,p}=\{ A_{s,v}(\zeta) \ \! {\bf :} \ \! (s,v,\zeta)\in
{\cal S}_1 \times {\cal V}_1\times \stackrel{\circ}{\Delta}_1\}$. But
$|\varpi|_S| \equiv 1$, $S\neq \emptyset$, so, if $|\varpi|=1$ at a
point of ${\cal W}_{A,p}$ (or if $|\varpi \circ A_{s,v}(\zeta_0) | =
1$ at some $\zeta_0\in \Delta$), then $\varpi \equiv 1$ in ${\cal
W}_{A,p}$, again because of the maximum principle, so $\varpi|_M\equiv
1$, contradiction. $\square$

\smallskip

{\sc Lemma 5.2.} (\cite{KYREA}, Lemma 3.)  {\it Let $L=\sum_{j=1}^d
a_j(x) \partial / \partial x_j$, $a_j(x) \in {\cal C}^1(\Omega,\C)$,
$\Omega\subset \R^r$ open, let $u\in {\cal C}^{\lambda}(\Omega)$,
$0<\lambda<1$, be a solution of $Lu=0$. Then for every $f\in
L_{loc}^1(\Omega)$ such that $Lf=0$ on $\Omega\backslash \{u=0\}$, one
has $L(u^{2/\lambda}f)=0$ in $\Omega$.} $\square$

\smallskip {\sc Lemma 5.3.} (\cite{KYREA}, Lemma 4.) {\it If $g\in
{\rm H}_{\rm a}^1(\Delta)$, $u\in {\cal H}(\Delta)$, ${\rm Re} \
u^{\lambda/2}>0$ in $\Delta$ and $g=uf$ on $b\Delta$ with $f\in
L^1(b\Delta)$, then $g/u\in {\rm H}_{\rm a}^1(\Delta)$.} $\square$

\smallskip Now, let $f\in L_{loc}^1(M)\cap L_{loc,CR}^1(M\backslash
S)$ where $S=\{\varpi=1\}$ is ${\cal C}^{\lambda}$ peak. Let $p\in M$,
$A_{s,v}$ regular at $p$.  Then $g:=(1-\varpi)^{2/\lambda}f \in
L_{loc,CR}^1(M)$ by Lemma 5.2, hence admits a holomorphic extension in
${\rm H}_{\rm a}^1({\cal W}_{A,p})$, thanks to Proposition 2.  But
$u:=(1-\varpi)^{2/\lambda}\neq 0$ in ${\cal W}_{A,p}$ and ${\rm Re} \
u^{\lambda/2}= {\rm Re} \ (1-\varpi)>0$ and also ${\rm Re} (u\circ
A_{s,v})^{\lambda/2}>0$ in $\Delta$ thanks to Lemma 5, hence
$(g/u)\circ A_{s,v}\in {\rm H}_{\rm a}^1(\Delta)$ by Lemma 7 and
finally $g/u \in {\rm H}_{\rm a}^1({\cal W}_{A,p})$. $\square$

\bigskip
\noindent
{\bf 6. Hypoanalytic structures.} In this section, 
$M$ will be a $d$-dimensional manifold, of class ${\cal C}^{k,\alpha}$, 
$k\geq 2$, $0<\alpha <1$, countable at infinity. 
A {\it hypoanalytic structure on $M$} (\cite{TREV1}, \cite{TREV2}) 
means the data of an open covering $\{U^j\}_{j\in J}$ of $M$ and, 
for each $j\in J$, the data of $n$ complex-valued ${\cal C}^{k,\alpha}$ functions
$Z_1^j,...,Z_n^j$ in $U^j$ {\it with 
linearly independent differentials over $U^j$}, with 
$n\geq 1$ independent of $j$, such that
the following is satisfied: whenever
$U^j \cap U^k \neq \emptyset$, there
is a holomorphic mapping
$F_k^j$ from a neighborhood of $Z^j(U^j \cap U^k)$
in $\C^{n}$ such that $Z^k=F_k^j \circ Z^j$ in $U^j \cap U^k$. The integer
$m=d-n$ is called the {\it codimension} 
of the hypoanalytic structure. Smooth 
{\it hypoanalytic functions} are complex-valued functions
which are locally pullbacks by $Z^j$ of holomorphic functions 
on open neighborhoods of $Z^j(U^j)$ in $\C^{n}$. 
We shall say that $f$ is a {\it RC function}
if $df \in \Gamma(T',U)$, where $T'$ is defined as follows.

If $(U,Z_1,...,Z_n)$ is a hypoanalytic 
local chart in $M$, the differentials $dZ_1,...,dZ_n$ span
a complex vector subbundle $T_U'$ of
$\C T^* U$. According to the compatibility 
condition, this defines a complex vector subbundle
$T'$ of $\C T^*M$ of rank $n$, called
the {\it structure bundle}. Under 
the duality between tangent
and cotangent vectors, it is equivalent to give the 
formally integrable complex vector subbundle 
${\cal V} =(T')^{\bot}$ of
$\C TM$ of rank $m$. 

Let $(U,Z_1,...,Z_n)$ be a hypoanalytic chart in $M$. Possibly after contracting
$U$ about one of its points $p_0$, we can find
$m$ ${\cal C}^{k,\alpha}$ real-valued functions
$(y_1,...,y_m)$ in $U$ such that 
$(dZ_1,...,dZ_n,dy_1,...,dy_m)$ span the whole cotangent space $\C T^*_pM$ at every 
point of $U$. We then define $m+n$ ${\cal C}^{k,\alpha}$ complex vector
fields in $U$, $L_1,...,L_m$ and
$M_1,...,M_n$ by
$M_j Z_k=\delta_j^k$,
$L_iy_l=\delta_l^i$, $M_j y_l=0$, 
$L_i Z_k=0$, 
$i,l=1,...,m$,
$j,k=1,...,n$. Then the commutation relations 
$[L_i,L_l]=[L_i,M_j]=[M_j,M_k]=0$ follow
at once. Put $x_1=\hbox{Re} \ 
Z_1,..., x_n=\hbox{Re} \ Z_n$. Then 
$(x_1,...,x_n,y_1,...,y_m)$ is a 
system of local real coordinates for $M$ in the 
neighborhood $U$ of $p_0$. We can 
assume that the hypoanalytic 
functions $Z_1,...,Z_n$ vanish at $p_0$.

Let $r=n-\hbox{dim} \ T^*M \cap T'[p_0]$. We can make
a $\C$-linear transformation on the 
functions $Z_1,...,Z_n$ and choose
coordinates $(x_1,...,x_n,y_1,...,y_m)$ 
on a neighborhood $U$ of $p_0$
with $x_j=\hbox{Re} \ Z_j$, $j=1,...,n$ and $y_k=\hbox{Im} \ Z_k$, 
$k=1,...,r$, such that
$$Z_k(x,y)=x_k+iy_k, \ \ k=1,...,r, \ \ \ \ \
Z_{r+k}(x,y)=x_{r+k}+ih_{r+k}(x,y), \ \ k=1,...,n-r,$$
with $h_k(0,0)=0$, $dh_k(0,0)=0$, $k=r+1,...,n$.

Therefore, a basis of $(T')^{\bot}={\cal V}$ is given by 
                    \begin{equation}
                                  L_j=\frac{1}{2}
\left( \frac{\partial }{\partial x_j}
+i\frac{\partial }{\partial y_j}\right)+
\sum_{k=r+1}^n a_{jk}(x,y) \frac{\partial }{\partial x_k}, 
\ \ \ \ \ j=1,...,r,                                         
                    \end{equation}
$$L_{r+j}=\frac{i}{2} 
\frac{\partial }{\partial y_{r+j}}+
\sum_{k=r+1}^n a_{r+j \ k}(x,y)\frac{\partial }{\partial x_k}, 
\ \ \ \ \ j=1,...,m-r,$$
with $a_{jk}(0)=0$, the $a_{jk}$ being determined by
$L_j Z_k=0$,
$j=1,...,m$, $k=1,...,n$, and we 
have $L_j y_l= \frac{i}{2} \delta_l^j$. 
We are now prepared to speak of the 
local CR lifting of a hypoanalytic manifold.

The intersection $T^0=T' \cap T^* M$, called
the {\it characteristic set} of the 
differential operators in 
${\cal V}=(T')^{\bot}$ is {\it not}, 
in general, a vector bundle, nor is
$\hbox{Re} \ {\cal V}=
(T'\cap T^*M)^{\bot}$. However, it has 
the property that the integer-valued function
$M\ni p \mapsto r(p) = n -\hbox{dim} \ T^* M \cap T'[p]$ is 
lower semicontinuous.

Let $r=n-\hbox{dim} \ T^* M \cap T' [p_0]$.
 Equivalently, $r=r(p_0)$ is the integer
which satisfies $\hbox{dim} \ 
\hbox{Re} \ T'[p_0]=n+r(p_0)$.

\medskip
\noindent
{\sc 6.1. Definition.} By a {\it regular 
hypoanalytic chart $(Z,y,U)$ at $p_0$}, we mean
the data of a hypoanalytic chart 
$Z: U \to \C^{n}$, $U\ni p_0$, together
with $m-r$ real-valued functions 
$y_{r+1},...,y_m: U \to \R$,
$r=\hbox{dim} \ \hbox{Re} \ T'[p_0]-n$, 
such that $(d\hbox{Re} \ Z_1,...,d\hbox{Re} \ Z_n,
d\hbox{Im} \ Z_1,...,d\hbox{Im} \ Z_n, dy_1,...,dy_m)$ generate
$T_p^*M$ at every point $p\in U$. 
To a regular hypoanalytic 
chart $(Z,y,U)$ is associated a
local embedding of $U$ in $\C^{n} \times \R^{m-r}$,
$p\mapsto (Z(p),y(p))$, which is 
also a real chart on $U$.

\medskip
If $(Z,y,U)$ and $(Z,\tilde{y},U)$ 
are two regular hypoanalytic charts
at $p_0$, there exists a local 
diffeomorphism $G$ at $0$ in $\C^{n}\times \R^{m-r}$
of the form $G(z,y)=(z,G_y(z,y))$ such that
$(Z(p), \tilde{y}(p))=(Z(p), G_y(Z(p), y(p)))$, so 
$(Z, \tilde{y})_* {\cal V}|_U=
G_*(Z,y)_* {\cal V}|_U$. 
Therefore, the subbundle $(Z,y)_* 
{\cal V}|_U ={\cal V}|_{U'}' \subset \C T'$ does
not depend on the choice of 
$(y_1,...,y_{m-r})$.

Then the push forward by $(Z,y)=Z'$ 
of the vector fields given by 
(20) are the restriction to the 
manifold $U'$ of the vector fields on 
$\C^{n} \times \R^{m-r}$, 
equipped with coordinates $(z_1,...,z_n,y_{r+1},...,y_m)$, 
given by
                    \begin{equation}
                  L_j'=\frac{\partial }{\partial \bar{z}_j}+
\sum_{k=r+1}^n b_{jk}(z,\bar{z},y) 
\frac{\partial }{\partial \bar{z}_k}, \ \ \ \ \ j=1,...,r,                                         
                    \end{equation}
$$L_{r+j}'=\frac{i}{2}
\frac{\partial }{\partial y_{r+j}}+
\sum_{k=r+1}^n b_{r+j \ k}(z,\bar{z},y) 
\frac{\partial }{\partial \bar{z}_k}, \ \ \ \ \ j=1,...,m-r,$$
that are $L_j'|_{U'}=Z_*'L_j$, $j=1,...,m$, with 
$b_{jk}$ of class ${\cal C}^{k-1,
\alpha}$ near $0$ satisfying $b_{jk}(0)=0$.

For every completion of the coordinate system 
$(y_{r+1},...,y_m)$ at a point 
$p_0\in M$, $r=r(p_0)$, for which 
$(d\hbox{Re} \ Z_1,...,d \hbox{Re} 
\ Z_n, d\hbox{Im} \ Z_1,...,d\hbox{Im} \ Z_n, 
dy_{r+1},...,dy_m)$ span $T^*M$ 
in a neighborhood of $p_0$, we introduce
extra real variables $(x_{n+1},...,x_{n+m-r})$ 
in $\R^{m-r}$ and a mapping
$\tilde{Z}:U \times V=\tilde{U} 
\ni (p,x) \mapsto \tilde{Z} (p,x)\in \C^{n+m-r}$,
$V$ a neighborhood of $0$ in $\R^{m-r}$, defined by
$$\tilde{Z}_k(p,x)=Z_k(p), \ \ k=1,...,n, \ \ \ \ \
\tilde{Z}_{n+k}(p,x)=x_{n+k}+iy_{r+k}(p), \ \ k=1,...,m-r.$$
Then the image $\tilde{Z}(\tilde{U})=
\tilde{M}$ is a germ
of a CR generic submanifold
in $\C^{n+m-r}$, equipped with 
coordinates $(z_1,...,z_n,w_{n+1},...,w_{n+m-r})$,
given by the equation
$$\hbox{Im} \ z_k=h_k(z_1,...,
z_r,w_{n+1},...,w_{n+m-r},\hbox{Re} \ z_{r+1},..., 
\hbox{Re} \ z_m), \ \ \ \ \ k=r+1,...,m.$$
We have a commutative diagram
$$\diagram
U \rto^{id} \dto_Z & U \rto^i 
\dto^{Z'} & U\times V^{m-r} \dto^{\tilde{Z}} \\
\C^{n} & \C^n \times \R^{m-r} 
\lto^{\pi_1} & \C^{n+m-r} \lto^{\pi_2}
\enddiagram $$
Now, the vector fields $i_*L_1,...,i_
*L_r, \frac{1}{2} 
\frac{\partial }{\partial x_{r+1}}+i_*L_{r+1},...,
\frac{1}{2}\frac{\partial }{\partial x_m}+
i_* L_m$ on $U \times V^{m-r}$ span a formally integrable
complex subbundle $\tilde{\cal V}$ 
of $\C T \tilde{U}$ such that the 
push-forward $\tilde{Z}_* 
\tilde{{\cal V}}$ becomes 
the bundle $T^{0,1} \tilde{M}$
of antiholomorphic tangent 
vector fields to $\tilde{M}$ which does not depend on the 
choice of $(y_1,...,y_{m-r})$. These vector fields
are the restriction to 
$\tilde{M}$ of the vector fields on 
$\C^n \times \C^{m-r}$ given by
                    \begin{equation}
\tilde{L}_j = \frac{\partial }{\partial \bar{z}_j}+
\sum_{k=r+1}^n b_{r+j \ k} (z,\bar{z},\hbox{Im} \ w) 
\frac{\partial }{\partial \bar{z}_k}, \ \ \ \ \ j=1,...,r.
         \end{equation}
$$\tilde{L}_{r+j}=\frac{\partial }{\partial \bar{w}_{n+j}}+ 
\sum_{k=r+1}^n b_{j \ k} (z,\bar{z},\hbox{Im} \ w) 
\frac{\partial }{\partial \bar{z}_k}, \ \ \ \ \ j=1,...,m-r.$$

\medskip
\noindent
{\sc 6.2. Proposition.} {\it Let 
$M$ be a smooth hypoanalytic manifold. 
Then, for every $p_0\in M$, there
 exists $(Z,y,U)$ a regular
hypoanalytic chart at $p_0$, there 
exists a neighborhood ${\cal X}$ of 0 in $\R^{m-r}$, 
$r=n-\hbox{dim} \ T' \cap T^*M[p_0]$, such that
the mapping $U \times {\cal X}\to 
\C^{n+m-r}$ given by $(p,x) 
\mapsto (Z(p), x+iy(p))$
realizes $U\times {\cal X}$ as 
a generic submanifold of $\C^{n+m-r}$. Moreover, a function
$f$ is a RC function in some 
regularity class ${\cal D}$ on $U$
if and only if there exists
a CR function 
$\tilde{F} \in {\cal D}_{CR}(\tilde{M})$ 
which is independent 
of $(\hbox{Re} \ w_{n+1},...,\hbox{Re} \ w_{n+m-r})$ and
$\tilde{F} \circ \tilde{Z} \circ i=f$.}

$$\diagram
U \times {\cal X} \rto^{\tilde{Z}} 
& \C^n \times \C^{m-r} \dto^{\pi_2} \\
U \dto_Z \uto_i \rto_{Z'} & \C^n\times \R^{m-r} \dto^{F'} \\
\C^n \rto_F & \C
\enddiagram $$

{\it Proof.} Prove the statement 
for the weakest class of functions. Namely, 
a distribution $f$ is a RC distribution
on $U$ if and only if $f$ is annihilated
by ${\cal V}'=Z_*'{\cal V}$. 
Since the function 
$\tilde{F}=f\circ (\tilde{Z} 
\circ i)^{-1}=F' \circ \pi_2$ is by definition
independent of the extra 
variables $(x_{n+1},...,x_{n+m-r})$, 
$\tilde{F}$ is a solution
of ${\cal V}$, so is CR on $\tilde{M}$.
The proof of Proposition 6.2 is complete. $\square$

\medskip
\noindent
{\sc 6.3. Definition.} By a 
{\it local CR lifting} of a hypoanalytic structure, we mean
datas like 6.2.

\medskip
Let $p_0\in M$. According to 
Sussmann, the $\hbox{\small Re} 
\ {\cal V}$-orbit
of $p_0$ can be equipped with a 
structure of a ${\cal C}^{k-1, 
\alpha}$ manifold, possibly
with a finer topology than that 
induced by $M$. If 
${\cal O}_{\hbox{\small Re} 
\ {\cal V}}(M,p_0)$
contains a neighborhood of 
$p_0$ in $M$, the propagation
of analyticity method yields 
that $\tilde{F}=f\circ 
(\tilde{Z} \circ i)^{-1}$ extends
holomorphically into a wedge 
of edge $M$ at $0$ in $\C^{n+m-r}$. 

Indeed, the crucial fact which 
reduces this last statement to the 
generalized Tumanov's theorem 
(\cite{JO3}, \cite{ME1}, 
\cite{TR2}, \cite{TU3}), is that
changes of local CR liftings 
at $p_0$ behave correctly and always are
possible, even if $r(p)$ 
changes, since all the extensions 
successively obtained will 
{\it never} depend on the extra variables.
Therefore, extension of 
RC functions
in the CR lifting sense, 
support propagation,
propagation of the wave 
front set, are straightforward
consequences of those results
which hold true for CR 
structures. This idea goes 
back to Treves \cite{TREV1}.
In conclusion, we can 
restate theorems 1 and 2
for hypoanalytic structures.

\medskip
\noindent
{\sc 6.4. Theorem.} {\it Theorems 1 and 2
hold for ${\cal C}^{2,\alpha}$ 
hypoanalytic manifolds and RC functions.} $\square$

\end{document}